\documentclass[reqno]{amsart}
\usepackage{amsmath,amssymb}
\usepackage[usenames,dvipsnames]{color}
\usepackage{hyperref}
\hypersetup{
    colorlinks=true,       			
    linkcolor=blue, 
    citecolor=blue,        			
    filecolor=blue,      				
    urlcolor=blue           			
}
\newcommand{\corb}[1]{\color{blue}{#1}}

\usepackage[all]{xy}
\xyoption{arc}
\newtheorem*{thm}{Theorem}
\newtheorem*{prop}{Proposition}
\newtheorem*{lem}{Lemma}
\newtheorem*{cor}{Corollary}
\newtheorem{conj}{Conjecture}

\newenvironment{pf}{\paragraph{{\sc Proof}}}{\qed\par\medskip}
\theoremstyle{definition}
\newtheorem*{defn}{Definition}
\newtheorem*{rem}{Remark}

\numberwithin{equation}{section}

\numberwithin{figure}{section}

\newcommand{\hplane}[2]{\mathbb{H}_{#1}^{#2}}
\newcommand{\spec}{\sigma}

\newcommand{\veps}{\varepsilon}
\newcommand{\lp}{\left(}
\newcommand{\rp}{\right)}

\newcommand{\lb}{\left[}
\newcommand{\rb}{\right]}

\newcommand{\g}{\mathfrak{g}}
\newcommand{\h}{\mathfrak{h}}
\newcommand{\gl}{\mathfrak{gl}}

\newcommand{\Sym}{\mathfrak{S}}

\newcommand{\vac}{\Omega}

\newcommand{\QL}{$\mathcal{QL}$}
\newcommand{\EQ}{$\mathcal{EQ}$}
\newcommand{\bfA}{\mathbf{A}}

\newcommand{\cC}{\mathcal{C}}

\newcommand{\E}{\mathcal{E}}

\newcommand{\calO}{\mathcal{O}}
\newcommand{\PP}{\mathcal{P}}

\newcommand{\V}{\mathcal{V}}
\newcommand{\W}{\mathcal{W}}
\newcommand{\X}{\mathcal X}
\newcommand{\calZ}{\mathcal Z}
\newcommand{\calP}{\mathcal P}
\newcommand{\C}{\mathbb{C}}
\newcommand{\nC}{\mathbb{C}^{\times}}
\newcommand{\N}{\mathbb{Z}_{\geq 0}}

\newcommand{\Z}{\mathbb{Z}}
\newcommand{\nZ}{\mathbb{Z}_{\neq 0}}

\newcommand{\ttheta}[3]{\frac{\theta(#1)}{\theta(#2)\theta(#3)}}
\newcommand{\thetaratio}[2]{\frac{\theta\lp #1\rp}{\theta\lp #2\rp}}
\newcommand{\thetaratiop}[2]{\frac{\theta^+\lp #1\rp}{\theta^+\lp #2\rp}}
\newcommand{\thetaration}[2]{\frac{\theta^-\lp #1\rp}{\theta^-\lp #2\rp}}
\newcommand{\weight}{\mathsf{h}}
\newcommand{\eL}{\mathcal{L}_{\hbar,\tau}(\g)}
\newcommand{\Eht}{\mathcal{E}_{\hbar,\tau}(\g)}
\newcommand{\ReL}{\Rep(\mathcal{E}_{\hbar,\tau}(\g))}
\newcommand{\RefdL}{\Rep_{\scriptscriptstyle{\operatorname{fd}}}(\mathcal{E}_{\hbar,\tau}(\g))}
\newcommand{\eX}{\mathfrak{X}}
\newcommand{\eP}{\Phi}
\newcommand{\prps}[1]{S_{\eV}^{(#1)}}

\newcommand{\Rloop}
{\mathcal{O}_{\scriptscriptstyle{\operatorname{int}}}(\qloop)}
\newcommand{\Rloopslice}
{\mathcal{O}^{\slice}_{\scriptscriptstyle{\operatorname{int}}}(\qloop)}

\newcommand{\Rloopnc}
{\mathcal{O}_{\scriptscriptstyle{\operatorname{int}}}^{\scriptscriptstyle{\operatorname{NC}}}
(\qloop)}
\newcommand{\dwtU}{\mathcal{P}_+^U}
\newcommand{\Rloopintro}
{\operatorname{Rep}_{\scriptscriptstyle{\operatorname{fd}}}(\qloop)}
\newcommand{\RloopSintro}
{\operatorname{Rep}^{\slice}_{\scriptscriptstyle{\operatorname{fd}}}(\qloop)}

\newcommand{\partialX}{{\mathsf X}}
\newcommand{\partialH}{{\mathsf H}}
\newcommand{\eV}{\mathbb{V}}

\newcommand{\eW}{\mathbb{W}}
\newcommand{\Thp}[1]{(#1; p)_{\infty}}
\newcommand{\qG}[1]{\theta^+(#1)}
\newcommand{\qGm}[1]{\theta^-(#1)}
\newcommand{\qGpm}[1]{\theta^{\pm}(#1)}
\newcommand{\sfTh}{\mathsf{\Theta}}
\newcommand{\sfhT}{\mathsf{\Xi}}
\newcommand{\Mer}{\mathcal{M}}
\newcommand{\lflat}{$\lambda$--flat }
\newcommand{\ulA}{\underline{\mathcal{A}}}
\newcommand{\ulB}{\underline{\mathcal{B}}}
\newcommand{\bfb}{{\mathbf b}}
\newcommand{\bfc}{{\mathbf c}}

\newcommand{\slice}{\operatorname{S}}
\newcommand{\eLskeleton}{\mathcal{L}^{\slice}_{\hbar,\tau}(\g)}
\newcommand{\eLClinear}{\mathsf{L}^{\slice}_{\hbar,\tau}(\g)}

\newcommand{\fv}{{\mathbf v}}

\newcommand {\Exp}[1]{\exp{(2\pi\iota #1)}}
\newcommand {\wh}[1]{\widehat{#1}}
\newcommand {\ol}[1]{\overline{#1}}
\newcommand {\ul}[1]{\underline{#1}}

\newcommand {\Aut}{\operatorname{Aut}}

\newcommand{\End}{\operatorname{End}}
\newcommand{\Hom}{\operatorname{Hom}}

\newcommand{\Ad}{\operatorname{Ad}}

\newcommand {\aand}{\qquad\text{and}\qquad}
\renewcommand{\dim}{\operatorname{dim}}

\newcommand {\Rep}{\operatorname{Rep}}

\newcommand {\fd}{finite--dimensional }
\newcommand {\lhs}{left--hand side }
\newcommand {\rhs}{right--hand side }
\newcommand {\wrt}{with respect to }
\newcommand {\hw}{highest weight }
\newcommand {\irr}{irreducible }
\newcommand {\Irr}{\operatorname{Irr}}

\newcommand{\ds}{\displaystyle}
\newcommand{\wt}[1]{\widetilde{#1}}

\newcommand {\Comment}[1]{{\footnote{\corb{{#1}}}}}
\newcommand {\Omit}[1]{}

\newcommand{\Th}{\Theta}

\newcommand {\qloop}{U_q(L\g)}
\newcommand {\qaffine}{U_q{\wh{\g}}}

\newcommand {\Yhg}{Y_\hbar(\g)}

\newcommand {\longisom}{\stackrel{\sim}{\longrightarrow}}

\newcommand{\qbin}[3]{\left[\begin{array}{c} #1 \\ #2\end{array}\right]_{#3}}

\renewcommand {\gl}{\mathfrak{gl}}
\renewcommand {\sl}{\mathfrak{sl}}

\newcommand {\sfA}{{\mathsf A}}

\newcommand {\eg}{{\it e.g., }}

\newcommand {\Id}{\operatorname{Id}}
\newcommand {\ie}{{\it i.e.}, }

\newcommand {\KZB}{Knizhnik--Zamolodchikov--Barnard }

\newcommand {\JKOS }{Jimbo--Konno--Odake--Shiraishi }

\newcommand{\bfI}{{\mathbf I}}

\newcommand{\hres}{\ol{\h^*}}
\newcommand{\lres}{\ol{P^{\scriptscriptstyle{\vee}}}}

\newcommand{\ee}[1]{\left(e^{2\pi\iota #1}\right)}

\renewcommand {\Comment}[1]{}

\begin{document}

\title[Elliptic Quantum Groups]
{Elliptic quantum groups and their finite--dimensional representations}
\author[S. Gautam]{Sachin Gautam}
\address{Department of Mathematics, The Ohio State University, 231 W 18th Avenue, Columbus OH 43210}
\email{gautam.42@osu.edu}
\author[V. Toledano Laredo]{Valerio Toledano Laredo}
\address{Department of Mathematics, Northeastern University, 360 Huntington Avenue, Boston, MA 02115}
\email{V.ToledanoLaredo@neu.edu}
\thanks{The first author was supported in part through the Simons Collaboration grant 526947. 
The second author was supported in part through the NSF grant DMS--1505305.}
\maketitle

\begin{abstract}
Let $\g$ be a complex semisimple Lie algebra, $\tau$ a point in
the upper half--plane, and $\hbar\in\C$ a deformation parameter
such that the image of $\hbar$ in the elliptic curve $\C/(\Z+\tau\Z)$
is of infinite order. 
In this paper, we give an intrinsic definition of the category of \fd
representations of the elliptic quantum group $\Eht$ associated
to $\g$. The definition is given in terms of Drinfeld half--currents,
and extends that given by Enriquez--Felder for $\g=\sl_2$ \cite
{felder-enriquez}. When $\g=\sl_n$, it reproduces Felder's RLL
definition \cite {felder-icm} via the Gau{\ss} decomposition obtained
in \cite{felder-enriquez} for $n=2$ and \cite{sachin-rtt} for $n\geq 3$.
We classify the irreducible representations of $\Eht$ in terms of
elliptic Drinfeld polynomials, in close analogy to the case of the
Yangian $\Yhg$ and quantum loop algebra $\qloop$ of $\g$. 
A crucial ingredient in the classification is a functor $\Theta$
from \fd representations of $\qloop$ to those of $\Eht$, which
is an elliptic analogue of the monodromy functor constructed
in \cite{sachin-valerio-2}, and circumvents the fact that
$\Eht$ does not admit Verma modules.
Our classification is new even for $\g=\sl_2$. It holds more
generally when $\g$ is a symmetrisable Kac--Moody algebra,
provided  finite--dimensionality is replaced by an integrability
and category $\calO$ conditions. 
\end{abstract}

\Omit{
arXiv title
Elliptic quantum groups and their finite-dimensional representations
arXiv abstract
Let g be a complex semisimple Lie algebra, tau a point in
the upper half-plane, and h a complex deformation parameter
such that the image of h in the elliptic curve E_tau is of infinite
order. In this paper, we give an intrinsic definition of the category
of finite-dimensional representations of the elliptic quantum group
E_{h,tau}(g) associated to g. The definition is given in terms of
Drinfeld half-currents and extends that given by Enriquez-Felder
for g=sl_2. When g=sl_n, it reproduces Felder's RLL definition
via the Gauss decomposition obtained by Enriquez-Felder for
n=2 and by the second author for n greater than 2. 
We classify the irreducible representations of E_{h,tau} in terms
of elliptic Drinfeld polynomials, in close analogy to the case of the
Yangian Y_h(g) and quantum loop algebra U_q(Lg) of g. 

An important ingredient in the classification, which circumvents
the fact that E_{h,tau} does not appear to admit Verma modules,
is a functor from finite-dimensional representations of U_q(Lg)
to those of E_{h,tau} which is an elliptic analogue of the monodromy
functor constructed in our previous work arXiv:1310.7318.

Our classification is new even for g=sl_2, and holds more
generally when g is a symmetrisable Kac-Moody algebra,
provided  finite-dimensionality is replaced by an integrability
and category O condition. 
}

\setcounter{tocdepth}{1}
\tableofcontents

\newpage

\section{Introduction}\label{sec: intro}

\subsection{}

Let $\g$ be a \fd simple Lie algebra over $\C$. The aim of this paper is to define and
study the category of \fd representations of the elliptic quantum group associated to
$\g$. The results of this paper hold more generally for any symmetrisable Kac--Moody
algebra, provided finite--dimensionality is replaced by an integrability and category
$\calO$ condition, and the main body of the paper deals with this generality. To simplify
the exposition, we restrict our attention to $\g$ of finite type in the introduction.

\subsection{}

Since Baxter's solution of the eight--vertex model \cite{baxter-8v,baxter}, and
Felder's introduction of the dynamical Yang--Baxter equations \cite{felder-icm},
there have been numerous proposals to define elliptic quantum groups (see
Sections \ref{isec: history}--\ref{ss:grt})of this Introduction for a more detailed history of the subject).
One of the contributions of this paper is such a definition which is intrinsic,
uniform for all Lie types, and valid for numerical values of the deformation
and elliptic parameters.

We fix throughout a complex number $\tau$ in the upper half plane $\hplane
{}{}$, and a deformation parameter $\hbar\in\C$ such that the image of $\hbar$
in the elliptic curve $E_\tau=\C/(\Z+\tau\Z)$ is of infinite order. Let $\h\subset\g$
be a fixed Cartan subalgebra, and $\bfI$ the set of vertices of the Dynkin
diagram of $\g$ relative to a choice of simple roots in $\h^*$. In Section
\ref{sec: ellqg}, we define the category $\RefdL$ of \fd representations of
the elliptic quantum group associated to $\g$.
An object of this category is a \fd semisimple $\h$--module $\eV$, together with
a collection of meromorphic, $\End(\eV)$--valued functions $\{\Phi_i(u,\lambda),
\eX_i^{\pm}(u,\lambda)\}_{i\in\bfI}$ of a spectral variable $u\in\C$ and dynamical
parameter $\lambda\in\h^*$. This data is subject to certain periodicity conditions and prescribed commutation
relations which are elliptic analogues of the relations satisfied by the half--currents
in Drinfeld's new realization of Yangians and quantum loop algebras
(see, \eg \cite[\S 2]{sachin-valerio-2} or Section \ref{sec: qla} below). 
Our definition
extends that given by Enriquez--Felder for $\g=\sl_2$ \cite{felder-enriquez} and,
for $\g=\sl_n$, reproduces Felder's RLL category \cite{felder-icm} via the Gau{\ss}
decomposition obtained in \cite{felder-enriquez} for $n=2$, and in \cite{sachin-rtt}
for $n\geq 3$.\footnote{A closely related decomposition, valid when both $\hbar$
and $p=e^{2\pi\iota\tau}$ are formal parameters, was recently obtained by Konno \cite{konno2}.}

\Omit{ obtained via the $RTT$--formalism
of Faddeev--Reshetikhin--Takhtajan \cite{frt-quantization}. More precisely, for $\g=\sl_n$,
our definition agrees with the one given by Felder. This is proved in \cite{sachin-rtt} by 
working out the Gauss decomposition of the $L$--matrix, in close analogy with
the one for the quantum loop algebra by Ding--Frenkel
\cite{ding-frenkel}. Recently this computation in the formal setting (when $\tau$ and $\hbar$
are formal parameters) was carried out by Konno \cite{konno2}.\\}

\subsection{}\label{th:intro class}

From Section \ref{sec: eqg-prep} onwards, 
we restrict our attention to the full subcategory $\eL\subset \ReL$ consisting
of those objects for which the endomorphisms $\Phi_i(u,\lambda)$ are independent
of the dynamical parameter $\lambda$. The main result of this paper (Theorem
\ref{thm: class-final}) is a classification of the irreducible objects in $\eL$ in terms
of elliptic analogues of Drinfeld polynomials. Let for this purpose 
\[SE_\tau=\bigcup_{n\geq 0}E_\tau^n/\Sym_n\]
be the union of the symmetric
powers of the elliptic curve $E_\tau$. A point in $SE_\tau$ may be thought of
as encoding the zeros of a monic elliptic polynomial.

\begin{thm}\label{ithm: main}
The set of isomorphism classes of irreducible objects in $\eL$ is in bijection with
$|\bfI|$ copies of $S E_\tau$ 
\begin{equation}\label{eq:label}
\Irr(\eL) \longisom 
\underbrace{S E_\tau\times\cdots\times SE_\tau}_{|\bfI|}
\end{equation}
\end{thm}

\Omit{Here, and for the entirety of this paper, by $\eV$ being {\em irreducible}
we mean that the only subspaces of $\eV$ stable under all the operators $\{\Phi_i(u),
\eX_i^{\pm}(u,\lambda)\}_{i\in\bfI}$ are $\{0\}$ and $\eV$. The reader might
rightfully wonder if this ad hoc notion of irreducibility agrees with the more
sophisticated one for any category with an initial object. This is indeed so, 
and we shall prove it in Corollary \ref{cor:generalinjective}.}

\subsection{}

The classification above is formally analogous to that of \fd irreducible
representations of Yangians and quantum loop algebra by Drinfeld
polynomials \cite{chari-pressley-qaffine,drinfeld-yangian-qaffine}.
Specifically, given an $\bfI$--tuple 
$\bfb=\left(\{b_i^{(j)}\}_{j=1}^{N_i}\right)_{i\in\bfI}\in
\prod_{i\in\bfI}\C^{N_i}/\Sym_{N_i}$,
there is an irreducible object $\eV(\bfb)$ in the category $\eL$
which contains a non--zero vector $\vac$ such that
\begin{itemize}
\item $\eX_i^+(u,\lambda)\,\vac = 0$ for any $i\in\bfI$, $u\in\C$ and $\lambda\in\h^*$.
\item For any $i\in\bfI$,
\begin{equation}\label{ieq: d-form}
\Phi_i(u)\,\vac = \prod_{j=1}^{N_i} \frac{\theta(u-b_i^{(j)} + d_i\hbar)}{\theta(u-b_i^{(j)})}\,\vac
\end{equation}
where $\theta(x)$ is the odd theta function (see Section \ref{ssec: ell-funs}), and $\{d_i\}$
are the symmetrising integers for the Cartan matrix of $\g$.
\item $\eV$ is spanned by vectors obtained by successively applying the lowering
operators $\{\eX_i^-(u,\lambda)\}_{i\in\bfI}$ to $\vac$.
\end{itemize}
Moreover, shifting each $b_i^{(j)}$ by an element of lattice $\Lambda_\tau=\Z+\tau\Z$
yields an isomorphic representation, and this condition is necessary and sufficient for two
irreducible objects $\eV(\bfb)$ and $\eV(\bfc)$ to be isomorphic.

\subsection{}
\Omit{
The purpose of this paper is to give a uniform and intrinsic definition of $\ReL$ for arbitrary $\g$.
As we said earlier our main theorem classifies irreducible objects of a full subcategory $\eL$ of 
$\ReL$. We conjecture that every object of $\ReL$ is isomorphic to some object of $\eL$.\\
}

The proof of Theorem \ref{ithm: main} can be broken down into the following steps
\begin{itemize}
\item[1.] A triangularity result (Theorem \ref{thm: eqg-pbw}), according to which
an irreducible object $\eV\in\eL$ contains a (unique up to scalar) non--zero weight
vector $\vac$ which is annihilated by the raising operators $\eX^+_i(u,\lambda)$, 
is an eigenvector for each $\Phi_i(u)$, and is such that $\eV$ is spanned by
successively applying lowering operators $\eX_i^-(u,\lambda)$ to $\vac$.
\item[2.] Theorem \ref{thm: class1} states that the eigenvalues of $\Phi_i(u)$ on $\vac$
are of the form \eqref{ieq: d-form}.
\end{itemize}
Steps 1.--2. give rise to a well--defined classifying map \eqref{eq:label}. For Yangians
and quantum loop algebras, the injectivity of the analogous map is proved using Verma
modules. For elliptic quantum groups, Verma modules do not seem to exist, however.
We circumvent this obstacle as follows.
\begin{itemize}
\item[3.] Let $\qloop$ be the quantum loop algebra of $\g$, where $q = e
^{\pi\iota\hbar}$. In Section \ref{sec: functor}, we extend the construction
of \cite{sachin-valerio-2} to the elliptic setting, and obtain
a monodromy functor
\[\sfTh:\Rloopintro \longrightarrow\RefdL\]
whose image lies in $\eL$. We prove further that $\sfTh$ is essentially
surjective, and maps irreducibles to irreducibles (Theorem \ref{thm:inverse}).
These results allow us to rely on the classification of \fd irreducible representations
of $\qloop$ \cite{chari-pressley-qaffine,drinfeld-yangian-qaffine}, and prove
that the map \eqref{eq:label} is in fact a bijection (Theorem \ref{thm:  class-final}).
\end{itemize}

\subsection{}\label{isec: eqg-ops}

We now sketch the construction of the functor $\sfTh$. 
Let $\Rloopintro$ be the category of \fd representations of the quantum loop
algebra $\qloop$. By using half--currents, an object of $\Rloopintro$ may be
viewed as a \fd vector space $\V$, together with a collection of rational, $\End
(\V)$--valued functions $\{\Psi_i(z), \X_i^{\pm}(z)\}_{i\in\bfI}$ of a complex variable
$z$, which are regular at $z=0,\infty$ and satisfy certain commutation relations
(see \cite{beck-kac}, or Section \ref{ssec: rloop}).

Let $\V\in\Rloopintro$. The action of the elliptic quantum group on $\sfTh(\V)=
\V$ is obtained as follows.\footnote{Strictly speaking, the functor $\sfTh$ is only
defined on a dense subcategory of $\Rloopintro$, but we will gloss over this point
in the introduction.}
\begin{itemize}
\item Set $p=e^{2\pi\iota\tau}$, so that $|p|<1$. Let $K_i = \Psi_i(\infty) = \Psi_i(0)
^{-1}\in GL(\V)$, and consider the infinite products
\[G_i^{+}(z) = \prod_{n\geq 1} \left(K_i \Psi_i(p^{n}z)\right)
\quad\text{and}\quad
G_i^{-}(z) = \prod_{n\geq 1} \left(K_i^{-1} \Psi_i(p^{-n}z)\right)
 \]
The $GL(\V)$--valued function $G_i^+(z)$ and $G_i^-(z)$ are meromorphic on $\C^\times$,
holomorphic in a neighbourhood of $z=0,\infty$ respectively, and such that $G_i^+(0)
=1=G_i^-(\infty)$. By construction, $G_i^{\pm}(z)$ are the canonical fundamental
solutions near $z=0,\infty$ of the $p$--difference equations
\Omit{\[
G_i^+(pz) = \left[K_i\Psi_i(pz)\right]^{-1}\cdot G_i^+(z)
\quad\text{and}\quad
G_i^-(pz) = \left[K_i^{-1}\Psi_i(z)\right]\cdot G_i^-(z)
\]}
\begin{align*}
G_i^+(pz)	&= \left[K_i\Psi_i(pz)\right]^{-1}\cdot G_i^+(z)\\
G_i^-(pz)	&= \left[K_i^{-1}\Psi_i(z)\right]\cdot G_i^-(z)
\end{align*}
determined by the commuting field $\Psi_i(z)$ of $\qloop$.
\item By definition, the action of the current $\Phi_i(u)$ of $\Eht$
on $\V$ is given by
\begin{equation}\label{ieq: Phi}
\Phi_i(u) = \left.G_i^+(z)\Psi_i(z)G_i^-(z)\right|_{z=e^{2\pi\iota u}}
\end{equation}
By construction, $\Phi_i(u)$ is a doubly quasi--periodic function of $u$,
that is, it satifies the following
\[\Phi_i(u+1)=\Phi_i(u)\aand
\Phi_i(u+\tau)=K_i^{-2}\Phi_i(u)\]
and is the monodromy of the $p$--difference equation defined by $\Psi_i(z)$.

\item The raising and lowering operators act on  $\V$ by the
following contour integral
\[
\eX_i^{\pm}(u,\lambda) = \frac{\theta^+(0)}{\theta^+(d_i\hbar)}
\int_{C_i^{\pm}} \ttheta{u-v+\lambda_i}{u-v}{\lambda_i}
G_i^{\pm}\lp e^{2\pi\iota v}\rp \X_i^{\pm}\lp e^{2\pi\iota v}\rp
\, dv
\]
where $\theta^+(x) = \prod_{n\geq 1} (1-p^n e^{2\pi\iota x})$ and 
$\lambda_i = (\lambda,\alpha_i)$. The choice of contours $C_i^{\pm}$ is explained
in Section \ref{ssec: contours}.
\end{itemize}

\begin{thm}\label{ithm: functor}
The operators $\{\Phi_i(u), \eX^{\pm}_i(u,\lambda)\}$ defined above
satisfy the defining relations of $\Eht$ and therefore give rise to an exact, 
faithful functor
\[\sfTh: \Rloopintro\longrightarrow\eL\subset\RefdL\]
which commutes with the forgetful functor to vector spaces.
\end{thm}


\subsection{}\label{issec: hw-calculation}

We now state a number of additional properties of the functor $\sfTh$,
which mirror those of the functor constructed in \cite{sachin-valerio-2}
and are used in the proof of Theorem \ref{ithm: main}.

The action of $\sfTh$ on highest weights is easily determined. Let $\V
\in\Rloopintro$ be irreducible, with Drinfeld polynomials $\{P_i(w)\}_{i
\in\bfI}$, written as
\[
P_i(w) = \prod_{k=1}^{N_i} (w-\beta_k^{(i)})
\]
where $\beta_k^{(i)}\in\nC$ for each $i\in\bfI$ and $1\leq k\leq N_i$. Thus, 
if $\vac\in\V$ is the (unique up to scaling) highest weight vector, we have
\[
\Psi_i(z)\,\vac = \prod_{k=1}^{N_i} \frac{q_iz-q_i^{-1}\beta_k^{(i)}}{z-\beta_k^{(i)}}\,\vac
\]
for every $i\in \bfI$, where $q_i = q^{d_i}$. 
By \eqref{ieq: Phi}, and Jacobi's triple product identity 
\cite[Chapter 21]{whittaker-watson},
\[\theta(u) = 
\frac{\theta^+(u)\sin(\pi u)\theta^-(u)}{\theta^+(0)^2}\]
where $\theta^{\pm}(u) = \prod_{n\geq 1} (1-p^n e^{\pm 2\pi\iota u})$, we get
\[
G_i^{\pm}(z)\,\vac = \prod_k \frac{\theta^{\pm}(u-b_k^{(i)}+d_i\hbar)}
{\theta^{\pm}(u-b_k^{(i)})}\,\vac \, 
\aand
\Phi_i(u)\,\vac =  \prod_k \frac{\theta(u-b_k^{(i)}+d_i\hbar)}{\theta(u-b_k^{(i)})}\,\vac
\]
where
$b_k^{(i)}\in\C$ are such that $e^{2\pi\iota b_k^{(i)}}=\beta_k^{(i)}$. Thus, $\sfTh$
maps the roots of Drinfeld polynomials for $\qloop$ to their image in $E_\tau$.

\subsection{}

Choose a subset $\slice\subset \C$ which is a fundamental domain for
the action of $\Lambda_{\tau}=\Z+\tau\Z$ on $\C$, and is stable under shifts
by $\hbar/2$. Define the full subcategories
\[\RloopSintro\subset\Rloopintro
\aand
\eLskeleton\subset\eL\]
as follows. $\RloopSintro$ 
consists of those representations
$\V$ for which the poles of $\{\Psi_i(z)^{\pm 1}\}$ lie in $e^{2\pi\iota \slice}$.
$\eLskeleton$ consists of representations $\eV$ for which the generalized
eigenvalues of $\{\Phi_i(u)\}$ are of the following form
\Omit{\begin{equation}\label{ieq: d-form-knight}
A_i(u) = \prod_k \frac{\theta(u-c_{i,k}+d_i\hbar)}{\theta(u-c_{i,k})}
\prod_l \frac{\theta(u-c_{i,l}'-d_i\hbar)}{\theta(u-c_{i,l}')}
\end{equation}}
\begin{equation}\label{ieq: d-form-knight}
A_i(u) =
\prod_k \frac{\theta(u-c_{i,k}+d_i\hbar)}{\theta(u-c_{i,k})}
\prod_l \frac{\theta(u-c_{i,l}')}{\theta(u-c_{i,l}'+d_i\hbar)}
\end{equation}
where $c_{i,k},c_{i,l}'\in\slice$.

\begin{thm}\label{ithm: functor-slice}\hfill
\begin{enumerate}
\item The functor $\sfTh$ restricts to a functor $\sfTh^{\slice}: \RloopSintro \to \eLskeleton$.
This restriction sets up a bijection between simple objects of the two categories.
\item Every object of $\eL$ is isomorphic to an object of $\eLskeleton$.
\end{enumerate}
\end{thm}
\noindent
For a given $\eV\in\eLskeleton$, we explicitly construct a representation $\V\in
\RloopSintro$ such that $\eV = \sfTh(\V)$ (Theorem \ref{thm:inverse}). This
inverse contruction is achieved
by solving a Riemann--Hilbert factorization problem for doubly--quasi periodic
(abelian, matrix--valued) functions in Section \ref{sec: factorization}, 
and is analogous to \cite[\S 6]{sachin-valerio-2}. 

\Omit{Theorem \ref{ithm: main} is obtained by combining Theorem \ref{ithm: functor-slice}
with the calculation given in Section \ref{issec: hw-calculation} and 
the classification theorem of \fd irreducible
representations of $\qloop$.}

\Omit{
The second part of the theorem involves the {\em second gauge transformation} given in Proposition
\ref{prop:lambda-independence} which is achieve by carefully analyzing the solution of the factorization
problem. In order to elaborate on that point, one observes that for any generalized eigenvalue
of the form \eqref{ieq: d-form-knight}, we can always assume that $c_{i,k},c_{i,l}'$ lie
in $\slice$. However $\text{Contant}_i$ need not equal $1$. As it turns
out this constant can be absorbed in an isomorphism of $\eL$ which results in a representation
lying in $\eLskeleton$.\\

This phenomenon can even be seen for the trivial object: $\eV = \C$ one--dimensional
vector space with $\eX_i^{\pm}(u,\lambda) = 0$ and $\Phi_i(u) = C_i\in\nC$. This object, denoted by
${\mathbf 1}_{\ul{C}}$, clearly
does not lie in $\eLskeleton$, unless all $C_i = 1$. However we can define
\[
\varphi(\lambda) = \exp\lp \frac{1}{\hbar} \sum_{j\in\bfI} (\lambda,\varpi_j^{\vee})\ln(C_j)\rp
\]
which is an isomorphism between ${\mathbf 1}_{\ul{C}}$ and ${\mathbf 1}_{\ul{1}}\in\eLskeleton$.\\
}

\subsection{}

One difference between Theorem \ref{ithm: functor-slice} and the analogous relation
between \fd representations of Yangians and quantum loop algebras obtained in \cite
{sachin-valerio-2} is that the functor $\sfTh^{\slice}: \RloopSintro \to \eLskeleton$ 
is not an equivalence. This is
because the two categories are defined over different fields. 
Indeed, $\Rloopintro$ is defined over $\C$, whereas the field of definition of $\eLskeleton$ is
\[
\End_{\eLskeleton}({\mathbf 1}) = \{ \varphi : \h^* \to \C \text{ meromorphic such that }
\varphi(\lambda+\hbar\alpha_i) = \varphi(\lambda)\}
\]
However, if $\mathsf{L}_{\hbar,\tau}^{\slice}(\g)$ is the category with the same objects
as $\eLskeleton$, and morphisms the $\C$--linear homomorphisms commuting
with the operators $\{\Phi_i(u), \eX_i^{\pm}(u,\lambda)\}$, the following is an
immediate corollary of the proof of Theorem \ref{ithm: functor-slice}.

\begin{cor}
The functor $\sfTh^{\slice} : \RloopSintro \to \mathsf{L}_{\hbar,\tau}^{\slice}(\g)$ is
an isomorphism of categories.
\end{cor}

\subsection{Elliptic $R$--matrices}\label{isec: history} 

To put the results of this paper in perspective, we now give some background
on elliptic quantum groups. Their origins can be traced to  Baxter's elliptic
$R$--matrix, which records the Boltzmann weights of the eight--vertex model
\cite{baxter-8v, baxter}, and its generalisation $R_B(u)$ to higher rank by
Belavin \cite{belavin}. The study of the quantum groups corresponding to
$R_B(u)$ was initiated by Sklyanin \cite{sklyanin,sklyanin2} in rank 1, and
extended by Cherednik \cite{cherednik-sklyanin} to higher rank. The \fd
representations of these {\em Sklyanin algebras} \cite{feigin-odesskii1, feigin-odesskii2},
which may be thought of as elliptic deformations of the loop algebra of $\sl_n$,
and are also known as vertex--type elliptic algebras, are also investigated
in these works.

\subsection{Dynamical Yang--Baxter equations} 
 
It is natural to look for elliptic solutions of the Yang--Baxter equations (YBE)
corresponding to any simple Lie algebra $\g$. Belavin and Drinfeld proved
that, for the classical YBE, such solutions only exist for $\g=\sl_n$ \cite
{belavin-drinfeld-book}. Partly to remedy this, Felder
\cite{felder-icm} introduced the {\it dynamical} quantum Yang--Baxter equation
(DQYBE)
\begin{multline*}\label{eq:dqybe}
R_{12}(u,\lambda-\hbar\weight^{(3)})R_{13}(u+v,\lambda)R_{23}(v,\lambda-\hbar\weight^{(1)})\\
=R_{23}(v,\lambda)R_{13}(u+v,\lambda-\hbar\weight^{(2)})R_{12}(u,\lambda)
\end{multline*}
where $\h$ is an abelian Lie algebra, $V$ a \fd semisimple $\h$--module, and
$R:\C\times \h^*\to\End(V\otimes V)$ a meromorphic function. The equation
holds in $\End(V\otimes V\otimes V)$, the subscript $R_{ij}$ indicates
which tensor factors $R$ acts on, and the symbol $\weight^{(j)}$ denotes
the weight of the $j^{\text{th}}$ tensor factor. For example, if $v_3\in V_3$
is a vector of weight $\mu_3\in\h^*$, then
\[
R_{12}(u,\lambda-\hbar\weight^{(3)})\,
v_1\otimes v_2\otimes v_3 = \left(R(u,\lambda-\hbar\mu_3)\,
v_1\otimes v_2 \right)\otimes v_3
\]

Felder obtained an elliptic solution $R_F(u,\lambda)$ of the above equations
when $\h$ is the Cartan subalgebra of diagonal matrices in $\g=\gl_n$, and $V=\C^n$
is the defining representation of $\g$. He also pointed out that the \KZB equations
satisfied by genus 1 conformal blocks \cite{felder-wieczerkowski} give rise to
solutions of the underlying classical dynamical Yang--Baxter equation for any
semisimple $\g$, thus bypassing in a sense Belavin and Drinfeld's no go theorem.
\Omit{
Some authors prefer the term ``face--type elliptic algebra" for the object
resulting from Felder's $R$--matrix \cite{abf}.
}

Starting from Felder's discovery, a large number of results have been obtained
for elliptic quantum groups. These can be roughly grouped into four related
threads described in \ref{ss:Felder RLL}--\ref{ss:dyn q gps}.

\subsection{Felder's elliptic representation category}\label{ss:Felder RLL}

The RTT formalism of Faddeev--Reshetikhin--Takhtajan \cite
{frt-quantization}, applied to Felder's $R$--matrix $R_F(u,\lambda)$
yields a tensor category, which can be thought of as the category of
representations of a face--type elliptic algebra associated to $\g=\sl_n$.
This study of this category for $\g=\sl_2$ was
initiated by Felder-Varchenko in \cite{felder-icm,fv-reps}. An extension
of this approach to other Lie types seems, however, less convenient
for other classical and for exceptional Lie algebras.


\subsection{Quasi--Hopf twists}

Following the suggestion of Babelon--Bernard--Billey to use 
dynamical twists to produce solutions of the DQYBE \cite{bbb},
Fr\o{}nsdal discovered that both the Baxter--Belavin and Felder
$R$--matrices 
can be obtained from the trigonometric $R$--matrix of the quantum
affine algebra $\wh{\sl_n}$ by using such a twist \cite{fronsdal1,fronsdal2}.
Fr\o{}nsdal's analysis was completed and extended by \JKOS \cite{jkos} and
Etingof--Schiffmann \cite{etingof-schiffmann-link,
etingofschiffmann-elliptic}, who showed that the $R$--matrix
of the quantum affine algebra $\wh{\g}$ of any simple Lie algebra
$\g$ can be canonically twisted to an elliptic solution of the DQYBE,
thus producing the first examples of such solutions outside
of type $\sfA$. In particular, this suggested defining the elliptic
quantum group $\Eht$ corresponding to $\g$ as the quantum affine
algebra of $\g$, with a twisted (dynamical) coproduct and $R$--matrix
\cite{etingofmoura-KL, etingof-schiffmann-link, etingofschiffmann-elliptic}. 
Theorems \ref{ithm: functor} and \ref{ithm: functor-slice} suggest
that such an extrinsic definition of 
$\Eht$ may only be appropriate when the parameter 
$p=e^{2\pi\iota\tau}$ is formal since, for numerical $\tau$, the
category of \fd representations of $\qloop$ is a covering of
that of those of $\Eht$.

\Omit{
In particular, in \cite{etingofschiffmann-elliptic} Felder's $R$--matrix
is obtained from a version of the fusion and exchange construction
for $U_q(L\sl_n)$ originally introduced in \cite{frenkel-reshetikhin,
etingof-frenkel-kirillov}, and it is pointed out that the latter leads to
elliptic solutions of the quantum dynamical Yang--Baxter equations
for any simple $\g$.
Thus elliptic quantum groups can be (extrinsically)
defined to be same as quantum loop algebras, but with coproduct
and $R$--matrix appropriately twisted. 
In a related vein, for $\g=
\sl_n$, the vector representation $\C^n$ can be used to construct
a functor from the category of \fd representations of $U_q(L\sl_n)$
to that of the elliptic quantum group of $\sl_n$ see \cite{etingofmoura-KL,
etingofschiffmann-elliptic}.
It is perhaps worth pointing out that the work of \JKOS \cite{jkos} was also motivated by the algebraic analysis of Jimbo--Miwa
\cite{jimbo-miwa}. For this purpose, one needs to incorporate the central extension,
which for the ``vertex--type elliptic algebra" of $\sl_n$ was carried out by
Foda et al in \cite{fijkmy}. In order to complete the ingredients required for
this analysis, one also needs the notion of Drinfeld currents, which for
$\sl_2$ was obtained almost simultaneously by Konno \cite{konno-first}
and Enriquez--Felder \cite{felder-enriquez} and was used extensively in 
\cite{jkos2}.
}

\subsection{The elliptic quantum group $U_{q,p}(\wh{\g})$}

A related elliptic quantum group $U_{q,p}(\wh{\g})$ was introduced by
Konno \cite{konno-first} for $\g=\sl_2$, and any simple $\g$ by \JKOS
\cite{jkos} as an extension of the quantum 
affine algebra $\qaffine$ by a one--dimensional Heisenberg algebra.
A presentation of $U_{q,p}(\wh{\g})$ in terms of elliptic Drinfeld (full)
currents was obtained by Kojima--Konno for $\g=\sl_n$ \cite{kojima-konno},
and more recently by Farghly--Konno--Oshima for an arbitrary $\g$ \cite{konno-farghly-oshima}.

The presentation in \cite{konno-farghly-oshima} is closely related
to that given in terms of half--currents in Section \ref {sec: ellqg}.
Indeed, if one regards the ratios of theta functions in Section \ref{sec: ellqg} 
as formal power series in the variable $p = e^{2\pi\iota\tau}$ whose coefficients are
Laurent polynomials in $z=e^{2\pi\iota u}$, the commutation relations
given in Section \ref{ssec: categoryL} can be formally expanded to
get those given in \cite[Defn. 2.1]{konno-farghly-oshima}. The converse
implication does not appear to be so straightforward, however.

lt is worth pointing out that the definition and presentation
of $U_{q,p}(\wh{\g})$ given in \cite{konno-first, konno2} 
are only valid when $\hbar=\log(q)$ and $p$ are formal,
whereas the classification of \fd irreducible representations
does not hold in the formal setting since the ground ring is
not algebraically closed. This is also the case for the quantum loop
algebra $\qloop$ when $q = e^{\pi\iota\hbar}$, with $\hbar$ formal,
where the classification via Drinfeld polynomials fails to hold.

\Omit{
While there are a few
differences between our work and Konno's, the most significant in our opinion
is that Konno's elliptic quantum group $U_{q,p}(\widehat{\g})$ is only defined
when the quantum and elliptic parameters $q,p$ are formal variables. For a
reader familiar with Konno's works, we would like to point out further relations
between his construction and the one we are proposing:

\begin{itemize}
\item The operators $\{\Phi_i(u,\lambda),\eX_i^{\pm}(u,\lambda)\}_{i\in\bfI}$ we
consider in this paper are to be thought of as ``half currents".
\item We merely wish to remark that getting to the commutation relations of
Section \ref{ssec: categoryL} from \cite[Definition 2.1]{konno-farghly-oshima}
is a non--trivial task. We derived them in an indirect fashion by carrying out
the computations of \cite{sachin-valerio-2} in elliptic setting.
\item Another relevant point is that classification of \fd irreducible representation
does not hold in the formal setting since $\C[[\hbar,\tau]]$ is not algebraically
closed. This is the case even for the quantum loop algebra $\qloop$ when
$q = e^{\pi\iota\hbar}$, with $\hbar$ formal, where the classification via Drinfeld
polynomials fails to hold.
\end{itemize}
}

\Comment{Things I have not mentioned (enough/yet): 1) Enriquez--Felder
give a half--current presentation for $\g=\sl_2$, 2) we do not address the
important question of what an elliptic quantum group is (say this in the
paragraph on algebroids) 3) central extensions.}

\subsection{Dynamical quantum groups}\label{ss:dyn q gps}

It is natural to ask about the algebraic structure of the dynamical quantum
groups defined by the DQBYE. These are {\em Hopf algebroids}, as defined 
by Etingof--Varchenko \cite{etingof-varchenko-classical,etingof-varchenko-quantum,
etingof-varchenko-exchange} (see also \cite{etingof-latour} for an exposition).
For $\g=\sl_2$, Konno has shown that $U_{q,p}(\widehat{\sl_2})$ has a 
tructure of Hopf algebroid \cite{konno-first}.
\Omit{
We note however the forthcoming work of Yang--Zhao
\cite{yang-zhao-elliptic} which characterises elliptic quantum groups as
algebra objects in a symmetric monoidal category of coherent sheaves
on symmetric powers of an elliptic curve.}
\Omit{
they are {\em Hopf algebroids}, and their study was taken up by Etingof and Varchenko
in \cite{etingof-varchenko-classical,etingof-varchenko-quantum, etingof-varchenko-exchange} (see also \cite
{etingof-latour} for an exposition). This theory gives a conceptual framework
which accounts for the dynamical
parameter, but is yet to incorporate the spectral variable. In other words, the
algebraic objects which should correspond to $R_F(u,\lambda)$ are still not
defined.
}

\subsection{Geometric representation theory}\label{ss:grt} 

There has been a resurgence of interest in elliptic quantum groups recently,
especially in connection to the geometry of Nakajima quiver varieties. Recall
that, for $\g$ simply--laced, the Yangian (resp. quantum loop algebra) of $\g$
acts on the (equivariant) cohomology (resp. $K$--theory) of the Nakajima quiver
varieties corresponding to the Dynkin diagram of $\g$ \cite{nakajima-qaffine,
varagnolo-yangian}. Maulik--Okounkov \cite{maulik-okounkov-qgqc} obtained
a new construction of Yangians via stable envelopes in the equivariant
cohomology of these quiver varieties.
This construction has been extended to the elliptic setting by Aganagic--Okounkov 
\cite{aganagic-okounkov}, and leads to a geometric definition of their
category of representations. 
We also note that a sheafified version of the elliptic quantum group
corresponding to a symmetric Kac--Moody algebra $\g$ has been
recently obtained by Yang--Zhao \cite{yang-zhao-elliptic}, and greatly
clarifies its algebraic structure. Yang--Zhao also proved this elliptic
quantum group acts on the equivariant elliptic cohomology of the
corresponding Nakajima quiver varieties.

\Omit{On a related note, elliptic quantum groups
can also be realized via the cohomological Hall algebra construction of Yang--Zhao
\cite{yang-zhao-elliptic}, which sheds light on the type of algebraic objects
they are.}

\subsection{Outline of the paper}

In Section \ref{sec: ellqg}, we define the category $\ReL$ and its full
subcategory $\eL$ for any symmetrisable Kac--Moody algebra $\g$. 
Section \ref{sec: eqg-prep} is devoted to proving a few basic properties
of $\eL$. Its main result is an analogue of the triangular decomposition
(Theorem \ref{thm: eqg-pbw}).
In Section \ref{sec: irr-class1}, we prove that the highest weight of an
\irr object of $\eL$ is of the form \eqref{ieq: d-form}, and also establish
a version of Knight's lemma. In Section \ref{sec: qla}, we review the
definition of the quantum loop algebra $\qloop$, and the classification
of the simple objects in $\Rloop$ in terms of Drinfeld polynomials. 
Section \ref{sec: functor} gives the construction of the functor
$\sfTh: \Rloop \to \eL$. In Section \ref{sec: factorization}, we give
a solution to  the factorization problem for (abelian) doubly quasi--periodic
functions. We use this factorization in Section \ref{sec: functor-surjective}
to construct a right inverse to the functor $\sfTh$. The final Section
\ref{sec: class2} gives the classification of irreducible objects in 
$\eL$. Appendix \ref{sec: serre-relns} contains the description of
the Serre relations for elliptic quantum groups. We also give a proof
of these relations in type $\mathsf{A}$. 

\subsection{Acknowledgments}

We are grateful to a number of people and institutions where parts of this
project were carried out. We would in particular like to thank Damien Calaque
and Giovanni Felder
for their invitation to ETH's Institute for Mathematical Research (FIM) in June
2012, and for numerous conversations on elliptic quantum groups, Olivier
Schiffmann and the Institut Henri Poincar\'e for their invitation to spend May
2014 at IHP as part of their research in Paris program, the Perimeter Institute
for inviting the second author in November 2015, Amritanshu Prasad and
the IMSc Chennai for inviting the first author in December 2015, and finally
David Ben--Zvi, Roman Bezrukavnikov, Alexander Braverman and
the Simons Center for Geometry and Physics for their invitation in January
2016.

\Omit{
Felder (FIM invitation in June, 2012)
Olivier Schiffmann - Research in Paris (May, 2014). Perimeter Institute (November 2015 - Valerio's visit). 
IMSc Chennai (December 2015: Amritanshu Prasad for inviting Sachin). 
Simons Center (January 2016). 
}

\section{Elliptic quantum groups}\label{sec: ellqg}

\subsection{}
Let $\bfA = (a_{ij})_{i,j\in\bfI}$ be a symmetrisable generalized Cartan
matrix \cite{kac}. Thus, $a_{ii}=2$ for any $i\in\bfI$, $a_{ij}\in\Z_{\leq
0}$ for any $i\neq j\in \bfI$, and there exists a diagonal matrix $D$
with positive integer entries $\{d_i\}_{i\in\bfI}$ such that $D\bfA$ is
symmetric. We assume that $(d_i)$ are relatively prime.

Let $(\h,\{\alpha_i\}_{i\in\bfI},\{\alpha_i^{\vee}\}_{i\in\bfI})$ be the unique
realization of $\bfA$. Thus, $\h$ is a complex vector space of dimension
$2|\bfI|-\text{rank}(\bfA)$, $\{\alpha_i\}_{i\in\bfI}\subset \h^*$ and $\{\alpha
_i^{\vee}\}_{i\in\bfI} \subset \h$ are linearly independent sets and, for any
$i,j\in\bfI$, $\alpha_j(\alpha_i^{\vee}) = a_{ij}$. Let $(.,.):\h^*\times\h^*
\to \C$ denote a non--degenerate symmetric bilinear form satisfying
\[
(\lambda,\alpha_i) = d_i\lambda(\alpha_i^{\vee}),\ \forall \lambda\in\h^*, i\in\bfI.
\]

Below, we will assume that our {\em dynamical variables} lie in
the subspace
$\ds \hres := \bigoplus_{i\in\bfI} \C\alpha_i \subset \h^*$. Our
periodicity conditions will involve the following ``lattice":

\[
\lres := \{\gamma\in\hres | (\gamma,\alpha_i)\in\Z \text{ for every } i\in\bfI\}.
\]

\begin{rem}\label{rm:hres}
Note that, if $\bfA$ is degenerate, then $\lres\subset\hres$ contains a linear
subspace, namely the radical of the bilinear form on $\hres$, defined by $D\bfA$.
The following fact is immediate from the definitions, and will be used in
Section \ref{sec: eqg-prep}: $\hres/(\lres+\tau\lres)$ is compact,
for any $\tau\in\C$ such that $\operatorname{Im}(\tau)\neq 0$.\\

In this paper, we will also need to shift the dynamical parameter by elements of 
the weight ``lattice".
\[
P := \{\gamma\in\h^* : (\gamma,\alpha_i)\in d_i\Z,\ \forall\ i\in\bfI \}
\]
This is achieved, as usual, by quotienting $\h^*$ by its subspace dual to
the radical of $D\bfA$ restricted to $\hres$.
\end{rem}

Let $\g$ be the Kac--Moody algebra associated to $\bfA$.

\subsection{}\label{ssec: ell-funs}

We fix a complex number $\tau\in\mathbb{H}$ in the upper half plane.
Let $\Lambda_{\tau} = \Z +\tau\Z\subset\C$ and
set $p= \exp(2\pi\iota\tau)$ so that $|p|<1$.
Let $\theta(u)$ be the odd theta function, namely the
unique holomorphic function of a complex variable $u$, such that:
\begin{itemize}
\item $\theta(u+1) = -\theta(u)$
\item $\theta(u+\tau) = -e^{-\pi\iota\tau}e^{-2\pi\iota u} \theta(u)$
\item $\theta(u) = 0 \iff u\in \Lambda_{\tau}$
\item $\theta'(0) = 1$
\end{itemize}

For this theta function, we have $\theta(-u)=-\theta(u)$.
Set $z = e^{2\pi\iota u}$. Explicitly, we have the following expression of $\theta(u)$
\begin{equation}\label{eq: Jacobi-triple}
\theta(u) = -\frac{e^{-\pi\iota u}}{2\pi\iota} \Thp{z}\Thp{pz^{-1}}\Thp{p}^{-2}
\end{equation}
where $\Thp{z} = \prod_{k\geq 0} (1-zp^k)$ is the standard $q$--Pochhamer notation.
Define:
\begin{equation}\label{eq:qGpm}
\qGpm{u} := \prod_{n\geq 1}\left(1-p^n e^{\pm 2\pi\iota u}\right).
\end{equation}
Note that $\qG{u}= \qGm{-u} = \Thp{pz}$. Moreover,
we have the following identity:
\begin{equation}\label{eq:thetaspm}
\theta(u) = \sin(\pi u)\frac{\qG{u}\qGm{u}}{\qG{0}^2}
\end{equation}

In subsequent calculations the following identity will be extensively used
(called Fay's trisecant identity)
\begin{multline}\label{eq: fti}
\tag{FTI}\theta(a-c)\theta(a+c)\theta(b-d)\theta(b+d) = \theta(a-b)\theta(a+b)\theta(c-d)\theta(c+d)
\\ + \theta(a-d)\theta(a+d)\theta(b-c)\theta(b+c)
\end{multline}

\subsection{The category $\ReL$}\label{ssec: categoryL}

Let $\hbar\in\C$ be such that $\Z\hbar\cap \Lambda_{\tau} = \{0\}$.

An object of $\ReL$ is an $\h$--diagonalizable module $\eV$ with \fd
weight spaces $\ds \eV = \bigoplus_{\mu\in\h^*} \eV_{\mu}$, together
with meromorphic $\End(\eV)$--valued functions $\{\eP_i(u,\lambda),
\eX^{\pm}_i(u,\lambda)\}_{i\in\bfI}$ of a {\em spectral variable} $u\in\C$
and a {\em dynamical variable} $\lambda\in\hres$,
satisfying the following axioms.\\

\noindent {\bf Category $\mathcal{O}$ and integrability condition.}
\begin{itemize}
\item There exist $\mu_1,\ldots,\mu_r\in\h^*$ such that $\eV_{\mu}\neq 0$ implies
that $\mu < \mu_k$ for some $k=1,\ldots,r$.
\item For each $\mu\in\h^*$ such that $\eV_{\mu}\neq 0$ and $i\in\bfI$, there
exists $N>0$ such that $\eV_{\mu-n\alpha_i}=0$ for all $n\geq N$.
\end{itemize}

\noindent {\bf Periodicity conditions.}
\begin{itemize}
\item $\eP_i(u+1,\lambda) = \eP_i(u,\lambda)$ and 
$\eP_i(u+\tau,\lambda) = e^{-2\pi\iota\hbar d_i \alpha_i^{\vee}}\eP_i(u,\lambda)$.\\

\item $\eX^{\pm}_i(u+1,\lambda) = \eX_i(u,\lambda)$ and
$\eX^{\pm}_i(u+\tau,\lambda) = e^{-2\pi\iota(\lambda,\alpha_i)}\eX_i^{\pm}(u,\lambda)$.\\

\item $\eX^{\pm}_i(u,\lambda+\gamma) = \eX^{\pm}_i(u,\lambda)$ 
and $\Phi_i(u,\lambda+\gamma) = \Phi_i(u,\lambda)$
for every $\gamma\in \lres$. 
\item Let $D_i^{\pm}\subset\C\times\hres$ be the set of poles of $\eX_i^{\pm}(u,\lambda)$.
Then we require that $D_i^{\pm}$ is stable under shifts by $\Lambda_{\tau}\times
(\lres +\tau \lres)$.
\end{itemize}


\noindent {\bf Commutation relations.}
\begin{itemize}
\item[(\EQ1)] For each $i,j\in\bfI$ and $h\in\h$ we have $[h,\Phi_i(u)]=0$ and
\[
\Phi_i\lp u,\lambda + \frac{\hbar}{2}\alpha_j\rp \Phi_j\lp v,\lambda-\frac{\hbar}{2}\alpha_i\rp=
\Phi_j\lp v,\lambda + \frac{\hbar}{2}\alpha_i\rp \Phi_i\lp v,\lambda-\frac{\hbar}{2}\alpha_j\rp
\]

Moreover, we assume that $\det(\Phi_i(u,\lambda)) \not\equiv 0$ for each $i\in\bfI$.

\item[(\EQ2)] For each $i\in\bfI$ and $h\in\h$ we have
\[
[h,\eX_i^{\pm}(u,\lambda)] = \pm \alpha_i(h)\eX_i^{\pm}(u,\lambda)
\]

\item[(\EQ3)] For each $i,j\in\bfI$ let $a=d_ia_{ij}\hbar/2$ and let $\lambda_j = (\lambda,\alpha_j)$.
Then the following relation holds on $\eV_{\mu}$.
\begin{multline*}
\Phi_i\lp u, \pm\lp \lambda-\frac{\hbar}{2}\mu\rp + \frac{\hbar}{2}\alpha_i\rp
\eX_j^{\pm}(v,\lambda) 
\Phi_i\lp u, \pm\lp \lambda-\frac{\hbar}{2}\mu\rp + \frac{\hbar}{2}(\alpha_i-\alpha_j)\rp^{-1}
\\= 
\frac{\theta(u-v\pm a)}{\theta(u-v\mp a)} \eX_j^{\pm}(v,\lambda\pm\hbar\alpha_i) 
\pm \frac{\theta(2a)\theta(u-v\mp a-\lambda_j)}{\theta(\lambda_j)\theta(u-v\mp a)}
\eX^{\pm}_j(u\mp a, \lambda\pm\hbar\alpha_i)
\end{multline*}

\item[(\EQ4)] For each $i,j\in\bfI$ and $\lambda\in\hres$, let $a = d_ia_{ij}\hbar/2$. 
Then we have
\begin{multline*}
\theta(\lambda_i+\lambda_j)\theta(u-v\mp a) \eX^{\pm}_i\lp u,\lambda\pm\frac{\hbar}{2}\alpha_j\rp
\eX_j^{\pm}\lp v,\lambda\mp\frac{\hbar}{2}\alpha_i\rp \\
-\theta(\lambda_i\pm a)\theta(u-v-\lambda_j)\eX^{\pm}_i\lp u,\lambda\pm\frac{\hbar}{2}\alpha_j\rp
\eX_j^{\pm}\lp u+\lambda_i,\lambda\mp\frac{\hbar}{2}\alpha_i\rp\\
-\theta(\lambda_j\mp a)\theta(u-v+\lambda_i)\eX_i^{\pm}\lp v+\lambda_j,\lambda\pm\frac{\hbar}{2}\alpha_j\rp
\eX_j^{\pm}\lp v,\lambda\mp\frac{\hbar}{2}\alpha_i\rp\\
=\theta(\lambda_i+\lambda_j)\theta(u-v\pm a) \eX^{\pm}_j\lp v,\lambda\pm\frac{\hbar}{2}\alpha_i\rp
\eX_i^{\pm}\lp u,\lambda\mp\frac{\hbar}{2}\alpha_j\rp \\
-\theta(\lambda_i\mp a)\theta(u-v-\lambda_j)\eX^{\pm}_j\lp u+\lambda_i,\lambda\pm\frac{\hbar}{2}\alpha_i\rp
\eX_i^{\pm}\lp u,\lambda\mp\frac{\hbar}{2}\alpha_j\rp\\
-\theta(\lambda_j\pm a)\theta(u-v+\lambda_i)\eX_j^{\pm}\lp v,\lambda\pm\frac{\hbar}{2}\alpha_i\rp
\eX_i^{\pm}\lp v+\lambda_j,\lambda\mp\frac{\hbar}{2}\alpha_j\rp
\end{multline*}
where $\lambda_i = (\lambda,\alpha_i)$ and $\lambda_j = (\lambda,\alpha_j)$.\\

\item[(\EQ5)] For each $i,j\in\bfI$ and $\lambda_1,\lambda_2\in\hres$ such that
$\lambda_1+\lambda_2 = \hbar(\mu+\alpha_i-\alpha_j)$, the following relation
holds on the weight space $\eV_{\mu}$, with $\lambda = \lambda_1 - \hbar(\mu-\alpha_i+\alpha_j)/2$.
\[
\theta(d_i\hbar)[\eX_i^+(u,\lambda_1), \eX^-_j(v,\lambda_2)] = \delta_{ij}\lp
\frac{\theta(u-v+\lambda_{1,i})}{\theta(u-v)\theta(\lambda_{1,i})}\Phi_i(v,\lambda)
+\frac{\theta(u-v-\lambda_{2,i})}{\theta(u-v)\theta(\lambda_{2,i})}\Phi_i(u,\lambda)\rp
\]
where $\lambda_{s,i} = (\lambda_s,\alpha_i)$ for $s=1,2$.
\end{itemize}

\subsection{Morphisms}\label{ssec:defn-morphisms}
A morphism between two objects $\eV$ and $\eW$
of $\ReL$ is a meromorphic function $\varphi : \hres\to\Hom_{\h}(\eV,\eW)$,
such that
\begin{itemize}
\item For each $i\in\bfI$ we have
\[
\varphi\lp\lambda+\frac{\hbar}{2}\alpha_i\rp\Phi_i(u,\lambda) = \Phi_i(u,\lambda)
\varphi\lp\lambda-\frac{\hbar}{2}\alpha_i\rp\\
\]

\item For each $i\in\bfI$ the following equation holds in 
$\Hom_{\C}(\eV_{\mu},\eW_{\mu\pm\alpha_i})$

\[
\varphi\lp \pm\lambda \mp \frac{\hbar}{2}(\mu\pm\alpha_i) + \frac{\hbar}{2}\alpha_i\rp
\eX^{\pm}_i(u,\lambda) = 
\eX^{\pm}_i(u,\lambda)
\varphi\lp \pm\lambda \mp \frac{\hbar}{2}\mu - \frac{\hbar}{2}\alpha_i\rp
\]

\end{itemize}

The composition of morphisms is defined as $(\varphi\circ\psi)(\lambda) = \varphi(\lambda)\psi(\lambda)$.
The constant function $\mathbf{1}_{\eV}(\lambda) = \operatorname{Id}_{\eV}$ is
the identity morphism in $\Hom_{\ReL}(\eV,\eV)$.

\begin{rem}
The definition of a morphism given above differs from the one usually given
for RLL type algebras, for instance in \cite{etingof-latour,etingof-varchenko-exchange,
fv-reps}. The reason is the freedom in choosing this notion:
given a set--theoretic map $\beta : \hres\to\hres$, and a zero--weight meromorphic function 
$\varphi(\lambda)\in\Aut(\eV)$ one can easily verify that
the following shifted conjugation respects the relations (\EQ 1)--(\EQ 5):

\begin{align*}
\wt{\Phi}_i(u,\lambda) &= \varphi\lp\lambda+\frac{\hbar}{2}\alpha_i\rp\Phi_i(u,\lambda)
\varphi\lp\lambda-\frac{\hbar}{2}\alpha_i\rp^{-1} \\
\wt{\eX_i}^{\pm}(u,\lambda) &= 
\varphi\lp \pm\lambda+\beta(\mu\pm\alpha_i)\mp\frac{\hbar}{2}(\mu\pm\alpha_i)+\frac{\hbar}{2}\alpha_i\rp
\eX_i^{\pm}(u,\lambda)\\
&\qquad .\varphi\lp \pm\lambda+\beta(\mu)\mp\frac{\hbar}{2}\mu-\frac{\hbar}{2}\alpha_i\rp^{-1}
\end{align*}

Clearly any two choices are related by an affine linear change of dynamical variables.
The one given above corresponds to $\beta\equiv 0$. From the RLL picture
one arrives at $\beta(\mu) = \pm\frac{\hbar}{2}\mu$,
depending on the choice of the Gauss decomposition of the $L$--matrix.
\end{rem}

\subsection{}\label{ssec: eqg-rks}

For $\g=\sl_2$ the category defined above is same as the one
defined by Felder \cite{felder-icm} and further studied by Felder and Varchenko
\cite{fv-reps}. The reader can consult \cite{felder-enriquez} for a proof of 
this assertion. More generally, the previous statement holds for $\g=\sl_n$.
A proof of this will be given in \cite{sachin-rtt}. Similar calculations were
carried out in \cite{konno2}, where the formal case is treated (\ie $p=e^{2\pi\iota\tau}$
is considered a formal variable).

We also note that we have not imposed any Serre-type relations in the
definition above. This is in accordance with the assertion that for Yangians
and quantum loop algebras, the Serre relations are a consequence of the
other relations and the category $\mathcal{O}$ and integrability axioms
\cite[\S 2]{sachin-valerio-2}. For the interested reader, we include the
elliptic analogue of the  Serre relations in Appendix \ref{sec: serre-relns},
together with a proof that they follow from the other relations in type
$\mathsf{A}$.

\subsection{}\label{ssec: categoryrestricted}

As mentioned in the introduction, we restrict our attention to the full subcategory
$\eL$ of $\ReL$ consisting of objects $(\eV,\Phi_i(u,\lambda),\eX_i(u,\lambda))$
such that $\Phi_i(u,\lambda)$ is independent of $\lambda$, for each $i\in \bfI$. The
relations (\EQ1) and (\EQ3) of the previous section become simpler for such objects
as follows.

\begin{itemize}

\item[(\EQ1)] For each $i,j\in\bfI$ and $h\in\h$ we have
\[
[\Phi_i(u),\Phi_j(v)]=0 \aand [h,\Phi_i(u)]=0
\]

\item[(\EQ3)] For each $i,j\in\bfI$ let $a=d_ia_{ij}\hbar/2$ and let $\lambda_j = (\lambda,\alpha_j)$.
Then
\begin{multline*}
\Phi_i(u)\eX_j^{\pm}(v,\lambda)\Phi_i(u)^{-1} = 
\frac{\theta(u-v\pm a)}{\theta(u-v\mp a)} \eX_j^{\pm}(v,\lambda\pm\hbar\alpha_i) \\
\pm \frac{\theta(2a)\theta(u-v\mp a-\lambda_j)}{\theta(\lambda_j)\theta(u-v\mp a)}
\eX^{\pm}_j(u\mp a, \lambda\pm\hbar\alpha_i)
\end{multline*}

\end{itemize}

\subsection{}\label{ssec: eq34-variant}

From the relation (\EQ 3) we obtain, by setting $v=u\pm a$
\[
\Ad(\Phi_i(u))\eX_j^{\pm}(u\pm a,\lambda) = \frac{\theta(\lambda_j\pm 2a)}{\theta(\lambda_j)}
\eX_j^{\pm}(u\mp a, \lambda\pm\hbar\alpha_i)
\]
Using this we obtain the following from (\EQ 3) by taking $\Ad(\Phi_i(u))^{-1}$ on the both
sides, and replacing $\lambda$ by $\lambda\mp\hbar\alpha_i$:
\begin{itemize}
\item[(\EQ $3^\prime$)] For every $i,j\in\bfI$ the following holds
\begin{multline*}
\Phi_i(u)^{-1}\eX_j^{\pm}(v,\lambda)\Phi_i(u) = 
\frac{\theta(u-v\mp a)}{\theta(u-v\pm a)} \eX_j^{\pm}(v,\lambda\mp\hbar\alpha_i) \\
\mp \frac{\theta(2a)\theta(u-v\pm a-\lambda_j)}{\theta(\lambda_j)\theta(u-v\pm a)}
\eX^{\pm}_j(u\pm a, \lambda\mp\hbar\alpha_i)
\end{multline*}
\end{itemize}

Similar to the case of Yangians and quantum loop algebras 
(see \cite[Remarks 2.3, 2.10]{sachin-valerio-2}), the relation (\EQ 4) can be 
simplified when $i=j$ as follows. Let $u=v$ in (\EQ 4) for $i=j$ and hence
$a=d_i\hbar$. We get
\begin{multline}\label{eq:EQ4-special}
\theta(\lambda_i\pm d_i\hbar)\eX_i^{\pm}\lp u,\lambda\pm\frac{\hbar}{2}\alpha_i\rp
\eX_i^{\pm}\lp u+\lambda_i,\lambda\mp\frac{\hbar}{2}\alpha_i\rp\\
-\theta(\lambda_i\mp d_i\hbar)\eX_i^{\pm}\lp u+\lambda_i,\lambda\pm\frac{\hbar}{2}\alpha_i\rp
\eX_i^{\pm}\lp u,\lambda\mp\frac{\hbar}{2}\alpha_i\rp\\
= 
\pm\frac{\theta(d_i\hbar)\theta(2\lambda_i)}{\theta(\lambda_i)}
\eX_i^{\pm}\lp u,\lambda\pm\frac{\hbar}{2}\alpha_i\rp
\eX_i^{\pm}\lp u,\lambda\mp\frac{\hbar}{2}\alpha_i\rp
\end{multline}

This implies the following relation again from (\EQ 4)

\begin{itemize}
\item[(\EQ $4^\prime$)] For each $i\in\bfI$ and $\lambda\in\hres$ we have
\begin{multline*}
\theta(u-v\mp d_i\hbar)\eX_i^{\pm}\lp u,\lambda\pm\frac{\hbar}{2}\alpha_i\rp
\eX_i^{\pm}\lp v,\lambda\mp\frac{\hbar}{2}\alpha_i\rp\\
\pm \frac{\theta(u-v-\lambda_i)\theta(d_i\hbar)}{\theta(\lambda_i)}
\eX_i^{\pm}\lp u,\lambda\pm\frac{\hbar}{2}\alpha_i\rp
\eX_i^{\pm}\lp u,\lambda\mp\frac{\hbar}{2}\alpha_i\rp
=\\
\theta(u-v\pm d_i\hbar)\eX_i^{\pm}\lp v,\lambda\pm\frac{\hbar}{2}\alpha_i\rp
\eX_i^{\pm}\lp u,\lambda\mp\frac{\hbar}{2}\alpha_i\rp\\
\pm \frac{\theta(u-v+\lambda_i)\theta(d_i\hbar)}{\theta(\lambda_i)}
\eX_i^{\pm}\lp v,\lambda\pm\frac{\hbar}{2}\alpha_i\rp
\eX_i^{\pm}\lp v,\lambda\mp\frac{\hbar}{2}\alpha_i\rp
\end{multline*}
\end{itemize}

In this form the relation (\EQ 4) appeared for $\g=\sl_2$
in \cite[Prop. 1.1]{felder-enriquez}.

\section{Some preparatory results}\label{sec: eqg-prep}

The aim of this section is to prove 
a weak triangularity property
for irreducible objects of 
$\eL$. In order to do so we will introduce the notion of a 
\lflat object and prove that every $\eV\in\eL$ is isomorphic
to some \lflat $\eV^{\flat}$. The triangularity property holds for $\eV^{\flat}$,
as stated in Theorem \ref{thm: eqg-pbw}.\\

\subsection{Poles of $\eX_j^{\pm}(u,\lambda)$}\label{ssec:xpm-poles}

As a consequence of the periodicity axioms of section
\ref{ssec: categoryL}, we have the following

\begin{prop}\label{prop:xpm-poles}
The set of poles of $\eX^{\pm}_j(u,\lambda)$ are contained in the union of
affine hyperplanes of the form $u=b$ ($b\in\C$) or
$\lambda_j=(\lambda,\alpha_j)\in\Lambda_{\tau}$.
Moreover the poles at $\lambda_j \in \Lambda_{\tau}$ are of order at most $1$.
\end{prop}

\begin{pf}
Let us restrict to the $+$ case for definiteness. Let $D_i^+\subset\C\times\hres$
be the set of poles of $\eX_i^+(u,\lambda)$. Consider a point $(v^0,\lambda^0)\in
D_i^+$.\\

\noindent {\bf Claim.} Either $(v^0,\lambda)\in D_i^+$ for every $\lambda\in\hres$
or $D_i^+\cap\{(v^0,\lambda^0+p\hbar\alpha_i)\}_{p\in\Z}$ is finite.\\

\noindent {\em Proof of the claim:} Assume that the set 
$D_i^+\cap\{(v^0,\lambda^0+p\hbar\alpha_i)\}_{p\in\Z}$
is infinite. Since $D_i^+$ is stable under shifts by $\lres + \tau \lres$
in the dynamical variable, and $p\hbar\alpha_i\in \lres+\tau \lres$
implies $p=0$ by the genericity assumption on $\hbar$, the set 
$\lambda^0+\Z\hbar\alpha_i+(\lres+\tau \lres)$ is dense in $\hres$.
Hence $(v^0,\lambda)$ is in $D_i^+$ for every $\lambda\in\hres$.\\

In order to prove the proposition we will restrict our attention to the later
case. Namely, there exist $n,m\in\N$ such that
\[
D_i^+\cap\{(v^0,\lambda^0+p\hbar\alpha_i)\}_{p\in\Z} \subset \{(v^0,\lambda^0-m\hbar\alpha_i),
\ldots,(v^0,\lambda^0+n\hbar\alpha_i)\}
\]

Let $\lambda^s = \lambda^0-m\hbar\alpha_i$ and $\lambda^t = \lambda^0+n\hbar\alpha_i$.
Without loss of generality we may assume $(v^0,\lambda^s),(v^0,\lambda^t)\in D_i^+$.
Consider the relation (\EQ 3) with $j=i$ and $v=v^0$:
\begin{multline*}
\Ad(\Phi_i(u))\eX_i^+(v^0,\lambda) = 
\frac{\theta(u-v^0+d_i\hbar)}{\theta(u-v^0-d_i\hbar)}
\eX_i^+(v^0,\lambda+\hbar\alpha_i) \\
+ \frac{\theta(d_i\hbar)\theta(u-v^0-\lambda_i-d_i\hbar)}
{\theta(u-v^0-d_i\hbar)\theta(\lambda_i)}
\eX_i^+(u-d_i\hbar,\lambda+\hbar\alpha_i)
\end{multline*}

Since the \lhs has a pole at $\lambda=\lambda^t$, so must the \rhs. However $(v^0,
\lambda^t+\hbar\alpha_i)\not\in D_i^+$. Therefore we have $\lambda^t_i\in\Lambda_{\tau}$.
A similar argument using (\EQ 3') from Section \ref{ssec: eq34-variant} implies
that $\lambda^s_i\in\Lambda_{\tau}$. Therefore we get
\[
(n+m)\hbar (2d_i) = (\lambda^s-\lambda^t,\alpha_i)\in\Lambda_{\tau}
\]
which implies that $n=m=0$. Therefore $\lambda^0_i = M-N\tau\in\Lambda_{\tau}$. 

Again consider
the relation (\EQ 3)
\begin{multline*}
\Ad(\Phi_i(u))\lp \theta(\lambda_i)\eX_i^+(v,\lambda)\rp = 
\frac{\theta(u-v+d_i\hbar)}{\theta(u-v-d_i\hbar)}
\theta(\lambda_i)\eX_i^+(v,\lambda+\hbar\alpha_i) \\
+ \frac{\theta(d_i\hbar)\theta(u-v-\lambda_i-d_i\hbar)}
{\theta(u-v-d_i\hbar)}
\eX_i^+(u-d_i\hbar,\lambda+\hbar\alpha_i)
\end{multline*}

The \rhs can be specialized at $\lambda=\lambda^0$ and gives
\[
\left.\text{R.H.S. } \right|_{\lambda=\lambda^0} = 
(-1)^{M+N}e^{-\pi\iota N^2\tau}e^{-2\pi\iota N(u-v-d_i\hbar)}\theta(d_i\hbar)
\eX_i^+(u-d_i\hbar,\lambda^0+\hbar\alpha_i)
\]

Thus we get

\begin{multline*}
\left.\theta(\lambda_i)\eX_i^+(v,\lambda)\right|_{\lambda=\lambda^0}
= e^{2\pi\iota Nv} \Ad(\Phi_i(u))^{-1}\cdot\\ \cdot\lp
(-1)^{M+N}e^{-\pi\iota N^2\tau}e^{-2\pi\iota N(u-d_i\hbar)}\theta(d_i\hbar)
\eX_i^+(u-d_i\hbar,\lambda^0+\hbar\alpha_i)
\rp
\end{multline*}

The \rhs is either identically zero or nowhere zero. The former case
will contradict the fact that $(v^0,\lambda^0)\in D_i^+$ since $\theta(x)$
has only simple zeroes. In the later
case we get that $(v,\lambda^0)\in D_i^+$ for every $v$ and the corresponding
pole is simple.

\end{pf}

\subsection{Partial fractions}\label{ssec: xpm-partial}

Let $S_j^{\pm}\subset\C$ be a choice of representatives modulo $\Lambda_{\tau}$
of poles of $\eX_j^{\pm}(u,\lambda)$ in the spectral variable $u$. Consider
the Laurent series expansion of these functions near $b\in S_j^{\pm}$:
\[
\eX^{\pm}_j(u,\lambda) = \sum_{n\in\N} \partialX_{j;b,n}^{\pm}(\lambda)\frac{1}{(u-b)^{n+1}}
+ \text{ part holomorphic near $b$}
\]

Note that the sum appearing above need not be finite (since $\eV$ is infinite--dimensional
in general). However, for each weight $\mu\in\h^*$, the sum is finite
in $\Hom_{\C}(\eV_{\mu},\eV_{\mu\pm\alpha_j})$.
We have the following analogue of
partial fraction decomposition of a rational function:

\begin{lem}\label{lem: xpm-partial}
The $\End(\eV)$--valued functions $\eX^{\pm}_j(u,\lambda)$ have the following
form, where $\partialX^{\pm}_{j;b,n}(\lambda)$ are holomorphic functions of $\lambda$
\[
\eX^{\pm}_j(u,\lambda) = 
\sum_{\begin{subarray}{c} b\in S_j^{\pm} \\ n\in\N\end{subarray}}
\partialX^{\pm}_{j;b,n}(\lambda)\frac{(-\partial_u)^n}{n!}
\left(\ttheta{u-b+\lambda_j}{u-b}{\lambda_j}\right)
\]
\end{lem}

\begin{pf}
The proof of this lemma is fairly standard. We take the difference of the
two sides of the equation given above

\[
F(u,\lambda) = 
\eX^{\pm}_j(u,\lambda) - 
\sum_{\begin{subarray}{c} b\in S_j^{\pm} \\ n\in\N\end{subarray}}
\partialX^{\pm}_{j;b,n}(\lambda)\frac{(-\partial_u)^n}{n!}
\left(\ttheta{u-b+\lambda_j}{u-b}{\lambda_j}\right)
\]

Then this function is holomorphic in $u$ and has the quasi--periodicity:
$F(u+1,\lambda) = F(u,\lambda)$ and $F(u+\tau,\lambda) = e^{-2\pi\iota\lambda_j}
F(u,\lambda)$. Since there are no such holomorphic functions, other than $0$,
we get the desired equation. Finally the claim that $\partialX^{\pm}_{j;b,n}(\lambda)$
are holomorphic functions of $\lambda$ follows from Proposition 
\ref{prop:xpm-poles} which states that $\theta(\lambda_j)\eX_j^{\pm}(u,\lambda)$
has no poles in $\lambda$ variable.
\end{pf}

\subsection{A difference equation}\label{ssec: xpm-difference}

As a consequence of the commutation relation (\EQ 3) we get the
following equation for the functions $\partialX^{\pm}_{j;b,n}(\lambda)$
appearing in Lemma \ref{lem: xpm-partial} above.

\begin{lem}\label{lem: xpm-difference}
For each $i,j\in\bfI$, $n\in\N$ and $b\in S_j^{\pm}$ we have
\begin{equation}\label{eq: xpm-difference}
\Ad(\Phi_i(u)).\partialX^{\pm}_{j;b,n}(\lambda) = 
\sum_{k\geq 0} \left.\frac{\partial_v^k}{k!}
\left(\frac{\theta(u-v\pm a)}{\theta(u-v\mp a)}\right)
\right|_{v=b} \partialX^{\pm}_{j;b,n+k}(\lambda\pm\hbar\alpha_i)
\end{equation}
where again $a = \hbar d_ia_{ij}/2$.
\end{lem}

\begin{pf}
Consider the commutation relation (\EQ 3):
\begin{multline*}
\Phi_i(u)\eX_j^{\pm}(v,\lambda)\Phi_i(u)^{-1} = 
\frac{\theta(u-v\pm a)}{\theta(u-v\mp a)} \eX_j^{\pm}(v,\lambda\pm\hbar\alpha_i) \\
\pm \frac{\theta(2a)\theta(u-v\mp a-\lambda_j)}{\theta(\lambda_j)\theta(u-v\mp a)}
\eX^{\pm}_j(u\mp a, \lambda\pm\hbar\alpha_i)
\end{multline*}

If we multiply both sides of this equation by $(v-b)^n$ and integrate over
a small circle around $b$, we get the equation \eqref{eq: xpm-difference}.
\end{pf}

\subsection{Generalized eigenspaces}\label{ssec: Phi-eigen}

Given $\mu\in\h^*$ define $\Mer(\mu)$ to be the set of $\bfI$--tuples
of meromorphic functions $\ulA = (A_i(u))$ such that for each $i\in\bfI$
we have
\[
A_i(u+1) = A_i(u) \qquad \text{and}\qquad A_i(u+\tau) = e^{-2\pi\iota\hbar(\mu,\alpha_i)}
A_i(u)
\]

According to the commutation relation (\EQ 1), we can decompose
each weight space $\eV_{\mu}$ into the generalized eigenspaces for the functions
$\{\Phi_i(u)\}$. Define:
\[
\eV_{\mu}[\ulA] := \{v\in\eV_{\mu} : \text{ for every } i\in\bfI, (\Phi_i(u)-A_i(u))^Nv=0 
\text{ for } N\gg 0\}
\]

Then we have
\[
\eV_{\mu} = \bigoplus_{\ulA\in\Mer(\mu)} \eV_{\mu}[\ulA]
\]

\subsection{}\label{ssec: important-prop}

\begin{prop}\label{prop: important-prop}
For a given $\mu\in\h^*$ and $\ulA\in\Mer(\mu)$, $\ulA^{\pm}\in\Mer(\mu\pm\alpha_j)$
consider the following composition, for $j\in\bfI$ and $n\in\N$:
\[
\xy
(0,25)*{\eV_{\mu}[\ulA]}="a";
(50,25)*{\eV_{\mu\pm\alpha_j}[\ulA^{\pm}]}="b";
(25,28)*{\partialX^{\pm}_{j;b,n}(\lambda)_{\ulA}^{\ulA^{\pm}}};
(0,0)*{\eV_{\mu}}="c";
(50,0)*{\eV_{\mu\pm\alpha_j}}="d";
(25,-3)*{\partialX^{\pm}_{j;b,n}(\lambda)};
{\ar@{^(->} (0,20)*{}; (0,5)*{}};
{\ar@{->>} (50,5)*{}; (50,20)*{}};
{\ar (5,0)*{}; (45,0)*{}};
{\ar@{-->} (5,25)*{}; (40,25)*{}};
\endxy
\]
If $\partialX^{\pm}_{j;b,n}(\lambda)_{\ulA}^{\ulA^{\pm}}\neq 0$ then there exist
$\alpha^\pm\in Q = \sum_{i\in\bfI} \Z\alpha_i$ such that for every $i\in\bfI$ we have

\[
A_i^{\pm}(u) = A_i(u)\frac{\theta(u-b\pm a)}{\theta(u-b\mp a)} e^{2\pi\iota\hbar (\alpha_i,\alpha^\pm)}
\]
where $a=d_ia_{ij}\hbar/2$.
Moreover, in this case, for every $m\in\N$ we have
\[
\exp\lp \mp2\pi\iota (\lambda,\alpha^{\pm})\rp
\partialX^{\pm}_{j;b,m}(\lambda)  \text{ is independent of } \lambda
\]
\end{prop}

\begin{pf}
For the purposes of the proof, let us keep $j\in\bfI$ and $b\in\C$ fixed.
Consider $W^{\pm} = \Hom_{\C}(\eV_\mu[\ulA],\eV_{\mu\pm\alpha_j}[\ulA^{\pm}])$ as
a \fd vector space over $\C$. We are given $Y_i(u) = \Ad(\Phi_i(u))\in\End(W^{\pm})$
which are $\End(W^{\pm})$--valued meromorphic function of $u$,
such that $Y_i(u)-\frac{A_i^{\pm}(u)}{A_i(u)}\Id$ is a nilpotent operator.\\

According to the equation \eqref{eq: xpm-difference}, we are looking for
a collection of holomorphic $W^{\pm}$--valued functions, $x^{\pm}_0(\lambda),\cdots,x^\pm_N(\lambda)$
where $\lambda\in\ol{\h^*}$ which are $\lres$--periodic.
These functions are required to satisfy: for each $i\in\bfI$,
\begin{equation}\label{eq:pf-diff1}
Y_i(u).x^\pm_n(\lambda) = \sum_{l\geq 0} E^\pm_{i,l}(u)x^\pm_{n+l}(\lambda\pm\hbar\alpha_i)
\end{equation}
where the functions $E^\pm_{i,l}(u)$
are given by
\[
E^\pm_{i,l}(u) = \left.\frac{\partial_v^l}{l!}
\left(\frac{\theta(u-v\pm a)}{\theta(u-v\mp a)}\right)
\right|_{v=b} \text{ where } a = d_ia_{ij}\hbar/2
\]

We prove the assertion of the proposition by descending (finite) induction,
as follows. Consider the function $x^\pm_N(\lambda)$ and its Laurent
series expansion. The role of Laurent monomials is played by $\exp(2\pi\iota (\lambda,\alpha))$ 
for $\alpha\in Q$. Therefore the Laurent series of $x^{\pm}_N(\lambda)$ has the following form:

\[
x_N^{\pm}(\lambda) = \sum_{\alpha\in Q} C_N^\pm(\alpha) \exp\lp 2\pi\iota (\lambda,\alpha)\rp
\]

Substituting this expansion in the equation \eqref{eq:pf-diff1},
and comparing the coefficients of the Laurent series, we get

\[
Y_i(u)C^\pm_N(\alpha) = E^\pm_{i,0}(u)\exp\lp \pm 2\pi\iota\hbar
(\alpha_i,\alpha)\rp C_N(\alpha)
\]

Thus, if for some $\alpha$ we have $C_N^\pm(\alpha)\neq 0$, then (since the only eigenvalue of $Y_i(u)$
on $W^{\pm}$ is $A_i^{\pm}(u)/A_i(u)$) we obtain:

\[
E^{\pm}_{i,0}(u)\exp\lp \pm 2\pi\iota\hbar(\alpha_i,\alpha)\rp 
= \frac{A_i^{\pm}(u)}{A_i(u)}
\]

This proves the first assertion of the proposition. Namely, there must
exist $\alpha^{\pm}\in Q$ such that

\[
\frac{A_i^{\pm}}{A_i} = \frac{\theta(u-b\pm a)}{\theta(u-b\mp a)}e^{2\pi\iota\hbar
(\alpha_i,\alpha^{\pm})} \text{ where } a = d_ia_{ij}\hbar/2
\]

Note that such $\alpha^{\pm}$ are not unique (unless
$\bfA$ is non--degenerate), but their classes modulo $Q^0 := \{\gamma\in Q : (\gamma,
\alpha_i)=0\}$ are. 
Thus we will have to make a choice of $\alpha^{\pm}\in Q$
so that the equation given above holds.
Then we get that $x_N^\pm(\lambda) = e^{\pm2\pi\iota (\lambda,\alpha^{\pm})}
C_N^{\pm}$.\\

Let us assume that we have proved
\[
\exp\lp \mp2\pi\iota(\lambda,\alpha^{\pm})\rp x^\pm_m(\lambda) = C^\pm_m \text{ is independent of }
\lambda
\]
for each $m$ with $n<m\leq N$. Now we will prove it for $x^\pm_n(\lambda)$. Again, let us
write the Laurent series expansion of $x^\pm_n(\lambda)$ as

\[
x^{\pm}_n(\lambda) = \sum_{\alpha\in Q/Q^0} C^{\pm}_n(\alpha)\exp\lp 2\pi\iota (\lambda,\alpha)\rp
\]

We substitute this in the equation \eqref{eq:pf-diff1} and compare
the coefficients. For each $i\in\bfI$, we get

\begin{multline*}
Y_i(u)C^\pm_n(\alpha) = E^\pm_{i,0}(u) \exp\lp \pm 2\pi\iota \hbar(\alpha_i,\alpha)\rp C^\pm_n(\alpha)
+\\ \delta_{\alpha, \alpha^{\pm} \text{ mod }Q^0}
e^{\pm2\pi\iota\hbar(\alpha_i,\alpha^{\pm})} \sum_{l\geq 1}E^\pm_{i,l}(u)C^\pm_{n+l}
\end{multline*}

The same argument as above applies and we are done.
\end{pf}

\subsection{\lflat representations and first gauge transformation}
\label{ssec: lflat}

Given a weight $\mu\in\h^*$ of $\eV$ consider the set of generalized eigenvalues
\[
\Mer(\mu;\eV) := \{\ulA\in\Mer(\mu) : \eV_{\mu}[\ulA]\neq 0\}
\]

In view of Proposition \ref{prop: important-prop} above, we introduce an
equivalence on the set $\Mer(\mu;\eV)$: $\ulA \sim \ulA'$ if, and only if
there exists $\alpha\in Q$
such that for every $i\in\bfI$ we have $A_i'(u) = A_i(u)\exp(2\pi\iota\hbar(\alpha_i,\alpha))$.

\begin{defn}\label{defn: lflat}
We say $\eV\in\eL$ is \lflat if for every $\ulA,\ulA'\in \Mer(\mu;\eV)$ such
that $\ulA\sim\ulA'$ we have $\ulA = \ulA'$.
\end{defn}

\begin{prop}\label{prop: gauge1}
Every $\eV\in\eL$ is isomorphic to a \lflat $\eV^{\flat}\in\eL$.
\end{prop}

\begin{pf}
We claim that there exists $\psi(\lambda)\in \Aut(\eV)$ such that upon conjugating
with $\psi$ the we obtain a \lflat object of $\eL$. This automorphism $\psi$
is weight preserving: $\psi(\lambda) = \oplus_{\mu} \psi_{\mu}(\lambda)$
where $\psi_\mu(\lambda)\in\Aut(\eV_{\mu})$. Below we give the construction
of $\psi_\mu$.\\

Choose representatives $\{\ulA^{(b)}\}_{b\in B}\subset\Mer(\mu;\eV)$ of their respective
equivalence classes modulo $\sim$. 
Let $\eV_\mu[\ulB]$ be a generalized eigenspace, where $\ulB\sim\ulA^{(b)}$ for some
$b\in B$. Then by definition we have $\alpha\in Q$ such that
\[
B_i(u) = A^{(b)}_i(u) e^{2\pi\iota\hbar (\alpha,\alpha_i)}
\]
The action of $\psi_\mu(\lambda)$ on the generalized weight space $\eV_\mu[\ulB]$
is then given by
\[
\left.\psi_\mu(\lambda)\right|_{\eV_\mu[\ulB]} = 
\exp\lp -2\pi\iota (\lambda,\alpha)\rp \Id_{\eV_\mu[\ulB]}
\]
Having defined $\psi(\lambda)$, we consider the object $\eV^{\flat}$ of $\eL$
obtained by twisting $\eV$ via $\psi$ in accordance with the notion of morphisms
introduced in Section \ref{ssec: categoryL}. Namely, $\eV^{\flat}=\eV$
as an $\h$--diagonalizable module and the operators $\{\Phi_i^{\flat},
\eX_i^{\pm, \flat}\}$ on a weight space $\eV^{\flat}_{\mu} = \eV_{\mu}$
are given by:

\begin{align*}
\Phi_i^{\flat}(u) &= \varphi\lp\lambda+\frac{\hbar}{2}\alpha_i\rp\Phi_i(u)
\varphi\lp\lambda-\frac{\hbar}{2}\alpha_i\rp^{-1}\\
\eX_i^{\pm,\flat}(u,\lambda)&= \varphi\lp \pm\lambda \mp \frac{\hbar}{2}(\mu\pm\alpha_i) + \frac{\hbar}{2}\alpha_i\rp
\eX^{\pm}_i(u,\lambda)\\
&\qquad .\varphi\lp \pm\lambda \mp \frac{\hbar}{2}\mu - \frac{\hbar}{2}\alpha_i\rp^{-1}
\end{align*}

The reader can easily verify that, by the construction of $\psi$, $\eV^{\flat}$ is 
a \lflat object.
\end{pf}

\subsection{}\label{ssec: eqg-pbw}

Assume $\eV\in\eL$ is a \lflat object. Then the following is a consequence
of Proposition \ref{prop: important-prop}.

\begin{lem}\label{lem: eqg-pbw}
Let $\eV\in\eL$ be a \lflat object. Then for every $\mu\in\h^*$ a weight of $\eV$
the following two subspaces of $\eV_{\mu\pm\alpha_j}$ are identical:
\[
\sum_{\begin{subarray}{c} u\in\C \\ \lambda\in\hres\end{subarray}}
\eX_j^{\pm}(u,\lambda)\eV_{\mu} = 
\sum_{u\in\C} \eX_j^{\pm}(u,\lambda^0)\eV_{\mu}
\]
where $\lambda^0\in\hres$ is chosen arbitrarily such that $\lambda^0_j\not\in\Lambda_{\tau}$.
\end{lem}

\begin{pf}
According to the partial fraction decompotion of $\eX^\pm_j(u,\lambda)$ 
given in Lemma \ref{lem: xpm-partial}, we have
\[
\text{Span of } \{\eX^\pm_j(u,\lambda)\}_{u\in\C}
 = \text{ Span of } \{\partialX^\pm_{j;b,n}(\lambda)\}_{b\in S_j^{\pm}, n\in\N}
\]
where as in Lemma \ref{lem: xpm-partial}, $S_j^{\pm}\subset\C$ are a choice of 
representatives of poles of $\eX_j^{\pm}(u,\lambda)$ in the spectral variable $u$,
modulo $\Lambda_{\tau}$.\\

Let $\ulA\in\Mer(\mu;\eV)$ and $b\in S_j^\pm$. Define $\ulA^\pm$ by:
\[
A^{\pm}_i(u) = A_i(u) \frac{\theta(u-b\pm a)}{\theta(u-b\mp a)}
\]

Proposition \ref{prop: important-prop} implies that 
\[
\partialX^\pm_{j;b,n}(\lambda)\lp \eV_\mu[\ulA]\rp \subset \bigoplus_{\ulB\sim\ulA^{\pm}}
\eV_{\mu\pm\alpha_j}[\ulB]
\]
The right--hand side has only one summand, by the definition of \lflat object. Thus
Proposition \ref{prop: important-prop} also gives that $\partialX^{\pm}_{j;b,n}(\lambda)$
on $\eV_\mu[\ulA]$ has the form $e^{2\pi\iota (\lambda,\alpha)}\partialX^{\pm}_{j;b,n}(0)$,
for some $\alpha\in Q$. Thus we get
\[
\text{ Span of } \{\partialX^\pm_{j;b,n}(\lambda)\}_{\lambda\in\hres,b\in S_j^{\pm}, n\in\N}
 = \text{ Span of } \{\partialX^\pm_{j;b,n}(0)\}_{b\in S_j^{\pm}, n\in\N}
\]
in $\Hom_{\C}(\eV_\mu, \eV_{\mu\pm\alpha_j})$. This proves the lemma.
\end{pf}

\subsection{}\label{ssec: weakpbw2}

Throughout this paper,
we say $\eV\in\eL$ is {\em irreducible} if there is no proper non--zero subspace of $\eV$
which is stable under $\{\Phi_i(u),\eX^{\pm}_i(u,\lambda)\}_{i\in\bfI}$. Later,
in Corollary \ref{cor:generalinjective},
we will show that this naive notion of irreducibility agrees with the more general
one, for any category with an initial object.\\

Now we are in the position to prove the following

\begin{thm}\label{thm: eqg-pbw}
Let $\eV\in\eL$ be an irreducible, \lflat object. Then there exists a
non--zero, unique
up to scalar, weight vector $\vac\in\eV_\mu$ such that
\begin{itemize}
\item $\eX_i^{+}(u,\lambda)\vac = 0$ for every $i\in\bfI, u\in\C, \lambda\in\hres$.
\item $\vac$ is an eigenvector for $\Phi_i(u)$ for every $i\in\bfI$.
\item $\eV$ is the span of the following vectors:
\[
\eX_{i_1}^-(u_1,\lambda+\hbar\alpha_{i_1})
\cdots\eX_{i_l}^-(u_l,\lambda+\hbar(\alpha_{i_1}+\cdots+\alpha_{i_l}))\vac
\]
where $u_1,\cdots,u_l\in\C$, $\lambda\in\hres$, $l\in\N$ and $i_1,\cdots,i_l\in\bfI$.
\end{itemize}
\end{thm}

\begin{pf}
By the category $\mathcal{O}$ condition, we can find a weight space $\eV_\mu$
such that $\eV_{\mu+\alpha_i}=0$ for every $i\in\bfI$. Let $\vac\in\eV_\mu$
be a non--zero eigenvector for $\{\Phi_i(u)\}_{i\in\bfI}$. Thus there exists
$\ulA\in\Mer(\mu)$ such that $\Phi_i(u)\vac = A_i(u)\vac$ for every 
$i\in\bfI$. Moreover $\eX_i^+(u,\lambda)\vac = 0$ since $\mu+\alpha_i$
is not a weight of $\eV$.\\

Let $W\subset \eV$ be the subspace spanned by the following vectors:
\[
\eX_{i_1}^-(u_1,\lambda+\hbar\alpha_{i_1})
\cdots\eX_{i_l}^-(u_l,\lambda+\hbar(\alpha_{i_1}+\cdots+\alpha_{i_l}))\vac
\]
over all $u_1,\cdots,u_l\in\C$, $\lambda\in\hres$ and $i_1,\cdots,i_l\in\bfI$.
Using the relation (\EQ 3) and the fact that $\vac$ is an eigenvector for $\Phi_i(u)$
we get that $W$ is stable under $\Phi_i(u)$ for every $i\in\bfI$.\\

Now we claim that $W$ is also stable under the raising and lowering operators.
Lemma \ref{lem: eqg-pbw} implies that we only need to prove this assertion
for the operator obtained by specializing the dynamical variable of our choice.
This makes it clear that $W$ is stable under $\eX^-_i(u,\lambda)$. For
the raising operators, again we have
\begin{multline*}
\eX_i^+(u,\lambda')\eX_{i_1}^-(u_1,\lambda+\hbar\alpha_{i_1})
\cdots\eX_{i_l}^-(u_l,\lambda+\hbar(\alpha_{i_1}+\cdots+\alpha_{i_l}))\vac
 = \\
\sum_{t=1}^l \cdots \left[ \eX_i^+(u,\lambda'),\eX^-_{i_t}(u_t,\lambda+\hbar(
\alpha_{i_1}+\cdots+\alpha_{i_t})\right]\cdots \vac
\end{multline*}
For $\lambda'$ such that $\lambda'+\lambda = \hbar(\mu-\alpha_{i_1}-\cdots-\alpha_{i_l})$
we can apply the relation (\EQ 5) to each term of the summation on the \rhs
above and conclude that $W$ is stable under $\eX_i^+(u,\lambda')$.\\

Thus $W$ is a subobject of $\eV$. Since $\eV$ is assumed to be irreducible,
this implies that $W = \eV$ as required.
\end{pf}

\subsection{Composition series}\label{ssec: kac-lemma}

We have the following analogue of the existence of composition series
in category $\mathcal{O}$ for a symmetrisable Kac--Moody algebra
\cite[Lemma 9.6]{kac}. We refer the reader to \cite[Prop. 15]{hernandez-affinizations}
for a similar statement for quantum loop algebras.

\begin{lem}\label{lem:kac}
Let $\eV\in\eL$ be \lflat and $\mu\in\h^*$. Then, there exists a filtration

\[0=\eV_0\subset\cdots\subset\eV_t=\eV\]

such that the following holds for any $\eV_j$, $j=1,\ldots,t$
\begin{itemize}
\item either $\eV_j/\eV_{j-1}$, is an irreducible object with 
highest weight $\mu_j$ for some
$\mu_j\geq \mu$,
\item or $\lp \eV_j/\eV_{j-1}\rp_{\nu} = 0$ for every $\nu\geq\mu$.
\end{itemize}
\end{lem}

\begin{pf}
The proof of this lemma is essentially the same as the one given
in \cite{kac}. Namely, we define 
\[
a(\mu,\eV) = \sum_{\nu\geq \mu} \dim(\eV_{\nu})
\]
Note that $a(\mu,\eV)$ is finite by the category $\mathcal{O}$
axiom from Section \ref{ssec: categoryL}.\\

We prove the lemma by induction on $a(\mu,\eV)$. Assuming
$a(\mu,\eV) = 0$, the filtration $0\subset\eV$ satisfies the conditions
of the lemma. Let $a(\mu,\eV)>0$ and choose a weight $\mu_1\geq\mu$ of $\eV$
such that $\eV_{\mu+\alpha_i}=0$ for every $i\in\bfI$. Let $\vac_1\in\eV_{\mu_1}$ be an 
eigenvector for $\{\Phi_i(u)\}_{i\in\bfI}$ and consider the subspace
$W$ spanned by 
\[
\eX_{i_1}^-(u_1,\lambda+\hbar\alpha_{i_1})
\cdots\eX_{i_l}^-(u_l,\lambda+\hbar(\alpha_{i_1}+\cdots+\alpha_{i_l}))\vac_1
\]
over all $u_1,\cdots,u_l\in\C$, $\lambda\in\hres$ and $i_1,\cdots,i_l\in\bfI$.

The argument given in the proof of Theorem \ref{thm: eqg-pbw} yields
that $W$ is stable under $\{\Phi_i(u),\eX_i^{\pm}(u,\lambda)\}$ and hence
defines a subobject of $\eV$. Note that $\dim(W_{\mu_1})=1$
and for any proper subspace $W'$
of $W$ which is stable under these operators, we must have $\dim(W'_{\mu_1})=0$. 
Let $W_1\subset W$ be the largest such proper subspace.
Then we arrive at a filtration
$0\subset W_1\subset W\subset \eV$ where $W/W_1$ is irreducible and
$a(\mu,W_1), a(\mu,\eV/W) < a(\mu,\eV)$. Hence we are done by induction.

\end{pf}

\section{Classification of irreducibles I: Necessary condition}\label{sec: irr-class1}

In this section we give a necessary condition for the eigenvalues of $\{\Phi_i(u)\}$
on the highest weight vector of an irreducible object in $\eL$.

\subsection{}\label{ssec:class1}

Let $\eV\in\eL$ be an irreducbile, \lflat object. By Theorem \ref{thm: eqg-pbw}, we 
have a unique (up to scalar) weight vector $\vac\in\eV_\mu$ such that
\begin{itemize}
\item $\eX_i^+(u,\lambda)\vac = 0$ for every $i\in\bfI, u\in\C, \lambda\in\hres$.\\

\item There exist $\ulA\in\Mer(\mu)$ such that $\Phi_i(u)\vac = A_i(u)\vac$
for every $i\in\bfI$.\\

\item $\eV$ is spanned by the following vectors
\[
\eX_{i_1}^-(u_1,\lambda+\hbar\alpha_{i_1})
\cdots\eX_{i_l}^-(u_l,\lambda+\hbar(\alpha_{i_1}+\cdots+\alpha_{i_l}))\vac
\]
over all $u_1,\cdots,u_l\in\C$, $\lambda\in\hres$ and $i_1,\cdots,i_l\in\bfI$.
\end{itemize}

\begin{thm}\label{thm: class1}
For every $i\in\bfI$, 
there exist $\ul{N} = (N_i)_{i\in\bfI}\in \N^{\bfI}$ and collection of complex numbers
$\{b^{(i)}_1,\cdots, b^{(i)}_{N_i}\}$ such that
\[
A_i(u) = C_i \prod_{k=1}^{N_i} \frac{\theta(u-b^{(i)}_k+d_i\hbar)}{\theta(u-b^{(i)}_k)}
\]
where $C_i$ is a non--zero constant.
In particular $\mu(\alpha_i^{\vee}) = N_i$ and hence $\mu$ is a dominant integral weight.
\end{thm}

The proof of this theorem is given in Section \ref{ssec: thm-class1-pf} below.

\subsection{}\label{ssec:class1-rems}

We remark that the statement of Theorem \ref{thm: class1} for Yangians and quantum loop algebras
is part of the classification of their irreducible representations by Drinfeld polynomials. For
$\g=\sl_2$ the classification theorem is proved in 
\cite{chari-pressley-yangian, chari-pressley-qaffine}.
In \cite{chari-pressley-qaffine-rep} and \cite{hernandez-affinizations}
the theorem is proved for an arbitrary \fd semisimple Lie algebra and symmetrisable
Kac--Moody algebra respectively (see also \cite[Chapter 12]{chari-pressley}, \cite[Chapter 3]{molev-yangian}
and references therein).\\

The proof of Theorem \ref{thm: class1} given here differs significantly from the ones that
exist in the literature due to the lack of standard constructions of PBW basis, Verma
modules, tensor product etc. Namely, we produce an infinite sequence of lowering
operators which never annihilate the highest--weight vector, if the contrary
of Theroem \ref{thm: class1} were true. Lemma \ref{lem:class1-crucial} is crucial
in achieving this. We merely wish to remark that the proofs of this lemma
and Theorem \ref{thm: class1} can be easily modified to work for Yangians and quantum loop
algebras. Our proof uses techniques which are familiar in the theory of 
quantum loop algebras, see \eg \cite{young}.\\

\Omit{It is worth noticing that Theorem \ref{thm: class1} has the following weakness. It
yields a set--theoretic map from the isomorphism classes of irreducible objects in $\eL$
to the set $\PP_+^{\E}$ defined below:

\[
\PP_+^{\E} := \left\{ \lp \mu, \{b^{(i)}_1,\cdots,b^{(i)}_{N_i}\}_{i\in\bfI}\rp : \mu\in P_+, N_i = \mu(\alpha_i^{\vee}),
b^{(i)}_1,\cdots,b^{(i)}_{N_i} \in \mathsf{P}\right\}
\]
where $\mathsf{P}\subset\C$ is a fundamental domain for the action of $\Lambda_\tau$ on $\C$ via translations.
However Theorem \ref{thm: class1} makes no claim as to injectivity or surjectivity of this map.
The main theorem of Section \ref{sec: class2} will precisely yield the bijection
between $\PP_+^{\E}$ and the set of isomorphism classes of irreducible objects of
$\eL$.}

\subsection{Choice of representatives of poles}\label{ssec:class1-poles}

For $i\in\bfI$ we define $\eV^{(i)} = \oplus_{n\geq 0} \eV_{\mu-n\alpha_i}$. Note that this is a \fd
vector space which is stable under the action of $\eX^\pm_i(u,\lambda)$ and $\{\Phi_j(u)\}_{j\in\bfI}$.
Let $\spec(\eV,i)$ denote the set of poles of $\{\Phi_i(u),\eX^{\pm}_i(u,\lambda)\}$ acting
on $\eV^{(i)}$ in the spectral variable $u$. 
We choose a (finite) subset $\prps{i}\subset\C$ satisfying the following conditions

\begin{itemize}
\item For every $a\in\spec(\eV,i)$ there is a unique $\ol{a}\in\prps{i}$ such that $a-\ol{a}\in\Lambda_{\tau}$.

\item If $a\in\spec(\eV,i)\cap\prps{i}$ and $a+n\hbar_i\in\spec(\eV)$ for some $n\in\Z^{\times}$ then $a+n\hbar_i\in\prps{i}$.
\end{itemize}

Here and below we use the notation $\hbar_i = d_i\hbar$.

\subsection{Partial fractions}\label{ssec:class1-partial}

Considering the principal parts of the Laurent series expansions of $\Phi_i(u), \eX_i^{\pm}(u,\lambda)$ near
a point $a\in\prps{i}$ we can write the partial fraction decomposition of these functions (see Lemma
\ref{lem: xpm-partial}) as follows:
\begin{align*}
\eX_i^{\pm}(u,\lambda) &= \sum_{\begin{subarray}{c} a\in\prps{i} \\ n\in\N\end{subarray}}
\partialX^{\pm}_{i;a,n}(\lambda) \frac{(-\partial_u)^n}{n!}\left(\ttheta{u-a+\lambda_i}{u-a}{\lambda_i}\right)
\end{align*}

For each $a\in\prps{i}$ let us define $\partialH_{i;a,n}$ using the Laurent series expansion of $\Phi_i(u)$
near $u=a$:
\[
\Phi_i(u) = \sum_{n\in\N} \frac{\partialH_{i;a,n}}{(u-a)^{n+1}} + \text{ part holomorphic near } a
\]

Then we can rewrite the relations between $\{\Phi_i(u), \eX_i^{\pm}(u,\lambda)\}$ 
considered as meromorphic functions taking values in $\End\lp\eV^{(i)}\rp$,
in the following manner:

\begin{lem}\label{lem:class1-partial}
\begin{enumerate}
\item For each $a\in\prps{i}$ and $n\in\N$, we have:
\begin{equation}\label{aeq:hx}
\Ad(\Phi_j(u)).\partialX^{\pm}_{i;a,n}(\lambda) = \sum_{l\geq 0}\left.
\frac{\partial_v^l}{l!}\lp\frac{\theta(u-v\pm a_{ij}\hbar_i/2)}{\theta(u-v\mp a_{ij}\hbar_i/2)}\rp
\right|_{v=a} \partialX^{\pm}_{i;a,n+l}(\lambda\pm \hbar\alpha_j)
\end{equation}

\item For each $a\in\prps{i}$ and $m,n\in\N$ we have

\begin{multline}\label{aeq:x++}
\sum_{k,l\geq 0} \frac{(-1)^l}{k!l!} 
\lp\theta^{(k+l)}(\mp\hbar_i)\partialX^{\pm}_{i;a,n+k}(\lambda)
\partialX^{\pm}_{i;a,m+l}(\lambda\mp\hbar\alpha_i) - \right.\\
\left.\theta^{(k+l)}(\pm\hbar_i)\partialX^{\pm}_{i;a,m+l}(\lambda)
\partialX^{\pm}_{i;a,n+k}(\lambda\mp\hbar\alpha_i)\rp = 0
\end{multline}

\item For each $a,b\in\prps{i}$, $m,n\in\N$ we have the following relation
on a weight space $\eV_{\nu}$:
\begin{equation}\label{aeq:x+-}
\theta(\hbar_i)\left[
\partialX^+_{i;a,m}(\lambda),\partialX^-_{i;b,n}(-\lambda+\hbar\nu)
\right]
 = \delta_{a,b}\partialH_{i;a,m+n}
\end{equation}
\end{enumerate}
\end{lem}

\begin{pf}
(1) was already proved in Lemma \ref{lem: xpm-difference}.\\

(2): Consider the relation (\EQ 4) or its equivalent form (\EQ $4'$) given 
in Section \ref{ssec: categoryL} for $i=j$ and $\lambda$ replaced by
$\ds\lambda-\frac{\hbar}{2}\alpha_i$ (for notational convenience only), 
multiply it by $(u-a)^n(v-a)^m$. Integrate with respect to both $u$ and $v$
over a small circle centered at $a$, say $\cC$. Note that the second and third terms
from both sides vanish and we are left with
\begin{align*}
\text{L.H.S. } &= \oint_{\cC}\oint_{\cC} (u-a)^n(v-a)^m\theta(u-v\mp\hbar_i)
\eX_i^{\pm}(u,\lambda)\eX_i^{\pm}(v,\lambda\mp\hbar\alpha_i)\, du\, dv \\
&= \sum_{l\geq 0}\oint_{\cC} (u-a)^n \frac{(-1)^l}{l!} \theta^{(l)}(u-a\mp\hbar_i)
\eX_i^{\pm}(u,\lambda)\, du\, \partialX_{i;a,m+l}^{\pm}(\lambda\mp\hbar\alpha_i)\\
&= \sum_{k,l\geq 0}\frac{(-1)^l}{k!l!} \theta^{(k+l)}(\mp\hbar_i)
\partialX^{\pm}_{i;a,n+k}(\lambda)\partialX^{\pm}_{i;a,m+l}(\lambda\mp\hbar\alpha_i)
\end{align*}

The same calculation with the \rhs gives \eqref{aeq:x++}.\\

(3): Again we consider (\EQ 5) with $i=j$ on a weight space $\eV_{\nu}$.
Multiply both sides by $(v-b)^n$ and integrate over a small circle $\cC$
around $b$, to get:
\[
\theta(\hbar_i)\left[\eX_i^+(u,\lambda),
\partialX^-_{i;b,n}(-\lambda+\hbar\nu)\right]
 = \oint_{\cC} (v-b)^n\ttheta{u-v+\lambda_i}{u-v}{\lambda_i} \Phi_i(v)\, dv
\]
Now we multiply this by $(u-a)^m$ and integrate around $\cC'$,
a small circle around $a$. Clearly the function on the \rhs of the equation
above only has poles at $u=a$ (modulo $\Lambda_{\tau}$). Therefore, if 
$a\not= b\in\prps{i}$, the integral will be zero. Assume $b=a$ and
that $\cC'$ is a slight enlargement of $\cC$. Then we get
\begin{align*}
\theta(\hbar_i)\left[\partialX^+_{i;a,m}(\lambda),\partialX^-_{i;a,n}
(-\lambda+\hbar\nu)\right] &=
\oint_{\cC}\oint_{\cC'}(u-a)^m\ttheta{u-v+\lambda_i}{u-v}{\lambda_i}\, du\, (v-a)^n \Phi_i(v)\, dv\\
&= \oint_{\cC} (v-a)^{n+m}\Phi_i(v)\, dv \\
&= \partialH_{i;a,n+m}
\end{align*}
\end{pf}

\subsection{}\label{ssec:class1-straight}

As a consequence of \eqref{aeq:x++} we have the following

\begin{cor}\label{cor:class1-straight}
For every $m,n\in\N$ with $m\geq n$, the following relation holds
\[
\partialX^{\pm}_{i;a,m}(\lambda)\partialX^\pm_{i;a,n}(\lambda\mp\hbar\alpha_i) = 
\sum_{n\leq r<s} C^{\pm}_{m,n}(r,s)\partialX^\pm_{i;a,r}(\lambda)
\partialX^\pm_{i;a,s}(\lambda\mp\hbar\alpha_i)
\]
where $C^{\pm}_{m,n}(r,s)\in\C$.
\end{cor}

\begin{pf}
The idea of the proof is to use \eqref{aeq:x++} and the fact that there
exists $N>0$ such that $\partialX^{\pm}_{i;a,l}(\lambda)=0$ for every 
$l>N$. Let us assume that $N$ is smallest such positive integer. Then by
\eqref{aeq:x++} for $m=n=N$ we have
\[
2\theta(\mp\hbar_i)\partialX^{\pm}_{i;a,N}(\lambda)\partialX^{\pm}
_{i;a,N}(\lambda\mp\hbar\alpha_i) = 0
\]

In order to prove the assertion of the corollary for $m=N$ and $n<N$,
we again use \eqref{aeq:x++} to get
\[
\sum_{k\geq 0} \frac{1}{k!}\theta^{(k)}(\mp\hbar)
\partialX_{i;a,n+k}^{\pm}(\lambda)\partialX_{i;a,N}^{\pm}(\lambda\mp\hbar\alpha_i)
 = 
\sum_{k\geq 0} \frac{1}{k!}\theta^{(k)}(\pm\hbar)
\partialX^{\pm}_{i;a,N}(\lambda)\partialX^{\pm}_{i;a,n+k}(\lambda\mp\hbar\alpha_i)
\]

An easy induction argument proves the corollary for $m=N$ and arbitrary $n$.
Now we assume that we have proved it for $m>m_1$ and every $n\leq m$, and 
continue with proving it for $m=m_1$. Proceeding just as above we have

\begin{multline*}
\sum_{k,l\geq 0} \frac{(-1)^l}{k!l!} 
\theta^{(k+l)}(\mp\hbar_i)\partialX^{\pm}_{i;a,m_1+k}(\lambda)
\partialX^{\pm}_{i;a,m_1+l}(\lambda\mp\hbar\alpha_i)\\
=\sum_{k,l\geq 0}
\theta^{(k+l)}(\pm\hbar_i)\partialX^{\pm}_{i;a,m_1+l}(\lambda)
\partialX^{\pm}_{i;a,m_1+k}(\lambda\mp\hbar\alpha_i) = 0
\end{multline*}

The induction hypothesis implies that the claim holds for $n=m=m_1$.
The general case with $n<m_1$ then follows just as before (for the 
case of $m=N$).

\end{pf}

\subsection{Main Lemma}\label{ssec:class1-crucial}

\begin{lem}\label{lem:class1-crucial}
Let $a\in\prps{i}$ and $N\in\N$. Assume there exists a non--zero weight vector
$v\in\eV^{(i)}_{\nu}$ satisfying the following two conditions.
\begin{itemize}
\item $\partialX^+_{i;a,k}(\lambda)v = 0$ for every $k\in\N$ and $\lambda\in\C$.
\item $\Phi_j(u)v = B_j(u)v$ for some $\ulB=(B_j(u))\in\Mer(\nu)$ such that $B_i(u)$
has a pole of order $N+1$ at $a$.
\end{itemize}
Then we have
\[
\partialX^-_{i;a,0}(\lambda)\cdots\partialX^-_{i;a,N}(\lambda+N\hbar\alpha_i)v\neq 0
\]
for every $\lambda\in\hres$.
\end{lem}

\begin{pf}
For the purposes of the proof, let us keep $i\in\bfI$ and $a\in\prps{i}$ fixed and hence omit these from
the subscript.\\

We begin by making the following observation which directly follows from the proofs
of Proposition \ref{prop: important-prop} and Lemma \ref{lem: eqg-pbw}:\\

There exist $\alpha^{(1)},\alpha^{(2)},\cdots\in Q$ such that, if we define, $\ulB^{(r)}\in\Mer(\nu-r\alpha_i)$
by:
\[
B^{(r)}_j(u) = B_j(u) \lp \frac{\theta(u-a-d_ia_{ij}\hbar/2)}{\theta(u-a+d_ia_{ij}\hbar/2)}\rp^r
\exp\lp-2\pi\iota\hbar \lp\alpha_j,\alpha^{(1)}+\cdots+\alpha^{(r)}\rp\rp
\]

Then we have

\begin{equation}\label{aeq:obs}
\partialX^-_{i;a,n}(\lambda) = e^{2\pi\iota (\lambda,\alpha^{(r)})}\partialX^-_{i;a,n}(0)
: \eV_{\nu-(r-1)\alpha_i}[\ulB^{(r-1)}] \to \eV_{\nu-r\alpha_i}[\ulB^{(r)}]
\end{equation}

Now define $v(k_1,\cdots,k_r;\lambda)\in\eV_{\nu-r\alpha_i}[\ulB^{(r)}]$ by
\[
v(k_1,\cdots,k_r;\lambda) := \partialX^-_{i;a,k_1}(\lambda)\cdots\partialX^-_{i;a,k_r}(\lambda+(r-1)\hbar\alpha_i)v
\]

Note that by \eqref{aeq:obs} above, we have

\begin{equation}\label{aeq:obs2}
v(k_1,\cdots,k_r;\lambda) = \exp\lp 2\pi\iota\lp\lambda,\alpha^{(1)}+\cdots+\alpha^{(r)}\rp\rp v(k_1,\cdots,k_r;0)
\end{equation}

Let us consider the subspace $\wt{W}$ of $\eV$ spanned by $v(k_1,\cdots,k_r;\lambda)$ for all
$k_1,\cdots,k_r\in\N$ and $\lambda\in\hres$ (for $r=0$, this vector is just $v$). By the
straightening relation given in Corollary \ref{cor:class1-straight} we have
\[
\wt{W} = \text{ Span of }\{v(k_1,\cdots,k_r;\lambda)\}_{0\leq k_1<\cdots<k_r}
\]
Let $W'\subset \wt{W}$ be the span of such vectors with $k_r>N$, and let $W = \wt{W}/W'$ be the
quotient space. Define, for each $n\in\N$ and $0\leq p\leq n$:

\[
v_{p,n}(\lambda) := v(N-n,\cdots,N-n+p-1,\wh{N-n+p},N-n+p+1,\cdots,N;\lambda)
\]
where as usual $\hat{x}$ means $x$ is skipped. We prove the following assertion by induction
on $n$ (for $0\leq n\leq N$):

\noindent ${\mathbf S}_n$: For any fixed $\lambda\in\hres$, the collection of vectors
$\{v_{p,n}(\lambda)\}_{0\leq p\leq n}$ is linearly independent in the quotient space $W$.\\

Note that the vector $v_{0,n}(\lambda)$ is an eigenvector for $\{\Phi_j(u)\}$ by
\eqref{aeq:hx} and Corollary \ref{cor:class1-straight}. That is,

\begin{equation}\label{aeq:obs3}
\Phi_j(u)v_{0,n}(\lambda) = B_j(u)\lp\frac{\theta(u-a-d_ia_{ij}\hbar/2)}{\theta(u-a+d_ia_{ij}\hbar/2)}\rp^n
v_{0,n}(\lambda-\hbar\alpha_j) \text{ modulo } W'
\end{equation}

In particular, for each $l\geq 0$ we have a complex number $\mathit{h}_{l,n}\in\C$ such that

\begin{equation}\label{aeq:obs4}
\partialH_{i;a,l}v_{0,n}(\lambda) = \mathit{h}_{l,n}v_{0,n}(\lambda-\hbar\alpha_i) \text{ modulo } W'
\end{equation}

The base case ${\mathbf S}_0$ just means that $v\neq 0$ which is true by our hypothesis and the fact
that the subspace $W'$ is spanned by weight vectors of weight strictly less than $\nu$ (hence
$v\not\in W'$).\\

Now we carry out the induction step. Let us first assume ${\mathbf S}_n$ for each $0\leq n\leq N$. We claim
that $v_{0,n+1}(\lambda)\neq 0$. Note that when $n=N$, this is exactly the statement we need to prove.
Assume that this vector is zero. Then applying $\partialX^+_{i;a,0}(-\lambda+\hbar(\nu-(n+1)\alpha_i))$ and
using the relation \eqref{aeq:x+-}, we get

\begin{multline*}
\sum_{p=0}^n \partialX^-_{i;a,N-n}(\lambda)\cdots\partialX^-_{i;a,N-n+p-1}(\lambda+(p-1)\hbar\alpha_i)
\partialH_{i;a,N-n+p} \\
\partialX^-_{i;a,N-n+p+1}(\lambda+(p+1)\hbar\alpha_i)\cdots \partialX^-_{i;a,N}(\lambda+n\hbar\alpha_i)v = 0
\end{multline*}

Note that 
$$\partialX^-_{i;a,N-n+p+1}(\lambda+(p+1)\hbar\alpha_i)\cdots \partialX^-_{i;a,N}(\lambda+n\hbar\alpha_i)v
= v_{0,N-n+p}(\lambda+(p+1)\hbar\alpha_i)$$ 
and this vector is an eigenvector for $\partialH$ (as given in \eqref{aeq:obs4}). Thus we get

\[
\sum_{p=0}^n \mathit{h}_{N-n+p,N-n+p} v_{p,n}(\lambda) = 0
\]

This is a non--trivial dependence relation, since $\mathit{h}_{N,N}\neq 0$ by our hypothesis on $B_i(u)$
(it has a pole of order $N+1$ at $a$), which contradicts ${\mathbf S}_n$.\\

Assuming $n<N$, we can now prove ${\mathbf S}_{n+1}$. The proof relies on the following
easy computation, using \eqref{aeq:hx} and Corollary \ref{cor:class1-straight}:

\begin{multline*}
\Phi_i(u)v_{p,n}(\lambda) = B^{(n)}_i(u)v_{p,n}(\lambda) + \\
B^{(n)}_i(u)\frac{\theta(u-a+\hbar_i)}{\theta(u-a-\hbar_i)}\lp\left.\partial_v\lp
\frac{\theta(u-v-\hbar_i)}{\theta(u-v+\hbar_i)}\rp\right|_{v=a}\rp v_{p-1,n}(\lambda) + \ldots
\end{multline*}

where $\ldots$ refers to terms $v_{j,n}$ with $j<p-1$. This triangularity property of $\Phi_i(u)$
together with the fact that $v_{0,n+1}\neq 0$ implies the required result, using the following
general fact from linear algebra.
\end{pf}

\subsection{}
\begin{lem}\label{alem:easy}
Let $V$ be a vector space over $\C$ and $\{v_1,\cdots, v_l\}\subset V$ such that $v_1\neq 0$.
Assume that there exists
$X\in\End(V)$ such that $Xv_i = \sum_{j\leq i} c_{ij} v_j$ where
\begin{itemize}
\item $c_{ii} = c$ is independent
of $i$.
\item $c_{i-1,i}\neq 0$ for each $i$.
\end{itemize}
Then $v_1,\cdots, v_l$ are linearly independent.
\end{lem}

\subsection{}\label{assc:crucial-cor}

Let us note an important consequence of Lemma \ref{lem:class1-crucial}. 

Define a subspace of $\eV$:
\[
V' := \text{ Span of } \{\partialX^-_{i;a,k_0}(\lambda)\cdots 
\partialX^-_{i;a,k_N}(\lambda+N\hbar\alpha_i)v\}_{k_0,\cdots,k_N\in\N}
\]
This subspace is stable under $\{\Phi_j(u)\}$ because of the relation \eqref{aeq:hx}
of Lemma \ref{lem:class1-partial}. Moreover
by Proposition \ref{prop: important-prop} we have the following inclusion
\[
V' \subset \eV^{(i)}_{\nu-(N+1)\alpha_i}\left[\wt{\ulB}\right]
\]
where $\wt{\ulB} = \ulB^{(N+1)}$ in the notation introduced in the proof of 
Lemma \ref{lem:class1-crucial}.\\

The fact that this subspace is non--zero follows from Lemma \ref{lem:class1-crucial}
and hence implies the following

\begin{cor}\label{cor:class1-crucial}
There exists a non--zero vector $v'\in V'$ such that
$\Phi_j(u)v' = \wt{B}_j(u)v'$ for every $j\in\bfI$.
\end{cor}

\subsection{Rephrasing Drinfeld polynomials}\label{ssec:dr-rephrase}

For us it will be important to identify meromorphic functions from
the location of their zeroes and poles, which are of the following form
\begin{equation}\label{aeq:right-form}
\text{C }\prod_k \frac{\theta(u-c_k+\hbar_i)}{\theta(u-c_k)}
\end{equation}
where $c_k\in\C$ and $\text{C }$ is a non--zero constant.
A little bit of terminology is in order.
A subset $S\subset\C$ is called a {\em $\hbar_i$--string} if for every $s_1,s_2\in S$, we have
$s_1-s_2\in\Z\hbar_i$. Every finite subset of $\C$ can clearly be broken into a finite
union of $\hbar_i$--strings. \\

Given a meromorphic function $B(u)$ such that
\[
B(u+1) = B(u) \text{ and } B(u+\tau) = e^{-2\pi\iota \nu_i}B(u)
\]
for some $\nu_i\in\C$,
let $\sigma(B(u))$ be the set of its zeroes and poles and let
$S_B\subset\C$ denote a choice of representatives of the zeroes and poles
as outlined in Section \ref{ssec:class1-poles}. Let
\[
S_B = S_1\sqcup \ldots \sqcup S_m
\]
be the decomposition of $S_B$ into $\hbar_i$--strings. To each string $S_k$ ($1\leq k\leq m$)
we associate an expression of parantheses as follows. Let $s\in S_k$ be such that $s+n\hbar_i\not\in S_k$
for every $n\geq 1$. If $s$ is a zero (resp. pole) of $B(u)$ of order $n$, we write $n$
right (resp. left) parantheses. Then we continue with $S_k\setminus\{s\}$ and so on until
the $\hbar_i$--string $S_k$ is exhausted. We say that $S_k$ is {\em balanced} if the resulting
expression is of balanced parantheses. The following is immediate.

\begin{lem}\label{lem:class1-dr}
$B(u)$ is of the form \eqref{aeq:right-form}
if, and
only if $S_k$ is balanced for each $1\leq k\leq m$.
\end{lem}

\subsection{Proof of Theorem \ref{thm: class1}}\label{ssec: thm-class1-pf}

Let us assume that there exists $i\in\bfI$ such that the eigenvalue of $\Phi_i(u)$ on $\vac$,
$A_i(u)$ as introduced in the statement of Theorem
\ref{thm: class1} is not of the form \eqref{aeq:right-form}. Decompose the finite set of
representative of zeroes and poles of $A_i(u)$, denoted in the previous section by $S_{A_i}$,
into $\hbar_i$--strings. According to Lemma \ref{lem:class1-dr}, one of these $\hbar_i$--strings
is unbalanced. Let $S$ be an unbalanced $\hbar_i$--string in which number of left parantheses
is greater than or equal to that of right parantheses (such $S$ must exist since number
of zeroes and poles of $A_i(u)$ are the same). We pick $b\in S$ according to the following
procedure:\\

Let $E$ be the expression of parantheses associated to $S$. By the standard algorithm to
determine whether $E$ is balanced or not, one defines a number $\text{Counter}=0$ and
reads the expression from right to left. If a paranthesis is right (resp. left) we
decrease (resp. increase) the number Counter by $1$. Then $E$ is unbalanced if, and
only if the Counter becomes positive at some stage (or if upon exhausing $E$ it is
non--zero, which will be positive for us, since we are assuming that $E$ has
more or equal number of left paratheses). Let $x$ be the paranthesis at which the
counter becomes positive for the first time, and let $b\in S$ be the correspoding
complex number
(necessarily a pole of $A_i(u)$, since left parantheses were associated to poles).\\

Now we can finish the proof of Theorem \ref{thm: class1} as follows. Let $N_0$ be the
order of pole of $A_i(u)$ at $b_0 = b$. Applying Lemma \ref{lem:class1-crucial} and Corollary
\ref{cor:class1-crucial} to $v = \vac$, $a=b_0$ and $N = N_0-1$ we obtain
another (non--zero) vector $\vac_1$ satisfying the following properties:
\begin{itemize}
\item
$\displaystyle \vac_1\in \text{Span of } \{\partialX^-_{i;b_0,k_1}(\lambda)\cdots\partialX^-_{i;b_0,k_{N_0}}(\lambda+\hbar(N_0-1)\alpha_i)v\}_
{k_1,\cdots,k_{N_0}\in\N}$. \\

\item $\vac_1$ is an eigenvector for $\{\Phi_j(u)\}$ with eigenvalue $\ulA^{(1)}$ where
\[
A^{(1)}_j(u) = A_j(u)\lp\frac{\theta(u-b_0-d_ia_{ij}\hbar/2)}{\theta(u-b_0+d_ia_{ij}\hbar/2)}\rp^{N_0}.e^{-2\pi\iota\hbar (\alpha_j,\alpha)}
\]
for some $\alpha\in Q$.
\end{itemize}

The first property in particular implies that $\partialX^+_{i;c,k}(\lambda)\vac_1 = 0$
for every $c\neq b_0$, $k\in\N$ and $\lambda\in\hres$. Moreover the function $A^{(1)}_i(u)$ has
a pole at $b_1 = b_0-\hbar_i$ of some order, say $N_1$, since otherwise $A_i(u)$ must
have a zero of order $\geq N_0$ at $b_0-\hbar_i$ thus contradicting our choice
of $b$ as outlined above. Hence we can apply Lemma \ref{lem:class1-crucial} and its Corollary
\ref{cor:class1-crucial} again to $a = b_1$, $N = N_1-1$ and $v = \vac_1$.\\

It remains to observe that, according to our choice of $b$, this procedure will never
terminate. Thus we will obtain infinitely many non--zero vectors $\{\vac_n\}_{n\geq 0}$
which belong to different weight spaces and hence are linearly independent. This contradicts
the finite--dimensionality of $\eV^{(i)}$ and completes the proof of Theorem \ref{thm: class1}.

\subsection{Knight's lemma}\label{ssec: eqg-knight}

The following result follows directly from Lemma \ref{lem:kac}, Proposition \ref{prop: important-prop} and 
Theorems \ref{thm: eqg-pbw}, \ref{thm: class1}. The analogous assertion for
Yangians appeared in \cite{knight}.

\begin{prop}\label{prop: eqg-knight}
Let $\mu\in\h^*$ and $\ulA\in\Mer(\mu)$. If there exists $\eV\in\eL$ such that $\eV_\mu[\ulA]\neq 0$
then there exists $\ol{\mu}\in P_+$ such that $\mu\leq \ol{\mu}$ and for each $i\in\bfI$ the function
$A_i(u)$ is of the following form:
\[
A_i(u) = C_i \prod_{k=1}^{N_1} \frac{\theta(u-c_{i,k}+d_i\hbar)}
{\theta(u-c_{i,k})} \prod_{l=1}^{N_2}\frac{\theta(u-c'_{i,l}-d_i\hbar)}
{\theta(u-c'_{i,l})}
\]
where $C_i, c_{i,k}, c'_{i,l}\in\C$, and
$N_1,N_2\in\N$ with $N_1-N_2 = \mu(\alpha_i^{\vee})$.
\end{prop}

\section{Quantum loop algebra}\label{sec: qla}

\subsection{The quantum loop algebra $\qloop$}\label{ssec: qla}

Let $q = \exp(\pi\iota\hbar)$.
For any $i\in\bfI$, set $q
_i=q^{d_i}$. We use the standard notation for Gaussian integers
\begin{gather*}
[n]_q = \frac{q^n - q^{-n}}{q-q^{-1}}\\[.5ex]
[n]_q! = [n]_q[n-1]_q\cdots [1]_q\qquad
\qbin{n}{k}{q} = \frac{[n]_q!}{[k]_q![n-k]_q!}
\end{gather*}

The quantum loop algebra $\qloop$ is the $\C$--algebra generated
by $\{\Psi_{i,\pm r}^\pm\}_{i\in\bfI,r\in\N}$, $\{\X_{i,k}^\pm\}
_{i\in\bfI,k\in\Z}$ and $\{K_h\}_{h\in\h}$, subject to the following relations
\begin{itemize}
\item[(QL0)] For any $i\in\bfI$, $\displaystyle \Psi_{i,0}^\pm=K_{\pm d_i\alpha^\vee_i}$. \\

\item[(QL1)] For any $i,j\in\bfI$, $r,s\in\N$ and $h,h'\in\h$,
\begin{gather*}
[\Psi_{i,\pm r}^\pm,\Psi_{j,\pm s}^\pm]=0
\qquad
[\Psi_{i,\pm r}^\pm,\Psi_{j,\mp s}^\mp]=0
\qquad
[\Psi_{i,\pm r}^\pm,K_h]=0\\[.7ex]
K_hK_{h'} = K_{h+h'}
\qquad
K_0=1
\end{gather*}
\item[(QL2)] For any $i\in\bfI$, $k\in\Z$ and $h\in\h$,
$\displaystyle K_h\X^\pm_{i,k}K_h^{-1}=q^{\pm\alpha_i(h)}\X^\pm_{i,k}$.\\

\item[(QL3)] For any $i,j\in\bfI$, $\veps\in\{\pm\}$ and $l\in\Z$
\[\Psi^\veps_{i,k+1}\X^\pm_{j,l} - q_i^{\pm a_{ij}}\X^\pm_{j,l}\Psi^\veps_{i,k+1}
=
q_i^{\pm a_{ij}}\Psi^\veps_{i,k}\X^\pm_{j,l+1}-\X^\pm_{j,l+1}\Psi^\veps_{i,k}\]
for any $k\in\Z_{\geq 0}$ if $\veps=+$ and $k\in\Z_{<0}$ if $\veps=-$
\item[(QL4)] For any $i,j\in\bfI$ and $k,l\in \Z$
\[\X^\pm_{i,k+1}\X^\pm_{j,l} - q_i^{\pm a_{ij}}\X^\pm_{j,l}\X^\pm_{i,k+1}=
q_i^{\pm a_{ij}}\X^\pm_{i,k}\X^\pm_{j,l+1}-\X^\pm_{j,l+1}\X^\pm_{i,k}\]
\item[(QL5)] For any $i,j\in\bfI$ and $k,l\in \Z$
\[[\X^+_{i,k},\X^-_{j,l}] = \delta_{ij} \frac{\Psi^+_{i,k+l} - \Psi^-_{i,k+l}}{q_i - q_i^{-1}}\]
where we set $\Psi^{\pm}_{i,\mp k}=0$ for any $k\geq 1$.
\item[(QL6)] For any $i\neq j\in\bfI$, $m=1-a_{ij}$, $k_1,\ldots, k_m\in\Z$
and $l\in \Z$
\[\sum_{\pi\in \Sym_m} \sum_{s=0}^m (-1)^s\qbin{m}{s}{q_i}
\X^\pm_{i,k_{\pi(1)}}\cdots \X^\pm_{i,k_{\pi(s)}} \X^\pm_{j,l}\X^\pm_{i,k_{\pi(s+1)}}\cdots \X^\pm_{i,k_{\pi(m)}} = 0\]
\end{itemize}

\subsection{Category $\Rloop$}\label{ssec: rloop}

We consider the representations of $\qloop$ whose restriction to $U_q(\g)$
is in category $\mathcal{O}$ and is integrable \cite{hernandez-affinizations}.
Let $\Rloop$ be the category of such representations. We briefly summarize
the results of \cite[\S 2.10--2.13, \S 3.6]{sachin-valerio-2}. (1)
was proved in \cite{beck-kac, hernandez-q-toroidal}.
\begin{enumerate}
\item Define
\begin{align*}
\Psi_i(z)^0 &= \sum_{k\geq 0} \Psi_{i,-k}^-z^k
& \Psi_i(z)^{\infty} &= \sum_{k\geq 0} \Psi^+_{i,k} z^{-k}\\
\X_i^{\pm}(z)^0 &= -\sum_{k\geq 1} \X^{\pm}_{i,-k}z^k
& \X_i^{\pm}(z)^{\infty} &= \sum_{k\geq 0} \X^{\pm}_{i,k} z^{-k}
\end{align*}

Then on a representation $\V\in\Rloop$, there exist rational $\End(\V)$--valued
functions $\{\Psi_i(z), \X_i^{\pm}(z)\}$, regular at $z=0,\infty$ such that
$\Psi_i(z)^{0/\infty}$ and $\X_i^{\pm}(z)^{0/\infty}$ are Taylor expansions
of $\Psi_i(z)$ and $\X_i^{\pm}(z)$ and $z=0/\infty$.

\item Relations (QL1)--(QL5) can be written in terms of these rational functions
(see below).

\item On a representation $\V\in\Rloop$ the relation (QL6) follows from the
rest of the relations.

\end{enumerate}

Thus we can describe objects of $\Rloop$ in the following manner.

\begin{defn}\label{defn: rloop}
An object of $\Rloop$ is an $\h$--diagonalizable module with \fd weight spaces
$\V = \oplus_{\mu\in\h^*} \V_{\mu}$, together with rational $\End(\V)$--valued functions
$\{\Psi_i(z), \X_i^{\pm}(z)\}_{i\in\bfI}$ regular at $z=0,\infty$ satisfying
the following set of axioms:\\

\noindent Category $\mathcal{O}$ and integrability condition.
\begin{itemize}
\item There exist $\mu_1,\cdots,\mu_r\in\h^*$ such that $\V_{\mu}\neq 0$ implies
that $\mu < \mu_k$ for some $k=1,\cdots,r$.
\item For each $\mu\in\h^*$ such that $\V_{\mu}\neq 0$ and $i\in\bfI$, there
exists $N>0$ such that $\V_{\mu-n\alpha_i}=0$ for all $n\geq N$.
\end{itemize}

\noindent Normalization condition.
\begin{itemize}
\item $\Psi_i(\infty) = \Psi_i(0)^{-1} = q_i^{\alpha_i^{\vee}} =: K_i$, and
$\X_i^{\pm}(0)=0$.
\end{itemize}

\noindent Commutation relations.
\begin{itemize}
\item[(\QL1)] For any $i,j\in\bfI$, and $h,h'\in\h$,
\begin{gather*}
[\Psi_i(z),\Psi_j(w)]=0
\qquad
[\Psi_i(z),h]=0
\end{gather*}
\item[(\QL2)] For any $i\in\bfI$, and $h\in\h$,
\[[h,\X_i^\pm(z)]=\pm\alpha_i(h)\X_i^\pm(z)\]
\item[(\QL3)] For any $i,j\in\bfI$
\begin{multline*}
(z-q_i^{\pm a_{ij}}w)\Psi_i(z)\X_j^\pm(w)\\
=(q_i^{\pm a_{ij}}z-w)\X_j^\pm(w)\Psi_i(z)-
(q_i^{\pm a_{ij}} - q_i^{\mp a_{ij}})q_i^{\pm a_{ij}}w\X_j^\pm(q_i^{\mp a_{ij}}z)\Psi_i(z)
\end{multline*}
\item[(\QL4)] For any $i,j\in\bfI$
\begin{multline*}
(z-q_i^{\pm a_{ij}}w)\X_i^\pm(z)\X_j^\pm(w)-
(q_i^{\pm a_{ij}}z-w)\X_j^\pm(w)\X_i^\pm(z)\\
=z\lp\X_{i}^\pm(\infty)\X_j^\pm(w)-q_i^{\pm a_{ij}}\X_j^\pm(w)\X_{i}^\pm(\infty)\rp
+ w\lp\X_{j}^\pm(\infty)\X_i^\pm(z)-q_i^{\pm a_{ij}}\X_i^\pm(z)\X_{j}^\pm(\infty)\rp
\end{multline*}
\item[(\QL5)] For any $i,j\in\bfI$
\[(z-w)[\X^+_i(z),\X^-_j(w)] =
\frac{\delta_{ij}}{q_i-q_i^{-1}}\left(z\Psi_i(w)-w\Psi_i(z)-(z-w)\Psi_i(0)\right)\]
\end{itemize}
\end{defn}

\subsection{}\label{ssec: poles-lemma}

\begin{lem}\label{lem: poles-lemma}
Let $\V\in\Rloop$ and let $i\in\bfI$. Then the set of poles of $\Psi_i(z)$ is contained
in the set of poles of $\X_i^+(z)$ (resp. $\X_i^-(z)$).
\end{lem}
\begin{pf}
Consider the relation (\QL5) and its limit as $w\to\infty$. Then we get
\[
(q_i-q_i^{-1})[\X_i^+(z),\X_i^-(\infty)] = \Psi_i(z) - \Psi_i(0)
\]
which proves the $+$ case of the assertion. The $-$ case is proved similarly
by taking $z\to\infty$ limit of (\QL5).
\end{pf}

\subsection{Classification of irreducible representations}\label{ssec: qla-irr}

The following result was obtained in \cite{chari-pressley-qaffine-rep} for the case
when $\g$ is \fd simple Lie algebra
(see also \cite[Chaper 12]{chari-pressley}).
For the general case of symmetrisable Kac--Moody algebras, see \cite{hernandez-affinizations}.\\

Let $\gamma=\{\gamma_{i,\pm m}^{\pm}\}_{i\in\bfI, m\in \N}$ be a collection
of complex numbers and $\mu\in\h^*$ such that
$\gamma^{\pm}_{i,0} = q_i^{\pm\mu(\alpha_i^{\vee})}$.
A  representation $\V$ of $\qloop$ is said to be an $l$--\hw representation
of $l$--\hw $(\mu,\gamma)$ if there exists
$\fv\in \V$ such that
\begin{enumerate}
\item $\V = \qloop \fv$.
\item $\X^+_{i,k}\fv = 0$ for every $i\in\bfI$ and $k\in\Z$.
\item $\Psi_{i,\pm m}^{\pm}\fv = \gamma_{i,\pm m}^{\pm}\fv$ and $K_h\fv=
q^{\mu(h)}\fv$ for any $i\in \bfI,m\in\N$ and $h\in\h$.
\end{enumerate}

For any $(\mu,\gamma)$, there is a unique irreducible
representation with $l$--highest weight $(\mu,\gamma)$.

\begin{thm}\label{thm: qla-dp}\hfill
\begin{enumerate}
\item Every irreducible representation in $\Rloop$ is a \hw representation
for a unique highest weight $(\mu,\gamma)$.
\item The \irr representation $\V(\mu,\gamma)$ is in $\Rloop$ if, and
only if there exist monic polynomials $\{\PP_i(w)\in \C[w]\}_{i\in\bfI}$, $\PP
_i(0)\neq 0$, such that
\[\sum_{m\geq 0} \gamma_{i,m}^{+} z^{-m}=
q_i^{-\deg(\PP_i)} \frac{\PP_i(q_i^{2}z)}{\PP_i(z)}=
\sum_{m\leq 0} \gamma_{i,m}^{-} z^m\]
\end{enumerate}
\end{thm}

Note that we have $\mu(\alpha_i^{\vee}) = \deg(\PP_i)$ and hence $\mu$
is a dominant integral weight. The polynomials $\PP_i$ are called Drinfeld
polynomials. The set of isomorphism classes of simple objects in $\Rloop$ is in
bijection with the set $\dwtU$ of pairs $(\mu\in\h^*,\{\PP_i\})$ such that $\PP_i$
are monic, $\PP_i(0)\neq 0$ and $\mu(\alpha_i^{\vee})=\deg(\PP_i)$. We
denote by $\V(\mu,\{\PP_i\})$ the \irr $\qloop$--module corresponding to
$(\mu,\{\PP_i\})\in\dwtU$.

\section{An elliptic monodromy functor}\label{sec: functor}

The aim of this section is to construct a functor $\sfTh$ from
a dense subcategory of $\Rloop$
to $\eL$.

\subsection{Non--congruent representations}\label{ssec: rnc}

A representation $\V$ of $\qloop$ is said to be non--congruent if for any
$\alpha\neq\beta$ poles of $\X_i^+(z)$ (resp. $\X_i^-(z)$)
$\alpha\beta^{-1}\not\in p^{\Z}$.
Let $\Rloopnc$ be the full subcategory of $\Rloop$ consisting of non--congruent
representations.

\subsection{Contours}\label{ssec: contours}

By a Jordan curve $\cC$ we shall mean a disjoint union of simple closed curves
in $\C$, the inner domains of which are pairwise disjoint. For $\cC$ a Jordan
curve and $f$ a continuous function on $\cC$, we set
\[
\oint_{\cC} f(u)\, du = \frac{1}{2\pi\iota}\int_{\cC} f(u)\, du
\]

The definition of the functor $\sfTh$ relies upon the following choice of
contours of integration. For $\V\in\Rloop$, $\mu\in\h^*$ a weight of $\V$
and $i\in\bfI$, we choose a contour $\cC^{\pm}_{i,\mu}$ such that
\begin{itemize}
\item $\cC^{\pm}_{i,\mu}$ encloses a representative (modulo $\Z$) of
each pole of $\X_i^{\pm}(\Exp{u})$ acting on weight space $\V_{\mu}$.
\item $\cC^{\pm}_{i,\mu}$ does not enclose any $\nZ\tau$ translates of the
poles of $\X_i^{\pm}(\Exp{u})$ acting on weight space $\V_{\mu}$
and $\V_{\mu\pm\alpha_i}$.
\item The contour is small enough in the following sense.
Let $\mathcal{D}_{i,\mu}^{\pm}$ be the interior of the contour
$\cC^{\pm}_{i,\mu}$, including the contour itself. Then $\mathcal{D}^{\pm}_{i,\mu}
+\Lambda_{\tau} \subsetneq \C$.
\end{itemize}

\subsection{Main construction}\label{ssec: functor-defn}

Given $\V\in\Rloopnc$ we set $\sfTh(\V) = \V$ as a diagonalizable $\h$--module.
For each weight $\mu\in\h^*$ of $\V$ and $i\in\bfI$ we define operators
$\Phi_i(u), \eX_i^{\pm}(u,\lambda)$ on $\V_{\mu}$ as follows:
\begin{enumerate}
\item Define
\begin{equation}\label{eq:G+-}
G_i^{\pm}(z) := \prod_{n\geq 1} (K_i^{\pm 1}\Psi_i(p^{\pm n}z))
\end{equation}
Note that $G_i^{\pm}(z)$ satisfy the following multiplicative difference equations
\[
G_i^+(pz) = \lb K_i\Psi_i(pz)\rb^{-1}G_i^+(z) \aand 
G_i^-(pz) = \lb K_i^{-1}\Psi_i(z) \rb G_i^-(z)
\]
Since $K_i\Psi_i(0) = K_i^{-1}\Psi_i(\infty) = 1$, these equations are regular
near $z=0,\infty$ respectively. Therefore, it follows (see \eg \cite{sauloy})
that $G_i^{\pm}(z)$ are holomorphic in a neighborhood of $z=0,\infty$ respectively
and $G_i^+(0) = G_i^-(\infty)=1$.\\

\begin{equation}\label{eq:functorPhi}
\Phi_i(u) = \left.G_i^+(z)\Psi_i(z)G_i^-(z)\right|_{z=\Exp{u}}
\end{equation}

Hence, by construction $\Phi_i(u)$ is $1$--periodic and we have
\[
\Phi_i(u+\tau) = K_i^{-2}\Phi_i(u) = e^{-2\pi\iota\hbar d_i\alpha_i^{\vee}}\Phi_i(u)
\]

\item Define
\begin{equation}\label{eq:ellX}
\eX_i^{\pm}(u,\lambda) = c_i^{\pm} \oint_{\cC_{i,\mu}^{\pm}}
\frac{\theta(u-v+\lambda_i)}{\theta(u-v)\theta(\lambda_i)}
G_i^{\pm}(\Exp{v})\X_i^{\pm}(\Exp{v})\, dv
\end{equation}
where $\lambda_i = (\lambda,\alpha_i)$. And we consider the \rhs as defining a holomorphic
function outside of the shifts of the contour $\cC^{\pm}_{i,\mu}$ by elements of the
lattice $\Lambda_{\tau}$. Moreover, $c_i^{\pm}$ are constants which are to be chosen
to satisfy the following
\begin{equation}\label{eq:constant-fixing}
c_i^+c_i^- = (2\pi\iota)^2 \frac{\qG{0}^2}{\qG{d_i\hbar}^2}
\end{equation}
\end{enumerate}

\begin{rem}\label{rem: contour-independence}
The operators defined by \eqref{eq:ellX} are independent of the choice
of the contours $\cC_{i,\mu}^{\pm}$ satisfying the conditions of \S\ref{ssec: contours}.
This is because, by \eqref{eq:G+-} and Lemma \ref{lem: poles-lemma},
the poles of $G_i^{\pm}(z)$ are contained in $p^{\Z_{\lessgtr 0}}$--multiples of
the poles of $\X_i^{\pm}(z)$. In particular, $G_i^{\pm}(\Exp{u})$ are holomorphic
within $\cC_{i,\mu}^{\pm}$.
\end{rem}

\begin{thm}\label{thm:first-main-theorem}
The operators constructed above satisfy the axioms of \S \ref{ssec: categoryL}. Hence
we obtain a functor $\sfTh : \Rloopnc \to \eL$ which is faithful and exact.
\end{thm}

\subsection{}\label{ssec: prep-prop}

Let $i,j\in\bfI$ and set $a = d_ia_{ij}\hbar/2$. Consider a contour $\cC$
with interior domain $D$ and $\Omega_1,\Omega_2\subset\C$ two open subsets
with $\ol{D}\subset\Omega_2$. Assume given a holomorphic function $f(u,v):
\Omega_1\times\Omega_2\to\End(\V)$ such that $[\Psi_i(\Exp{u}),f(u,v)]=0$
for any $u,v$. We have the following analog of \cite[Proposition 5.5]
{sachin-valerio-2}.

\begin{prop}\label{prop: prep}
For each $\epsilon\in\{\pm 1\}$ we have:

\begin{enumerate}
\item If $u\not\in \ol{D}\pm \epsilon a+\Z$, then
\[
\Ad(\Psi_i\ee{u})^{\pm 1}\oint_{\cC} f(u,v)\X_j^{\epsilon}\ee{v}\, dv = 
\oint_{\cC} \left(\frac{e^{2\pi\iota(u+ \epsilon a)}-e^{2\pi\iota v}}
{e^{2\pi\iota u}-e^{2\pi\iota(v+\epsilon a)}}\right)^{\pm 1} f(u,v)\X_j^{\epsilon}\ee{v}\, dv
\]

\item For $u\not\in \ol{D} + \Z \pm \epsilon a + \Z_{<0}\tau$, we have
\[
\Ad(G^+_i\ee{u})^{\pm 1}\oint_{\cC} f(u,v)\X_j^{\epsilon}\ee{v}\, dv = 
\oint_{\cC} \lp  
\frac{\qG{u-v+\epsilon a}}{\qG{u-v-\epsilon a}}
\rp^{\pm 1} 
f(u,v)\X_j^{\epsilon}\ee{v}\, dv
\]

\item For $u\not\in \ol{D} + \Z \pm \epsilon a+ \Z_{>0}\tau$, we have
\[
\Ad(G^-_i\ee{u})^{\pm 1}\oint_{\cC} f(u,v)\X_j^{\epsilon}\ee{v}\, dv = 
\oint_{\cC} \lp  
\frac{\qGm{u-v+\epsilon a}}{\qGm{u-v-\epsilon a}}
\rp^{\pm 1} f(u,v)\X_j^{\epsilon}\ee{v}\, dv
\]
\item For $u\not\in \ol{D} \pm \epsilon a + \Lambda_{\tau}$, we have
\[
\Ad(\Phi_i(u))^{\pm 1}\oint_{\cC}f(u,v)\X_j^{\epsilon}\ee{v}\, dv = 
\oint_{\cC} \left(\frac{\theta(u-v+\epsilon a)}{\theta(u-v-\epsilon a)}\right)^{\pm 1}
f(u,v)\X_j^{\epsilon}\ee{v}\, dv
\]
\end{enumerate}
\end{prop}

\subsection{Proof of (\EQ3)}\label{ssec: pf-EQ3}

Let $i,j\in\bfI$ and let $a=d_ia_{ij}\hbar/2$. For $\lambda\in\hres$, we
write $\lambda_j = (\lambda,\alpha_j)$.\\

Using the definition of $\X_j^{\pm}(v,\lambda)$ given in \eqref{eq:ellX}
and Proposition \ref{prop: prep}, the left--hand side of (\EQ3) is given by
\[
c_j^{\pm}\oint_{\cC_j^{\pm}} \frac{\theta(u-v'\pm a)}{\theta(u-v'\mp a)}
\frac{\theta(v-v'+\lambda_j)}{\theta(v-v')\theta(\lambda_j)}
G_j^{\pm}\ee{v'}\X_j^{\pm}\ee{v'}\, dv'
\]
and by definition the right--hand side of (\EQ3) is given by a similar contour
integral
\[
c_j^{\pm}\oint_{\cC_j^{\pm}} \mathcal{K}(u,v,v') G_j^{\pm}\ee{v'}\X_j^{\pm}\ee{v'}\, dv'
\]
where
\[
\mathcal{K} = \frac{\theta(u-v\pm a)\theta(v-v'+\lambda_j\pm 2a)}
{\theta(u-v\mp a)\theta(\lambda_j\pm 2a)\theta(v-v')}+
\frac{\theta(\pm 2a)\theta(u-v\mp a-\lambda_j)\theta(u-v'+\lambda_j\pm a)}
{\theta(\lambda_j)\theta(u-v\mp a)\theta(\lambda_j\pm 2a)\theta(u-v'\mp a)}
\]
Thus we have to prove the following identity
\begin{multline*}
\frac{\theta(u-v'\pm a)}{\theta(u-v'\mp a)}
\frac{\theta(v-v'+\lambda_j)}{\theta(v-v')\theta(\lambda_j)}
=
\frac{\theta(u-v\pm a)\theta(v-v'+\lambda_j\pm 2a)}
{\theta(u-v\mp a)\theta(\lambda_j\pm 2a)\theta(v-v')}+\\
\frac{\theta(\pm 2a)\theta(u-v\mp a-\lambda_j)\theta(u-v'+\lambda_j\pm a)}
{\theta(\lambda_j)\theta(u-v\mp a)\theta(\lambda_j\pm 2a)\theta(u-v'\mp a)}
\end{multline*}
Clearing the denominator, the required identity takes the following form
(where for notational convenience, we write $b = \pm a$).
\begin{multline*}
\theta(u-v'+b)\theta(v-v'+\lambda_j)\theta(u-v-b)\theta(\lambda_j+2b) =\\ 
\theta(u-v+b)\theta(v-v'+\lambda_j+2b)\theta(u-v'-b)\theta(\lambda_j) + 
\theta(2b)\theta(u-v-\lambda_j-b)\theta(u-v'+\lambda_j+b)\theta(v-v')
\end{multline*}
which is precisely \eqref{eq: fti}.

\subsection{Proof of (\EQ4)}\label{ssec: pf-EQ4}

Let $i,j\in\bfI$ and let $\lambda\in\hres$.
Let us write $b = \pm a$ where $a=d_ia_{ij}\hbar/2$ and
$\lambda_k = (\lambda,\alpha_k)$ for $k=i,j$. With
these notations in mind, the \lhs of (\EQ4) can be written as
\[
\text{L.H.S.} = c_i^{\pm}c_j^{\pm} \oint_{\cC_j^{\pm}}\oint_{\cC_i^{\pm}} 
I(b) G_i^{\pm}\ee{u'}
\X_i^{\pm}\ee{u'}G_j^{\pm}\ee{v'}\X_j^{\pm}\ee{v'}\, du'\, dv'
\]
where
\begin{align*}
I(b) &= \frac{\theta(\lambda_i+\lambda_j)\theta(u-v-b)\theta(u-u'+\lambda_i+b)\theta(v-v'+\lambda_j-b)}
{\theta(u-u')\theta(\lambda_i+b)\theta(v-v')\theta(\lambda_j-b)} \\
&-\frac{\theta(u-v-\lambda_j)\theta(u-u'+\lambda_i+b)\theta(u-v'+\lambda_i+\lambda_j-b)}
{\theta(u-u')\theta(u+\lambda_i-v')\theta(\lambda_j-b)} \\
&-\frac{\theta(u-v+\lambda_i)\theta(v-u'+\lambda_i+\lambda_j+b)\theta(v-v'+\lambda_j-b)}
{\theta(v+\lambda_j-u')\theta(v-v')\theta(\lambda_i+b)}
\end{align*}
Similarly, the \rhs of (\EQ4) is equal to
\[
\text{R.H.S.} = c_i^{\pm}c_j^{\pm} \oint_{\cC_j^{\pm}}\oint_{\cC_i^{\pm}} 
I(-b) G_j^{\pm}\ee{v'}
\X_j^{\pm}\ee{v'}G_i^{\pm}\ee{u'}\X_i^{\pm}\ee{u'}\, du'\, dv'
\]

\noindent {\bf Claim:} We have the following identity. In particular,
$I(b)=0$ for $u'-v'-b\in\Lambda_{\tau}$.
\[
\theta(u'-v'+b)I(b) = \theta(u'-v'-b)I(-b)
\]

The proof of this claim is given in the next section. Assuming this, we can
finish the proof of (\EQ4) as follows. Using Proposition \ref{prop: prep}
we have
\begin{multline}\label{eq:EQ4-1}
\oint_{\cC} I(b)G_i^{\pm}\ee{u'}\X_i^{\pm}\ee{u'}\, du' \cdot G_j^{\pm}\ee{v'} = \\
G_j^{\pm}\ee{v'} \oint_{\cC} I(b) \frac{\qGpm{v'-u'-b}}{\qGpm{v'-u'+b}}
G_i^{\pm}\ee{u'}\X_i^{\pm}\ee{u'}\, du'
\end{multline}
as long as $v'$ does not lie in $\ol{D}+\Z \mp a \mp \Z_{>0}\tau$. However, using
the claim above, the integrand on the \rhs is regular on $\cC +\Z\mp a\mp \Z_{>0}\tau$.
This is because the simple poles of $1/\qGpm{v'-u'+b}$ are cancelled by the zeroes
of $I(b)$.\\

Thus multiplying \eqref{eq:EQ4-1} by $\X_j^{\pm}\ee{v'}$ and integrating $v'$
over the contour $\cC_j^{\pm}$, we obtain the following expression for the \lhs
of (\EQ4) divided by $c_i^{\pm}c_j^{\pm}$.

\begin{multline*}
\oint_{\cC_j^{\pm}}\oint_{\cC_i^{\pm}} I(b)\frac{\qGpm{v'-u'-b}}
{\qGpm{v'-u'+b}}G_i^{\pm}\ee{u'}G_j^{\pm}\ee{v'}\X_i^{\pm}\ee{u'}
\X_j^{\pm}\ee{v'}\, du'\, dv'
\end{multline*}

Similarly, the \rhs of (\EQ4) divided by $c_i^{\pm}c_j^{\pm}$ is equal to

\begin{multline*}
\oint_{\cC_j^{\pm}}\oint_{\cC_i^{\pm}} I(-b)\frac{\qGpm{u'-v'-b}}
{\qGpm{u'-v'+b}}G_i^{\pm}\ee{u'}G_j^{\pm}\ee{v'}\X_j^{\pm}\ee{v'}
\X_i^{\pm}\ee{u'}\, du'\, dv'
\end{multline*}

Using (\QL4) this expression can be rewritten as

\begin{multline*}
\oint_{\cC_j^{\pm}}\oint_{\cC_i^{\pm}} I(-b)\frac{\qGpm{u'-v'-b}}
{\qGpm{u'-v'+b}}\\
\cdot \frac{e^{2\pi\iota(v'+b)}-e^{2\pi\iota u'}}
{e^{2\pi\iota v'}-e^{2\pi\iota(u'+b)}}
G_i^{\pm}\ee{u'}G_j^{\pm}\ee{v'}\X_i^{\pm}\ee{u'}
\X_j^{\pm}\ee{v'}\, du'\, dv'
\end{multline*}

Thus, it suffices to check that
\[
I(b)\frac{\qGpm{v'-u'-b}}{\qGpm{v'-u'+b}} = 
I(-b)\frac{\qGpm{u'-v'-b}}{\qGpm{u'-v'+b}}
\frac{e^{2\pi\iota(v'+b)}-e^{2\pi\iota u'}}
{e^{2\pi\iota v'}-e^{2\pi\iota(u'+b)}}
\]

Using the claim, this equation is equivalent to the following

\[
\theta(u'-v'-b)\frac{\qGpm{v'-u'-b}}{\qGpm{v'-u'+b}} = 
\theta(u'-v'+b)\frac{\qGpm{u'-v'-b}}{\qGpm{u'-v'+b}}
\frac{e^{2\pi\iota(v'+b)}-e^{2\pi\iota u'}}
{e^{2\pi\iota v'}-e^{2\pi\iota(u'+b)}}
\]

which follows directly from \eqref{eq:thetaspm}.

\subsection{Proof of the claim}\label{ssec: pf-EQ4-claim}

Recall that $I(b) = T_1(b) - T_2(b) - T_3(b)$ where we have

\begin{align*}
T_1(b) &= \frac{\theta(\lambda_i+\lambda_j)\theta(u-v-b)\theta(u-u'+\lambda_i+b)\theta(v-v'+\lambda_j-b)}
{\theta(u-u')\theta(\lambda_i+b)\theta(v-v')\theta(\lambda_j-b)} \\
T_2(b) &= \frac{\theta(u-v-\lambda_j)\theta(u-u'+\lambda_i+b)\theta(u-v'+\lambda_i+\lambda_j-b)}
{\theta(u-u')\theta(u+\lambda_i-v')\theta(\lambda_j-b)} \\
T_3(b) &= \frac{\theta(u-v+\lambda_i)\theta(v-u'+\lambda_i+\lambda_j+b)\theta(v-v'+\lambda_j-b)}
{\theta(v+\lambda_j-u')\theta(v-v')\theta(\lambda_i+b)}
\end{align*}

Using \eqref{eq: fti} we can easily verify that

\begin{align*}
\theta(u'-v'+b)T_2(b)-\theta(u'-v'-b)T_2(-b) &=
\frac{\theta(2b)\theta(u-v-\lambda_j)\theta(u'-v'+\lambda_j)\theta(u-u'+\lambda_i+\lambda_j)}
{\theta(u-u')\theta(\lambda_j-b)\theta(\lambda_j+b)}\\
\theta(u'-v'+b)T_3(b)-\theta(u'-v'-b)T_3(-b) &=
-\frac{\theta(2b)\theta(u-v+\lambda_i)\theta(u'-v'-\lambda_i)\theta(v-v'+\lambda_i+\lambda_j)}
{\theta(v-v')\theta(\lambda_i-b)\theta(\lambda_i+b)}
\end{align*}

Using these, we get the following:
\begin{multline*}
\theta(u-u')\theta(v-v')\theta(\lambda_i-b)\theta(\lambda_i+b)\theta(\lambda_j-b)\theta(\lambda_j+b)
(\theta(u'-v'+b)I(b) - (u'-v'-b)I(-b)) = \\
\theta(u-v-b)\theta(u'-v'+b)\theta(u-u'+\lambda_i+b)\theta(v-v'+\lambda_j-b)
\theta(\lambda_i+\lambda_j)\theta(\lambda_i-b)\theta(\lambda_j+b)\\
-\theta(u-v+b)\theta(u'-v'-b)\theta(u-u'+\lambda_i-b)\theta(v-v'+\lambda_j+b)
\theta(\lambda_i+\lambda_j)\theta(\lambda_i+b)\theta(\lambda_j-b)\\
-\theta(u-v-\lambda_j)\theta(u'-v'+\lambda_j)\theta(u-u'+\lambda_i+\lambda_j)\theta(u-u')
\theta(\lambda_j+b)\theta(\lambda_j-b)\theta(2b)\\
+\theta(u-v+\lambda_i)\theta(u'-v'-\lambda_i)\theta(v-v'+\lambda_i+\lambda_j)\theta(v-v')
\theta(\lambda_i+b)\theta(\lambda_i-b)\theta(2b)
\end{multline*}

We need to prove that the \rhs of the equation written above is zero. The following
is a standard argument involving elliptic functions. Let us keep all the variables,
except for $u$, fixed and define:
\begin{multline*}
F_1(u) = \theta(u-v-b)\theta(u'-v'+b)\theta(u-u'+\lambda_i+b)\theta(v-v'+\lambda_j-b)
\theta(\lambda_i+\lambda_j)\theta(\lambda_i-b)\theta(\lambda_j+b)\\
-\theta(u-v-\lambda_j)\theta(u'-v'+\lambda_j)\theta(u-u'+\lambda_i+\lambda_j)\theta(u-u')
\theta(\lambda_j+b)\theta(\lambda_j-b)\theta(2b)\\
+\theta(u-v+\lambda_i)\theta(u'-v'-\lambda_i)\theta(v-v'+\lambda_i+\lambda_j)\theta(v-v')
\theta(\lambda_i+b)\theta(\lambda_i-b)\theta(2b)
\end{multline*}

\[
F_2(u) = \theta(u-v+b)\theta(u'-v'-b)\theta(u-u'+\lambda_i-b)\theta(v-v'+\lambda_j+b)
\theta(\lambda_i+\lambda_j)\theta(\lambda_i+b)\theta(\lambda_j-b)
\]

Note that both these functions have the same quasi--periodicity properties:
\[
F_s(u+1) = -F_s(u) \qquad \text{and}\qquad 
F_s(u+\tau) = -e^{-3\pi\tau}e^{-2\pi\iota(3u-v-u'-v'+2\lambda_i)}F_s(u)
\]

Also $F_2(u)$ has zeroes at $u=v-b$ and $u=u'-\lambda_i+b$ (modulo $\Lambda_{\tau}$).
The fact that $F_1(v-b) = 0$ and $F_1(u'-\lambda_i+b)=0$ is an easy 
application of \eqref{eq: fti}. Thus we deduce that $F_1/F_2$ is a holomorphic
doubly--periodic function of $u$ implying that it is a constant, say $C$ (constant
here means independent of $u$). To get that constant we use \eqref{eq: fti}
again to conclude that $F_1(v+b) = F_2(v+b)$, hence this constant must be $1$.
This proves that $F_1 = F_2$ which implies the desired claim.

\subsection{Proof of (\EQ5)}\label{ssec: pf-EQ5}

Let $i,j\in\bfI$ and $\lambda_1,\lambda_2\in\hres$ such that
$\lambda_1+\lambda_2 = \hbar(\mu+\alpha_i-\alpha_j)$. 
Let us write $c = (c_i^+c_j^-)^{-1}$, so that
\begin{align*}
c[\eX_i^+(u,\lambda_1),\eX_j^-(v,\lambda_2)] &= 
\oint_{\cC_j^-}\oint_{\cC_i^+} 
\mathcal{T}G_i^+\ee{u'}\X_i^+\ee{u'}G_j^-\ee{v'}\X_j^-\ee{v'}\, du'\, dv'\\
&-\oint_{\cC_j^-}\oint_{\cC_i^+} 
\mathcal{T}G_j^-\ee{v'}\X_j^-\ee{v'}G_i^+\ee{u'}\X_i^+\ee{u'}\, du'\, dv'
\end{align*}

where

\[
\mathcal{T} = \frac{\theta(u-u'+\lambda_{1,i})}
{\theta(u-u')\theta(\lambda_{1,i})} \frac{\theta(v-v'+\lambda_{2,j})}
{\theta(v-v')\theta(\lambda_{2,j})}
\]

We would like to apply Proposition \ref{prop: prep} in order to permute
the factors $\X_i^+$ and $G_j^-$ in the first term, $X_j^-$ and $G_i^+$
in the second. This cannot be done directly however, since $v'\in\cC_j^-$
may lie inside $\ol{D_i}^+-a+\Z+\Z_{>0}\tau$ (and similarly for $u'$). This
problem was encountered in \cite[\S 5.9]{sachin-valerio-2} and we employ
the same method as in there to circumvent this. Namely, let $\delta\in\C$
and define
\begin{align*}
J(\delta) &= \oint_{\cC_j^-}\oint_{\cC_i^+} 
\mathcal{T}G_i^+\ee{(u'+\delta)}\X_i^+\ee{u'}G_j^-\ee{(v'-\delta)}\X_j^-\ee{v'}\, du'\, dv'\\
&-\oint_{\cC_j^-}\oint_{\cC_i^+} 
\mathcal{T}G_j^-\ee{(v'-\delta)}\X_j^-\ee{v'}G_i^+\ee{(u'+\delta)}\X_i^+\ee{u'}\, du'\, dv'
\end{align*}

We consider this integral for $\delta$ in a disc $|\delta|<R$, where $R$ is such
that $G_i^+$ and $G_j^-$ are holomorphic in $\ol{D_i}^++\delta'$ and
$\ol{D_j}^--\delta'$ for any $\delta'$ with $|\delta'|<R$.\\

Moreover, if the contours $\cC_i^+$ and $\cC_j^-$ are small enough, there is
an $r<R$ such that if $|\delta|>r$, then $\ol{D_j}^--\delta$ is disjoint
from $\ol{D_i}^+-a+\Z+\Z_{>0}\tau$, and $\ol{D_i}^++\delta$ is disjoint
from $\ol{D_j}^-+a+\Z+\Z_{<0}\tau$.\\

Assuming $r<|\delta|<R$, we can apply Proposition \ref{prop: prep} to find that
$J(\delta)$ equals
\begin{gather*}
\oint_{\cC_j^-}\oint_{\cC_i^+} 
\mathcal{T} \frac{\qG{u'-v'+a+\delta}}{\qG{u'-v'-a+\delta}}G_i^+\ee{(u'+\delta)}
G_j^-\ee{(v'-\delta)} \X_i^+\ee{u'}\X_j^-\ee{v'}\, du'\, dv'\\
-\oint_{\cC_j^-}\oint_{\cC_i^+} 
\mathcal{T} \frac{\qG{u'-v'+a+\delta}}{\qG{u'-v'-a+\delta}}
G_j^-\ee{(v'-\delta)}G_i^+\ee{(u'+\delta)}
\X_j^-\ee{v'}\X_i^+\ee{u'}\, du'\, dv'
\end{gather*}
where we have used that $\qG{u} = \qGm{-u}$.
Using (\QL5) we get that the \rhs is zero for $i\neq j$. From now on we assume
that $i=j$ and hence $a = d_i\hbar$. Then we get

\begin{multline*}
(q_i-q_i^{-1})J(\delta) = \oint_{\cC_j^-}\oint_{\cC_i^+} 
\mathcal{T} \frac{\qG{u'-v'+d_i\hbar+\delta}}{\qG{u'-v'-d_\hbar+\delta}}
G_i^+\ee{(u'+\delta)}
G_i^-\ee{(v'-\delta)}\\
\cdot\frac{e^{2\pi\iota u'}\Psi_i\ee{v'}-e^{2\pi\iota v'}\Psi_i\ee{u'}}
{e^{2\pi\iota u'}-e^{2\pi\iota v'}}\, du'\, dv'
\end{multline*}

Writing
\begin{multline*}
\frac{e^{2\pi\iota u'}\Psi_i\ee{v'}-e^{2\pi\iota v'}\Psi_i\ee{u'}}
{e^{2\pi\iota u'}-e^{2\pi\iota v'}} = \\
\frac{e^{-2\pi\iota v'}\Psi_i\ee{v'}-e^{-2\pi\iota u'}\Psi_i\ee{u'}}
{u'-v'}.\frac{u'-v'}
{e^{-2\pi\iota v'}-e^{-2\pi\iota u'}}
\end{multline*}
and taking the two contours to be equal, say $\cC$,
we can apply \cite[Lemma 5.8]{sachin-valerio-2} to conclude that
\begin{equation}\label{eq:pf-EQ5}
(q_i-q_i^{-1})J(\delta) =\frac{1}{2\pi\iota} \oint_{\cC}
\ol{\mathcal{T}}\frac{\qG{\delta+d_i\hbar}}{\qG{\delta-d_i\hbar}}
G_i^+\ee{(u'+\delta)}
G_i^-\ee{(u'-\delta)}
\Psi_i\ee{u'}\, du'
\end{equation}
where $\ol{\mathcal{T}} = \left.\mathcal{T}\right|_{v'=u'}$, that is,
\[
\ol{\mathcal{T}} = \frac{\theta(u-u'+\lambda_{1,i})}
{\theta(u-u')\theta(\lambda_{1,i})} \frac{\theta(v-u'+\lambda_{2,j})}
{\theta(v-u')\theta(\lambda_{2,j})}
\]
Now both sides of \eqref{eq:pf-EQ5} are defined for $|\delta|<R$. Thus
we conclude that $J(0)$ is given by
\begin{align*}
(q_i-q_i^{-1})J(0) &= \frac{1}{2\pi\iota}\frac{\qG{d_i\hbar}}{\qGm{d_i\hbar}}
\oint_{\cC} \ol{\mathcal{T}}G_i^+\ee{u'}G_i^-\ee{u'}\Psi_i\ee{u'}
\, du'\\
&= \frac{1}{2\pi\iota}\frac{\qG{d_i\hbar}}{\qGm{d_i\hbar}}
\oint_{\cC} \ol{\mathcal{T}}\Phi_i(u')\, du'
\end{align*}
Thus, we have the following
\begin{align*}
\theta(d_i\hbar)[\eX_i^+(u,\lambda_1),\eX_i^-(v,\lambda_2)] &= \frac{c_i^+c_i^-}{2\pi\iota}
\frac{\qG{d_i\hbar}}{\qGm{d_i\hbar}}\frac{\theta(d_i\hbar)}{q_i-q_i^{-1}}
\oint_{\cC} \ol{\mathcal{T}}\Phi_i(u')\, du'\\
&= \frac{(2\pi\iota)^2}{2\pi\iota}\frac{\qG{0}^2}{\qG{d_i\hbar}^2}
\frac{\qG{d_i\hbar}}{\qGm{d_i\hbar}}\frac{\theta(d_i\hbar)}{q_i-q_i^{-1}}
\oint_{\cC} \ol{\mathcal{T}}\Phi_i(u')\, du'\\
&= \oint_{\cC} \ol{\mathcal{T}}\Phi_i(u')\, du'
\end{align*}
where we use the equation \eqref{eq:constant-fixing} in the second equation
and \eqref{eq:thetaspm} in the third.\\

In order to carry out the computation, we note the following identity
which follows from \eqref{eq: fti} upon clearing denominators:

\begin{multline*}
\frac{\theta(u-u'+x)}
{\theta(u-u')\theta(x)} \frac{\theta(v-u'-x+t)}
{\theta(v-u')\theta(-x+t)} = \\
\frac{\theta(u-v+x)}{\theta(u-v)\theta(x)}
\frac{\theta(v-u'+t)}{\theta(v-u')\theta(t)} - 
\frac{\theta(u-v+x-t)}{\theta(u-v)\theta(x-t)}
\frac{\theta(u-u'+t)}
{\theta(u-u')\theta(t)}
\end{multline*}

and letting $t\to 0$ in the equation above, we get

\begin{multline*}
\frac{\theta(u-u'+x)}
{\theta(u-u')\theta(x)} \frac{\theta(v-u'-x)}
{\theta(v-u')\theta(-x)} = \\
\frac{\theta(u-v+x)}{\theta(u-v)\theta(x)}
\frac{\theta'(v-u')}{\theta(v-u')} - 
\frac{\theta(u-v+x)}{\theta(u-v)\theta(x)}
\frac{\theta'(u-u')}
{\theta(u-u')} \\
-\ttheta{u-v+x}{u-v}{x}\frac{\theta'(x)}{\theta(x)}
+\frac{\theta'(u-v+x)}{\theta(u-v)\theta(x)}
\end{multline*}

Thus, using the fact that $\lambda_1+\lambda_2=\hbar\mu$, we obtain the following expression
for $\ol{\mathcal{T}}$:

\begin{itemize}
\item If $\mu_i\neq 0$ then
\begin{multline*}
\ol{\mathcal{T}} = 
\frac{\theta(u-v+\lambda_{1,i})}{\theta(u-v)\theta(\lambda_{1,i})}
\frac{\theta(v-u'+\hbar\mu_i)}{\theta(v-u')\theta(\hbar\mu_i)} + 
\frac{\theta(u-v-\lambda_{2,i})}{\theta(u-v)\theta(\lambda_{2,i})}
\frac{\theta(u-u'+\hbar\mu_i)}
{\theta(u-u')\theta(\hbar\mu_i)}
\end{multline*}

\item If $\mu_i=0$ then
\begin{multline*}
\ol{\mathcal{T}} = 
\frac{\theta(u-v+\lambda_{1,i})}{\theta(u-v)\theta(\lambda_{1,i})}
\frac{\theta'(v-u')}{\theta(v-u')} - 
\frac{\theta(u-v+\lambda_{1,i})}{\theta(u-v)\theta(\lambda_{1,i})}
\frac{\theta'(u-u')}
{\theta(u-u')} \\
-\ttheta{u-v+\lambda_{1,i}}{u-v}{\lambda_{1,i}}\frac{\theta'(\lambda_{1,i})}{\theta(\lambda_{1,i})}
+\frac{\theta'(u-v+\lambda_{1,i})}{\theta(u-v)\theta(\lambda_{1,i})}
\end{multline*}

\end{itemize}

We substitute these expressions for $\ol{\mathcal{T}}$ in 
\[
\theta(d_i\hbar)[\eX_i^+(u,\lambda_1),\eX_i^-(v,\lambda_2)] 
= \oint_{\cC} \ol{\mathcal{T}}\Phi_i(u')\, du'
\]

Note that the last two terms in the formula obtained for $\ol{\mathcal{T}}$,
for the case when $\mu_i=0$, will not contribute to the contour integral, since
in this case $\Phi_i(u)$ is a doubly--periodic function and hence the 
sum of its residues will be zero.\\

Finally the relation (\EQ5) follows from the following general result 
about doubly quasi--periodic functions.

\subsection{}
\begin{lem}\label{lem:complexanalysis}
Let $f(u)$ be a meromorphic function such that
\[
f(u+1) = f(u) \qquad\text{and}\qquad f(u+\tau) = e^{-2\pi\iota a}f(u)
\]
Choose $S\subset\C$ a (finite) set of representatives modulo $\Lambda_{\tau}$
of the poles of $f(u)$ and let $\cC$ be a contour enclosing $S$ and no other
poles of $f(u)$. 

\begin{itemize}
\item If $a=0$ then
\[
f(u) = \oint_{\cC}\frac{\theta'(u-u')}{\theta(u-u')}f(u')\, du' + K
\]
where $K$ is a contant.\\

\item If $a\not\in\Lambda_{\tau}$ then
\[
f(u) = \oint_{\cC}\ttheta{u-u'+a}{u-u'}{a} f(u')\, du'
\]
\end{itemize}
\end{lem}

\begin{pf}
This lemma is the analogue of partial fractions for doubly (quasi) periodic
functions. The proof is a standard exercise in complex analysis, see \eg
\cite[Section 21.50]{whittaker-watson}, and is given here solely for 
completeness.\\

For each $b\in S$ consider the Laurent series expansion of $f(u)$ near
$u=b$:
\[
f(u) = \sum_{n\in\N} \frac{f_{b,n}}{(u-b)^{n+1}} + \text{ regular part}
\]

Define a new function $\wt{f}(u)$ as:

\begin{itemize}
\item If $a=0$ then
\[
\wt{f}(u) = \sum_{\begin{subarray}{c} b\in S \\ n\in \N\end{subarray}}
f_{b,n} \frac{(-\partial_u)^n}{n!} \lp \frac{\theta'(u-b)}{\theta(u-b)}
\rp
\]
\item If $a\not\in\lambda_{\tau}$ then
\[
\wt{f}(u) = \sum_{\begin{subarray}{c} b\in S \\ n\in \N\end{subarray}}
f_{b,n} \frac{(-\partial_u)^n}{n!} \lp \ttheta{u-b+a}{u-b}{a}\rp
\]
\end{itemize}

Note that $\wt{f}(u)$ is the \rhs of the required equations. It remains to observe
that both $f(u)$ and $\wt{f}(u)$ have the same periodicity properties. Moreover
since both $\ds \frac{\theta'(x)}{\theta(x)}$ and $\ds\ttheta{x+a}{x}{a}$
are of the form $x^{-1} + O(1)$ near $x=0$ (here $O(1)$ stands for an 
element of $\C[[x]]$), we get that $f(u)-\wt{f}(u)$ is holomorphic. The
result now follows from the fact that every doubly periodic holomorphic
function has to be a constant, and the only double quasi--periodic holomorphic
function is zero.

\end{pf}

\section{Factorization problem}\label{sec: factorization}

\subsection{Set up and statement of the problem}\label{ssec: fp}

Let $V$ be a \fd vector space over $\C$ and $K\in\Aut(V)$ be an invertible operator.
Assume given a meromorphic function $\Phi : \nC \to \End(V)$ such that
\[
[K,\Phi(z)]=0 \aand \Phi(pz) = K^{-2}\Phi(z)
\]
where recall that $p = e^{2\pi\iota\tau}$ is a non--zero complex number with $|p|<1$.\\

\noindent {\bf Factorization problem.} Find $H^{\pm}(z)$ meromorphic $\End(V)$--valued
functions of a complex variable $z\in\nC$ such that
\begin{itemize}
\item[(F1)] $\ds \Phi(z) = H^+(z)^{-1}H^-(z)$.
\item[(F2)] $\ds [K,H^{\pm}(z)]=0$.
\item[(F3)] $H^{\pm}(z)$ are holomorphic and invertible in a neighborhood of $z=0,\infty$
respectively.
\item[(F4)] $H^-(\infty)=1$.
\end{itemize}

Note that we do not impose any normalization condition on $H^+(z)$.

\subsection{Coefficient matrix}\label{ssec: fact-coefficient}

Assuming $H^{\pm}$ is a solution to the factorization problem stated in section
\ref{ssec: fp}, we define $\ol{A}(z):= H^-(pz)H^-(z)^{-1}$. By the periodicity
condition on $\Phi(z)$ and (F1), we have
\[
\ol{A}(z) = K^{-2}H^+(pz)H^+(z)^{-1} = H^-(pz)H^-(z)^{-1}
\]

Since $H^{\pm}$ are regular near $z=0,\infty$ respectively, $\ol{A}(z)$ is a meromorphic
function on $\mathbb{P}^1$ and hence a rational function of $z$. Moreover, we have
\[
\ol{A}(0) = K^{-2} \aand \ol{A}(\infty) = 1
\]

Let $A(z) = K\ol{A}(z)$ so that $A(\infty) = K = A(0)^{-1}$. We refer to $A(z)$ as
the coefficient matrix.

\subsection{Isomonodromy transformation}\label{ssec: fact-iso}

If $H_1^{\pm}$ and $H_2^{\pm}$ are two solutions of the factorization problem
(section \ref{ssec: fp}), then we have
\[
H_2^+(z)H_1^+(z)^{-1} = H_2^-(z)H_1^-(z)^{-1}
\]
Since the \lhs of this equation is regular at $0$ and the \rhs is regular at $\infty$,
the resulting function is again a rational function. Let us denote it by
$G(z) = H_2^{\pm}(z)H_1^{\pm}(z)^{-1}$. This rational function is regular at
both $0$ and $\infty$ and $G(\infty)=1$.\\

Let $A_1(z)$ and $A_2(z)$ be the cofficient matrices defined using the
solutions $H_1^{\pm}$ and $H_2^{\pm}$ of the factorization problem, respectively.
Then we get
\begin{align*}
A_2(z) &= K H_2^-(pz)H_2^-(z)^{-1} = G(pz)KH_1^-(pz)H_1^-(z)^{-1}G(z)^{-1}\\
&= G(pz)A_1(z)G(z)^{-1}
\end{align*}

Given two subsets $S_1,S_2\in\nC$, we call the pair $(S_1,S_2)$
{\em non--congruent} if $\alpha_1\alpha_2^{-1}\not\in p^{\Z_{\neq 0}}$
for any $\alpha_{s}\in S_{s}$, $s=1,2$. We shall also say that $S\subset\nC$
is non--congruent if it is non--congruent to itself.
We have the following analog of \cite[Proposition 4.11]{sachin-valerio-2}.

\begin{prop}\label{pr:fact-unique}
Let $\mathcal{P}_s,\mathcal{Z}_s\subset\nC$ be the set of poles $A_s(z)$,
$A_s(z)^{-1}$ respectively ($s=1,2$). If the pairs
\[
(\mathcal{Z}_1,\mathcal{P}_1), \,
(\mathcal{Z}_2,\mathcal{P}_2), \,
(\mathcal{Z}_1,\mathcal{Z}_2), \,
(\mathcal{P}_1,\mathcal{P}_2)
\]
are non--congruent, then $A_1=A_2$. And therefore, $H_1^{\pm}=H_2^{\pm}$.
\end{prop}

\subsection{Existence of a factorization in the abelian case}\label{ssec: fact-existence}

Now we turn our attention to finding a solution to the factorization problem,
under the assumption that $[\Phi(z),\Phi(w)]=0$, for every $z,w\in\nC$.
Using Lemmas 4.12 and 4.13 of \cite{sachin-valerio-2}, we have the Jordan
decomposition of $\Phi(z) = \Phi_S(z)\Phi_U(z)$. Let $\calP, \calP_S, \calP_U$
be the sets of poles of $\Phi(z),\Phi_S(z),\Phi_U(z)$ respectively, and
$\calZ, \calZ_S, \calZ_U$ those of $\Phi(z)^{-1},\Phi_S(z)^{-1},\Phi_U(z)^{-1}$
respectively.
Then we have the following (see \cite[Lemma 4.13]{sachin-valerio-2})

\begin{enumerate}
\item $\calP = \calP_S\cup\calP_U$ and $\calZ = \calZ_S\cup\calZ_U$.
\item $\calP_U = \calZ_U$.
\end{enumerate}

We provide a solution to the factorization problem in the semisimple and unipotent
cases respectively in Sections \ref{ssec: fact-ss} and \ref{ssec: fact-uni} below.

\subsection{Semisimple case}\label{ssec: fact-ss}

Assuming both $K$ and $\Phi(z)$ are semisimple, we can restrict to their
joint eigenspace which reduces the problem to the scalar case. Let
$\eta\in\nC$ and $\varphi(z)$ a meromorphic function be a joint eigenvalue
of $K,\Phi(z)$. Then by general theory of elliptic functions we have
\[
\varphi(z) = C\left.\prod_l \frac{\theta(u-a_l)}{\theta(u-b_l)}\right|_{z=e^{2\pi\iota u}}
\]
where by the quasi--periodicity of $\Phi(z)$ we have
\begin{equation}\label{eq:consistency}
\eta = \exp\lp\pi\iota\sum_l (b_l-a_l)\rp
\end{equation}

Now using the expression \eqref{eq:thetaspm} for the theta function, we have
the following solution to the factorization problem:
\begin{align*}
H^+(z)^{-1} &= C\prod_l\lp \frac{z-\alpha_l}{z-\beta_l}\rp\frac{\theta^+(u-a_l)}
{\theta^+(u-b_l)} \\
H^-(z) &= \prod_l\frac{\theta^-(u-a_l)}{\theta^-(u-b_l)} 
\end{align*}
where the change of variable $\ds z=e^{2\pi\iota u}$ is understood.

\subsection{Unipotent case}\label{ssec: fact-uni}

Now we assume that both $\Phi(z)$ and $K$ are unipotent. This enables us to
take logarithm in order to convert the problem to an additive one. Namely,
we have the following periodicity property 
\begin{equation}\label{eq:fact-uni-add}
\log(\Phi(pz)) - \log(\Phi(z)) = \log(K^{-2})
\end{equation}
and we are required to find $h^{\pm}(z)$, regular at
$0$ and $\infty$ respectively, such that $h^-(\infty)=0$ and
\[
\log(\Phi(z)) = h^-(z) - h^+(z)
\]

The problem again reduces to the one for a single function. Let $\varphi(z)$
be an entry of the matrix $\log(\Phi(z))$ and $k\in\C$ be the correspoding
entry of $\log(K^{-2})$, so that we have $\varphi(pz)-\varphi(z) = k$.
Again, from the general theory of elliptic functions, we have the following
expression for $\varphi(z)$:
\[
\varphi(z) = \sum_{\begin{subarray}{c} a\in \PP \\ n\in\N \end{subarray}}
f_{a,n} \frac{\partial_u^{n+1}}{(n+1)!} \log(\theta(u-a))
\]
where $\PP\subset\C$ is a choice modulo $\Lambda_{\tau}$ of representatives
of poles of $\varphi\ee{u}$, and $-2\pi\iota \sum_a f_{a,0} = k$.\\

Using the equation \eqref{eq:thetaspm} we obtain the following solution
\begin{align*}
h^+(z) &= -\sum_{\begin{subarray}{c} a\in \PP \\ n\in\N \end{subarray}}
f_{a,n} \frac{\partial_u^{n+1}}{(n+1)!} 
\log((e^{\pi\iota u}-e^{-\pi\iota u})\theta^+(u-a))\\
h^-(z) &= \phantom{-}\sum_{\begin{subarray}{c} a\in \PP \\ n\in\N \end{subarray}}
f_{a,n} \frac{\partial_u^{n+1}}{(n+1)!} \log(\theta^-(u-a))
\end{align*}

Note that $h^+(0) = -\pi\iota\sum_{a\in\PP} f_{a,0} = k/2$.

\subsection{}\label{ssec:fact-cor}

Let us define $G^-(z) = H^-(z)$ and $G^+(z) = H^+(z)^{-1}A(z)^{-1}$, where $A(z)$
was defined in Section \ref{ssec: fact-coefficient}, so that
$\Phi(z) = G^+(z)A(z)G^-(z)$. As a consequence of the explicit factorization
given in the previous two sections we have the following:

\begin{cor}\label{cor:fact-cor}
$G^+(0)$ is semisimple. More precisely, if $V'\subset V$ is a generalized
eigenspace for $\Phi(z)$ with the following generalized eigenvalue
\[
\varphi(z) = C\prod_l \frac{\theta(u-a_l)}{\theta(u-b_l)}
\]
Then for the solution given in Section \ref{ssec: fact-ss}, $G^+(0) = C$
on $V'$.
\end{cor}

\section{The inverse functor}\label{sec: functor-surjective}

\subsection{}\label{ssec: choice-slice}

From now on we fix a subset $\slice\subset\C$ subject to the following conditions
\begin{itemize}
\item $a\neq b\in \slice \Rightarrow a-b\not\in\Lambda_{\tau}$.
\item $\slice\pm \frac{\hbar}{2} = \slice$.
\item $\slice$ is in bijection with $\C/\Lambda_{\tau}$ under the natural surjection
$\C\to\C/\Lambda_{\tau}$.
\end{itemize}

\begin{rem}
The existence of such a subset $\slice$ relies upon Zorn's lemma
(see \cite[\S 6.2]{sachin-valerio-2}). Note that the argument involving
Zorn's lemma is only needed to ensure that the last condition on $\slice$,
imposed above, holds. Theorem \ref{thm:inverse} and
Corollary \ref{cor:isocats} will continue to hold, even if we drop
the last condition on $\slice$. However, Corollary \ref{cor:secondgauge}
will no longer be true. Thus, if we were to avoid the use of Zorn's lemma,
we will still be able to prove the main classification theorem 
(Theorem \ref{thm: class-final}) by appealing to certain finiteness properties.
For instance,
\begin{itemize}
\item finite--dimensionality of $\eV$, in case $\g$ is of finite type.
\item finitely many $\hbar$--strings of poles (as proved in Section 
\ref{sec: irr-class1}), in case $\g$ is arbitrary and $\eV$ is irreducible.
\end{itemize}
\end{rem}

Let us assume that $\eV\in\eL$ is a \lflat object (see Section \ref{ssec: lflat}). 
We consider the factorization problem
for the functions $\Phi_i(u)\in\End(\eV)$. Using Proposition \ref{pr:fact-unique},
Section \ref{ssec: fact-existence} and Proposition \ref{prop: eqg-knight} we get

\begin{prop}\label{prop:flat-factorization}
There exist unique functions $G_i^{\pm}(z), \Psi_i(z)\in\End(\eV)$ such that

\begin{enumerate}
\item $\Phi_i(u) = G_i^+(z)\Psi_i(z)G_i^-(z)$ where $z = e^{2\pi\iota u}$.\\

\item $G_i^{\pm}(z)$ are holomorphic and invertible near $z=0, \infty$ respectively.
$\Psi_i(z)$ is a rational function of $z$ regular at both $z=0$ and $\infty$.\\

\item $G_i^-(\infty)=1$.\\

\item $$\sigma(G_i^{\pm}(z))\subset \bigcup_{n\geq 1} p^{\mp n}\slice^*
\text{ and } \sigma(\Psi_i(z))\subset\slice^*$$
where $\slice^* = \{\exp(2\pi\iota s)\}_{s\in\slice}$.
\end{enumerate}
\end{prop}

\subsection{}

Combining the results Proposition \ref{prop:flat-factorization}, \ref{prop: eqg-knight}
and Corollary \ref{cor:fact-cor} we get that for every generalized eigenspace
$\eV[\ulB]$ the functions $B_i(u)$ have the following form:

\begin{equation}\label{eq:almost-right}
B_i(u) = C_i(\ulB)\prod_k \frac{\theta(u-c_{i,k}+d_i\hbar)}{\theta(u-c_{i,k})}
\prod_l \frac{\theta(u-c'_{i,l}-d_i\hbar)}{\theta(u-c'_{i,l})}
\end{equation}
where the numbers $c_{i,k}, c'_{i,l}\in \slice$.\\

And the matrix $G_i^+(0)$ is given by:

\begin{equation}\label{eq:g+0-exp}
G_i^+(0) = \sum_{\ulB} C_i(\ulB) \Id_{\eV[\ulB]}
\end{equation}

\subsection{Second gauge transformation}

Consider the following zero--weight $\End(\eV)$--valued function, which is 
given by the following expression on a weight space $\eV_{\mu}$
\[
\varphi(\lambda) = \exp\lp\frac{1}{\hbar}\sum_{j\in\bfI}
\lp \lambda+\frac{\hbar}{2}\mu,\varpi_j^{\vee}\rp\ln(G_j^+(0))\rp
\]

We conjugate $\eV$ by the automorphism given above to obtain $\eV^{\varphi}$ which is manifesly
isomorphic to $\eV$ with isomorphism given by $\varphi(\lambda)$.
Again $\eV^{\varphi} = \eV$ as $\h$--diagonalizable module and the action 
of $\{\Phi_i(u),\eX^{\pm}_i(u,\lambda)\}$ on a weight space $\eV^{\varphi}_{\mu}$
is given by:

\begin{align*}
\Phi_i^{\varphi}(u) &=
\varphi\lp\lambda+\frac{\hbar}{2}\alpha_i\rp^{-1} \Phi_i(u)
\varphi\lp\lambda-\frac{\hbar}{2}\alpha_i\rp\\
\eX^{\pm,\varphi}_i(u,\lambda) &=
\varphi\lp \pm\lambda \mp \frac{\hbar}{2}(\mu\pm\alpha_i) + \frac{\hbar}{2}\alpha_i\rp^{-1}
\eX^{\pm}_i(u,\lambda)\\
& \qquad.\varphi\lp \pm\lambda \mp \frac{\hbar}{2}\mu - \frac{\hbar}{2}\alpha_i\rp
\end{align*}

\begin{prop}\label{prop:lambda-independence}
$\eV^{\varphi}\in\eL$ and we have the following
\begin{enumerate}
\item The generalized eigenvalues of $\Phi_i^{\varphi}(u)$ have the following form
\[
\prod_k \frac{\theta(u-c_{i,k}+d_i\hbar)}{\theta(u-c_{i,k})}
\prod_l \frac{\theta(u-c'_{i,l}-d_i\hbar)}{\theta(u-c'_{i,l})}
\]
where the numbers $c_{i,k}, c'_{i,l}\in \slice$.\\

\item $G_i^{+,\varphi}(0) = 1 = G^{-,\varphi}(\infty)$.\\

\item For each $i\in\bfI$, $b\in\slice$ and $n\in\N$, $\partialX^{\pm,\varphi}_{i;b,n}(\lambda)$
is independent of $\lambda$.
\end{enumerate}
\end{prop}

\begin{pf}
We begin by computing the action of $\Phi_i^{\varphi}(u)$. Consider a generalized eigenspace
$\eV[\ulB]$ where $\ulB$ consists of $\bfI$--tuple of functions $(B_i(u))$ given
by \eqref{eq:almost-right}. On this generalized eigenspace $\Phi_i^{\varphi}(u)$
acts by:
\begin{align*}
\left.\Phi_i^{\varphi}(u)\right|_{\eV[\ulB]} &=
\exp\lp\frac{1}{\hbar}\sum_{j\in\bfI} (-\hbar\alpha_i,\varpi_j^{\vee})\ln(C_i(\ulB))\rp
\left.\Phi_i(u)\right|_{\eV[\ulB]} \\
&= C_i(\ulB)^{-1}\left.\Phi_i(u)\right|_{\eV[\ulB]}
\end{align*}

This implies that the shifted conjugation by $\varphi$ of $\Phi_i(u)$ is independent
of $\lambda$, as we require. Moreover we obtain (1) and (2) of the Proposition.\\

Let us consider the raising and lowering operators now. Using Proposition
\ref{prop: important-prop} we know that between two generalized eigenspaces
$\eV_\mu[\ulB]$ and $\eV_{\mu\pm\alpha_j}[\ulB^{\pm}]$ these operators are of the form:
\[
\partialX^{\pm}_{j;b,n}(\lambda)
 = e^{\pm 2\pi\iota (\lambda,\alpha^{\pm})}\partialX^{\pm}_{j;b,n}(0)
\]
and in this case the constants from \eqref{eq:almost-right} are related by
\[
C_i(\ulB^{\pm}) = C_i(\ulB)e^{2\pi\iota\hbar (\alpha_i,\alpha^{\pm})}
\]

With this information in mind, we can compute the shifted conjugation.
\begin{align*}
\partialX^{\pm,\varphi}_{j;b,n}(\lambda) &=
C^{\pm}\exp\lp\frac{1}{\hbar}\sum_{k\in\bfI} 2\pi\iota\hbar (\alpha^{\pm},\alpha_k)
(\mp\lambda,\varpi_k^{\vee})\rp e^{\pm 2\pi\iota (\lambda,\alpha^{\pm})}
\partialX^{\pm}_{j;b,n}(0) \\
&= C^{\pm} e^{\mp 2\pi\iota (\lambda,\alpha^{\pm})}e^{\pm 2\pi\iota (\lambda,\alpha^{\pm})}
\partialX^{\pm}_{j;b,n}(0) \\
&= C^{\pm}\partialX^{\pm}_{j;b,n}(0)
\end{align*}
where we have used the fact that $\lambda\in\ol{\h^*} = \sum_{i\in\bfI}\C\alpha_i$
to conclude that 
\[
\sum_{k\in\bfI} (\lambda,\varpi_k^{\vee})(\alpha,\alpha_k) = (\lambda,\alpha)
\]
for every $\alpha\in Q$. The constants $C^{\pm}$ appearing above are independent of $\lambda$.
They can be easily computed from the definition:

\begin{align*}
C^+ &= C_j(\ulB)^{-1}\exp\lp -\pi\iota\hbar (\alpha_j,\alpha^+)\rp\\
C^- &= C_j(\ulB)^{-2}\exp\lp -3\pi\iota\hbar (\alpha_j,\alpha^-)\rp
\exp\lp -2\pi\iota\hbar \sum_{k\in\bfI} (\mu,\varpi_k^{\vee})(\alpha^-,\alpha_k)\rp
\end{align*}

This calculation together with the partial fraction expression of the
raising and lowering operators from Lemma \ref{lem: xpm-partial}
implies both (3) and the fact that the resulting representation 
$\eV^{\varphi}$ is in $\eL$.
\end{pf}

\subsection{Definition of category $\eLskeleton$}\label{ssec: categoryskeleton}

Let $\eLskeleton$ be the full subcategory of $\eL$ consisting of $\eV$ which satisfies (1) 
of Proposition \ref{prop:lambda-independence}.
That is, if for each $i\in \bfI$, the generalized eigenvalues of $\Phi_i(u)$ acting on $\eV$
are of the form:
\[
\prod_k \frac{\theta(u-c_{i,k}+d_i\hbar)}{\theta(u-c_{i,k})}
\prod_l \frac{\theta(u-c_{i,l}'+d_i\hbar)}{\theta(u-c_{i,l}')}
\]
where the numbers $c_{i,k}, c_{i,l}'$ are in $\slice$. For future reference, we record
the following consequence of Proposition \ref{prop:lambda-independence}

\begin{cor}\label{cor:secondgauge}
Every object of $\eL$ is isomorphic to some object of $\eLskeleton$.
\end{cor}

Later we will also need a $\C$--linear analogue of $\eLskeleton$ which we define
as follows. Let $\eLClinear$ be the category whose objects are same as those of
$\eLskeleton$ while its morphisms are the usual $\C$--linear homomorphisms commuting
with the action of $\h$ and $\{\Phi_i(u),\eX_i^{\pm}(u,\lambda)\}$.\\

Following \cite[Section 3.5]{sachin-valerio-2} we consider the full category $\Rloopslice$
of $\Rloop$ consisting of $\V$ for which the poles of $\{\Psi_i(z)^{\pm 1}\}$ lie in 
$\exp(2\pi\iota \slice)$. We refer the reader to \cite[Theorem 3.8]{sachin-valerio-2}
for several equivalent characterizations of $\Rloopslice$.

\subsection{Choice of contours}\label{ssec:inv-contours}

Since $\eV$ and $\eV^{\varphi}$ are isomorphic, we identify the two and thus omit the
superscript $\varphi$.
In the next section we define an action of $\qloop$ on $\eV$ which relies upon a choice
of contours. Let $\mu\in\h^*$ be a weight of $\eV$ and $i\in\bfI$. Let $\cC_{i,\mu}^{\pm}$
be a Jordan curve such that

\begin{itemize}
\item $\cC_{i,\mu}^{\pm}$ encloses all the poles of $\eX_i^{\pm}(u,\lambda)$ in the spectral
variable $u$, which are contained in $\slice$, and none of their translates under $\Lambda_{\tau}
\setminus\{0\}$. \\

\item $\cC_{i,\mu}^{\pm}$ does not enclose any $\Lambda_{\tau}\setminus\{0\}$ translates of
the poles of $\Phi_i(u)_{\mu}^{-1}$.
\end{itemize}

Note that such a curve exists since the poles of $\Phi_i(u)^{-1}$ are contained in
$\slice$ by Proposition \ref{prop: eqg-knight}, and the latter is non--congruent.

\subsection{Inverse construction}

Let $\Psi_i(z), G_i^{\pm}(z)$ be the $\End(\eV)$--valued
functions given by Proposition \ref{prop:flat-factorization}.
Let us write $\V = \sfhT(\eV) = \eV$
as $\h$--diagonalizable module and define operators $\Psi_i(z), \X^{\pm}_i(z)$
on the weight space $\V_{\mu}$ as:\\

\begin{enumerate}
\item $\Psi_i(z)$ as given by Proposition \ref{prop:flat-factorization}.\\

\item For each $i\in\bfI$
\begin{equation}\label{eq:inv-Xpm}
\X_i^{\pm}(z) := \frac{1}{c_i^{\pm}}\int_{\cC^{\pm}_{\mu}}
\frac{z}{z-e^{2\pi\iota u}} G_i^{\pm}\ee{u}^{-1}\eX^{\pm}_i(u,\lambda)\, du
\end{equation}
where the contours $\cC^{\pm}_{i,\mu}$ are chosen as in section
\ref{ssec:inv-contours} and the constants $c_i^{\pm}$ are as in \eqref{eq:constant-fixing}.
Note that by (3) of Proposition \ref{prop:lambda-independence}
the resulting operators are independent of $\lambda$.
\end{enumerate}

\begin{thm}\label{thm:inverse}
The operators $\Psi_i(z),\X^{\pm}_i(z)$ defined above satisfy the relations
of $\qloop$ given in Section \ref{ssec: rloop} and hence define a representation
$\sfhT(\eV)\in\Rloopslice$. Moreover $\sfTh(\sfhT(\eV)) = \eV$ and for any
$\V\in\Rloopslice$ we have $\sfhT(\sfTh(\V)) = \V$. In particular, the categories
$\eLskeleton$ and $\Rloopslice$ are isomorphic.
\end{thm}

\begin{pf}
We need to verify the axioms of Section \ref{defn: rloop}. Note that the category $\mathcal{O}$
and integrability condition holds by the corresponding assumption on $\eV$. The normalization
condition follows from Section \ref{ssec: fact-coefficient} and the definition of $\X^{\pm}_i(z)$.\\

The relations (\QL 1) and (\QL 2) are immediate. After working out the commutation relations
between $G_i^{\pm}(z)$ and operators defined using contour integration involvoing $\eX_i^{\pm}(u,\lambda)$
in Section \ref{ssec:inv-pf-comm}, we prove the relations (\QL 3), (\QL 4) and (\QL 5)
in Sections \ref{ssec:inv-pf-ql3}, \ref{ssec:inv-pf-ql4} and \ref{ssec:inv-pf-ql5}
respectively.\\

The fact that the commuting operators on $\eV$ and $\V$ are the same follows
from the uniqueness statement of Proposition \ref{pr:fact-unique}.
The proof for the raising and lowering operators is given in Sections \ref{ssec:inv-pf-c1}
and \ref{ssec:inv-pf-c2}.
\end{pf}

Restricting to the category $\eLClinear$ defined in Section \ref{ssec: categoryskeleton},
we obtain the following.

\begin{cor}\label{cor:isocats}
The functors $\sfTh^{\slice} : \Rloopslice \to \eLClinear$ and
$\sfhT^{\slice} : \eLClinear \to \Rloopslice$ are inverse to each other and hence
the two categories are isomorphic.
\end{cor}

\subsection{Commutation relations}\label{ssec:inv-pf-comm}

Now let $i,j\in\bfI$ and $a = d_ia_{ij}\hbar/2$. Consider a contour $\cC$
with interior domain $D$ and $\Omega_1,\Omega_2\subset\C$ two open subsets
with $\ol{D}\subset\Omega_2$. Assume given a holomorphic function $f(u,v):
\Omega_1\times\Omega_2\to\End(\eV)$ such that $[\Phi_i(u),f(u,v)]=0$
for any $u,v$.

\begin{prop}\label{prop: eqg-prep}
For each $\epsilon\in\{\pm 1\}$ we have:

\begin{enumerate}
\item If $u\not\in \ol{D}\pm \epsilon a+\Z$, then
\[
\Ad(\Psi_i\ee{u})^{\pm 1}\oint_{\cC} f(u,v)\eX_j^{\epsilon}(v,\lambda)\, dv = 
\oint_{\cC} \left(\frac{e^{2\pi\iota(u+ \epsilon a)}-e^{2\pi\iota v}}
{e^{2\pi\iota u}-e^{2\pi\iota(v+\epsilon a)}}\right)^{\pm 1} f(u,v)\eX_j^{\epsilon}(v,\lambda)\, dv
\]

\item For $u\not\in \ol{D} + \Z \pm \epsilon a + \Z_{<0}\tau$, we have
\[
\Ad(G^+_i\ee{u})^{\pm 1}\oint_{\cC} f(u,v)\eX_j^{\epsilon}(v,\lambda)\, dv = 
\oint_{\cC} \lp  
\frac{\qG{u-v+\epsilon a}}{\qG{u-v-\epsilon a}}
\rp^{\pm 1} 
f(u,v)\eX_j^{\epsilon}(v,\lambda)\, dv
\]

\item For $u\not\in \ol{D} + \Z \pm \epsilon a+ \Z_{>0}\tau$, we have
\[
\Ad(G^-_i\ee{u})^{\pm 1}\oint_{\cC} f(u,v)\eX_j^{\epsilon}(v,\lambda)\, dv = 
\oint_{\cC} \lp  
\frac{\qGm{u-v+\epsilon a}}{\qGm{u-v-\epsilon a}}
\rp^{\pm 1} f(u,v)\eX_j^{\epsilon}(v,\lambda)\, dv
\]
\end{enumerate}
\end{prop}

\subsection{Proof of (\QL 3)}\label{ssec:inv-pf-ql3}

Let $i,j\in\bfI$ and set $a = \hbar d_ia_{ij}/2$. We need to prove that
\begin{multline*}
\Ad(\Psi_i(z))\X_j^{\pm}(w) = \frac{q_i^{\pm a_{ij}}z-w}{z-q_i^{\pm a_{ij}}w}\X_j^{\pm}(w)
 - \frac{(q_i^{\pm a_{ij}}-q^{\mp a_{ij}})q^{\pm a_{ij}}w}{z-q^{\pm a_{ij}} w}
\X_j^{\pm}(q_i^{\mp a_{ij}}z)
\end{multline*}

Using the definition \eqref{eq:inv-Xpm} and (1) of Proposition \ref{prop: eqg-prep} the \lhs
of this equation becomes
\[
\text{L.H.S. } = \frac{1}{c_j^{\pm}}\int_{\cC_j}
\frac{q_i^{\pm a_{ij}} z - w'}{z-q_i^{\pm a_{ij}}w'} \frac{w}{w-w'}
G_j^{\pm}\ee{u}^{-1}\eX_j^{\pm}(u,\lambda)\, du
\]
where for notational convenience we wrote $w' = e^{2\pi\iota u}$.\\

Similarly the \rhs takes the following form
\[
\text{R.H.S. } = \frac{1}{c_j^{\pm}}\int_{\cC_j} \mathcal{K}(z,w,u)
G_j^{\pm}\ee{u}^{-1}\eX_j^{\pm}(u,\lambda)\, du
\]

where
\[
\mathcal{K}(z,w,u) = \frac{q_i^{\pm a_{ij}} z - w}{z-q_i^{\pm a_{ij}}w} \frac{w}{w-w'}
- \frac{(q_i^{\pm a_{ij}}-q^{\mp a_{ij}})q^{\pm a_{ij}}w}{z-q^{\pm a_{ij}} w}
\frac{z}{z-q_i^{\pm a_{ij}}w'}
\]

Thus the relation (\QL 3) reduces to the following algebraic identity which can
be checked easily:
\[
\mathcal{K}(z,w,u) = \frac{q_i^{\pm a_{ij}} z - w'}{z-q_i^{\pm a_{ij}}w'} \frac{w}{w-w'}
\]

\subsection{Proof of (\QL 4)}\label{ssec:inv-pf-ql4}

For $i,j\in\bfI$, let $a=\hbar d_ia_{ij}/2$. We need to prove the following relation:

\begin{multline*}
(z-q_i^{\pm a_{ij}}w)\X_i^\pm(z)\X_j^\pm(w)
-z\X_{i}^\pm(\infty)\X_j^\pm(w)
+q_i^{\pm a_{ij}}w\X_i^\pm(z)\X_{j}^\pm(\infty)
=\\
(q_i^{\pm a_{ij}}z-w)\X_j^\pm(w)\X_i^\pm(z)
-q_i^{\pm a_{ij}}z\X_j^\pm(w)\X_{i}^\pm(\infty)
+ w\X_{j}^\pm(\infty)\X_i^\pm(z)
\end{multline*}

Let us take $\lambda\in\hres$ and let
\[
\lambda_1 = \lambda\pm\frac{\hbar}{2}\alpha_j \text{ and }
\lambda_2 = \lambda\mp\frac{\hbar}{2}\alpha_i
\]

Again using (2) and (3) of Proposition \ref{prop: eqg-prep} and the definition
\eqref{eq:inv-Xpm} we find that the \lhs multiplied by $c_i^{\pm}c_j^{\pm}$ is given by
\[
\iint \mathcal{I}_1(z,w,u,v)
G_i^{\pm}\ee{u}^{-1}G_j^{\pm}\ee{v}^{-1}\eX_i^{\pm}(u,\lambda_1)
\eX_j^{\pm}(v,\lambda_2)\, du\, dv
\]
where (we write $z'=e^{2\pi\iota u}$ and $w'=e^{2\pi\iota v}$):
\[
\mathcal{I}_1 = \frac{zw(z'-q_i^{\pm a_{ij}}w')}{(z-z')(w-w')}\frac{\theta^{\pm}(v-u\pm a)}
{\theta^{\pm}(v-u\mp a)}
\]

Similarly the \rhs multiplied by $c_i^{\pm}c_j^{\pm}$ equals
\[
\iint \mathcal{I}_2(z,w,u,v)
G_i^{\pm}\ee{u}^{-1}G_j^{\pm}\ee{v}^{-1}\eX_j^{\pm}(v,\lambda_2)
\eX_i^{\pm}(u,\lambda_1)\, du\, dv
\]
where
\[
\mathcal{I}_2 = \frac{zw(q_i^{\pm a_{ij}}z'-w')}{(z-z')(w-w')}\frac{\theta^{\pm}(u-v\pm a)}
{\theta^{\pm}(u-v\mp a)}
\]

Using relation (\EQ 4) we can flip the $\eX_j$ and $\eX_i$ factors to get
that the \rhs (multiplied by $c_i^{\pm}c_j^{\pm}$) is equal to
\[
\iint \mathcal{I}_2(z,w,u,v)\frac{\theta(u-v\mp a)}{\theta(u-v\pm a)}
G_i^{\pm}\ee{u}^{-1}G_j^{\pm}\ee{v}^{-1}\eX_i^{\pm}(u,\lambda_1)
\eX_j^{\pm}(v,\lambda_2)\, du\, dv
\]

Hence, in order to prove the relation (\QL 4), we need to check the following identity
\[
\frac{\theta(u-v+b)}{\theta(u-v-b)} = \frac{e^{2\pi\iota u} - e^{2\pi\iota (v+b)}}
{e^{2\pi\iota (u+b)} - e^{2\pi\iota v}}
\frac{\theta^+(u-v+b)}{\theta^+(u-v-b)}\frac{\theta^-(u-v+b)}{\theta^-(u-v-b)}
\]
which follows directly from \eqref{eq:thetaspm}.\\

\subsection{Proof of (\QL 5)}\label{ssec:inv-pf-ql5}

We will prove the following relation, for each $i,j\in\bfI$ and $k,l\in\Z$:
\[
[\X^+_{i,k}, \X^-_{j,l}] = \delta_{ij}
\frac{\Psi_{i,k+l}^+ - \Psi_{i,k+l}^-}{q_i-q_i^{-1}}
\]

on a weight space $\V_{\mu}$. Let $\lambda\in\hres$ and let us take
$\lambda_1 = \lambda$ and $\lambda_2 = -\lambda+\hbar\mu$. Again from Proposition
\ref{prop: eqg-prep} and definition \ref{eq:inv-Xpm} we get that the \lhs
is equal to
\begin{multline*}
\text{L.H.S = } \frac{1}{c_i^+c_j^-}\iint
e^{2\pi\iota (ku+lv)}\frac{\theta^+(u-v-a)}{\theta^+(u-v+a)}
G_i^+\ee{u}^{-1}G_j^-\ee{v}^{-1}.\\
\left[\eX_i^+(u,\lambda_1),\eX_j^-(v,\lambda_2)\right]\, du\, dv
\end{multline*}

Therefore, it is zero for $i\neq j$ by relation (\EQ 5). We now assume that $i=j$
and continue with the proof.
Let $c = 1/(c_i^+c_i^-\theta(d_i\hbar))$. Then

\begin{multline*}
\text{L.H.S = } c\int_{\cC}\int_{\cC}
e^{2\pi\iota (ku+lv)}\frac{\theta^+(u-v-d_i\hbar)}{\theta^+(u-v+d_i\hbar)}
G_i^+\ee{u}^{-1}G_i^-\ee{v}^{-1}.\\
\lp
\frac{\theta(u-v+\lambda_{1,i})}{\theta(u-v)\theta(\lambda_{1,i})}\Phi_i(v)
+\frac{\theta(u-v-\lambda_{2,i})}{\theta(u-v)\theta(\lambda_{2,i})}\Phi_i(u)
\rp
\, du\, dv
\end{multline*}

To continue with the calculation we use the same technique as in
\cite[Lemma 5.8]{sachin-valerio-2}. Namely let $\cC_<$ be a small deformation
of $\cC$ contained in $\cC$. Then we integrate with respect to $u$ first over
the contour $\cC_<$, to get:
\begin{multline*}
\text{L.H.S = } c\int_{\cC}\int_{\cC_<}
e^{2\pi\iota (ku+lv)}\frac{\theta^+(u-v-d_i\hbar)}{\theta^+(u-v+d_i\hbar)}
G_i^+\ee{u}^{-1}G_i^-\ee{v}^{-1}.\\
\lp
\frac{\theta(u-v-\lambda_{2,i})}{\theta(u-v)\theta(\lambda_{2,i})}\Phi_i(u)
\rp
\, du\, dv
\end{multline*}

Now we integrate this term with respect to the variable $v$. Note that
the integrand only has simple poles at $v=u$ within $\cC$. Using the
fact that $\theta$ is an odd function, we get
\begin{align*}
\text{L.H.S} &= c(2\pi\iota)\int_{\cC_<}
e^{2\pi\iota (k+l)u}\frac{\theta^+(-d_i\hbar)}{\theta^+(+d_i\hbar)}
G_i^+\ee{u}^{-1}G_i^-\ee{u}^{-1}\Phi_i(u)
\, du \\
&= c'\int_{\cC_<}e^{2\pi\iota (k+l)u}\Psi_i\ee{u}\, du \\
&= c'\oint_{\cC_<} z^{k+l-1}\Psi_i(z)\, dz \\
&= \frac{\Psi_{i,k+l}^+ - \Psi_{i,k+l}^-}{q_i-q_i^{-1}}
\end{align*}

where we have used the fact that $\Phi_i = G_i^+\Psi_iG_i^-$ in the second line,
the substitution $z=e^{2\pi\iota u}$ in the third line
and we have written $c' = 2\pi\iota c \theta^-(d_i\hbar)/\theta^+(d_i\hbar)$.
Note that by \eqref{eq:constant-fixing} we get
\[
c' = \frac{2\pi\iota}{\theta(d_i\hbar)} 
\frac{\theta^+(d_i\hbar)^2}{(2\pi\iota)^2 \theta^+(0)^2}
\frac{\theta^-(d_i\hbar)}{\theta^+(d_i\hbar)}
 = \frac{1}{q_i-q_i^{-1}}
\]
by \eqref{eq:thetaspm}.

\subsection{Proof of $\sfTh\circ\sfhT = \Id$}\label{ssec:inv-pf-c1}

Let $\eV\in\eL$. By definition \eqref{eq:ellX}, $\eX_i^{\pm}(u,\lambda)$
acts on $\sfTh\circ\sfhT(\eV)$ by the following

\begin{multline*}
c_i^{\pm}\oint_{\cC} \ttheta{u-v+\lambda_i}{u-v}{\lambda_i}
G_i^{\pm}\ee{v}\X_i^{\pm}\ee{v}\, dv \\
 = \oint_{\cC}\int_{\cC'} 
\ttheta{u-v+\lambda_i}{u-v}{\lambda_i}
G_i^{\pm}\ee{v}\frac{e^{2\pi\iota v}}{e^{2\pi\iota v}-e^{2\pi\iota v'}}
G_i^{\pm}\ee{v'}^{-1}\eX_i^{\pm}(v',\lambda)\, dv\, dv'
\end{multline*}

Assuming $\cC'$ is contained in $\cC$, we integrate with respect to $v$ first
to get:

\[
\frac{1}{2\pi\iota} \int_{\cC'}
\ttheta{u-v'+\lambda_i}{u-v'}{\lambda_i}
\eX_i^{\pm}(v',\lambda)\, dv'
\]

This is same as $\eX_i^{\pm}(u,\lambda)$, since both functions have same quasi--periodicity
with respect to shifts under $\tau$ and their difference has no poles by our choice of
the contour $\cC'$ in the definition \eqref{eq:inv-Xpm}. Therefore this difference must
be zero.

\subsection{Proof of $\sfhT\circ\sfTh = \Id$}\label{ssec:inv-pf-c2}

Let $\V\in\Rloop$. By definition \eqref{eq:inv-Xpm} the operators $\X_i^{\pm}(z)$
act on $\sfhT\circ\sfTh(\V)$ by the following
\begin{multline*}
\frac{1}{c_i^{\pm}}\int_{\cC}\frac{z}{z-e^{2\pi\iota v}} G_i^{\pm}\ee{v}^{-1}
\eX_i^{\pm}(v,\lambda)\, dv \\
 = \int_{\cC}\oint_{\cC'}\frac{z}{z-e^{2\pi\iota v}} G_i^{\pm}\ee{v}^{-1}
\ttheta{v-v'+\lambda_i}{v-v'}{\lambda_i} G_i^{\pm}\ee{v'}\X^{\pm}\ee{v'}
\, dv'\, dv
\end{multline*}

Again we integrate with respect to $v$ first, assuming $\cC'$ is contained in $\cC$,
to get

\[
\oint_{\cC'} \frac{z}{z-e^{2\pi\iota v'}} \X_i^{\pm}\ee{v'} (2\pi\iota)dv'
\]

Setting $w = e^{2\pi\iota v'}$, the last expression equals
\[
\oint_{\wt{\cC'}} \frac{z}{z-w} \X_i^{\pm}(w) \frac{dw}{w} = \X_i^{\pm}(z)
\]
where $\wt{\cC'} = \exp\lp 2\pi\iota\cC'\rp$.

\section{Classification of irreducibles II: Sufficient condition}
\label{sec: class2}

\subsection{Classification theorem}\label{ssec: class2-thm}

Consider $\eV\in\eL$ an irreducible object. Let $\mu\in\h^*$ be a weight of $\eV$
such that $\eV_{\mu+\alpha_i}=0$ for every $i\in\bfI$. Choose $\vac\in\eV_{\mu}$
a non--zero eigenvector for $\{\Phi_i(u)\}_{i\in\bfI}$. That is, there exists
$\ulA = (A_i(u)) \in \Mer(\mu)$ (see Section \ref{ssec: Phi-eigen} for definition
of $\Mer(\mu)$) such that $\Phi_i(u)\vac = A_i(u)\vac$.

\begin{thm}\label{thm: class-final}\hfill
\begin{enumerate}
\item For each $i\in\bfI$ there exists $N_i\in\N$, $b^{(i)}_1,\cdots,b^{(i)}_{N_i}
\in\C$, and a constant $C_i\in\nC$ such that 
\[
A_i(u) = C_i \prod_{k=1}^{N_i} \frac{\theta(u-b^{(i)}_k+d_i\hbar)}{\theta(u-b^{(i)}_k)}
\]
(in particular $\mu(\alpha_i^{\vee}) = N_i\in\N$ and hence $\mu$ is a dominant
integral weight). We say that the data of such $(\mu,\bfb)$ is the $l$--highest
weight of $\eV$, where 
$$\bfb=(\bfb_i=\{b^{(i)}_1,\cdots,b^{(i)}_{N_i}\})_{i\in\bfI}$$

\item Let $\eV$ and $\eW$ be two irreducible objects of $\eL$ with $l$--highest
weights $(\mu,\bfb)$ and $(\nu,\bfc)$. Then $\eV$ and $\eW$ are isomorphic
if, and only if $\mu=\nu$ and for every $i\in\bfI$ there exists $\sigma_i
\in\Sym_{N_i}$ such that $b^{(i)}_k\equiv c^{(i)}_{\sigma_i(k)}$
modulo $\Lambda_{\tau}$.\\

\item For every $(\mu,\bfb)$ such that $\mu\in\h^*$ is a dominant integral weight
and $\bfb=(\bfb_i)$ where each $\bfb_i$ is an unordered $\mu(\alpha_i^{\vee})$--tuple
of complex numbers, there exists an irreducible object $\eV(\mu,\bfb)$ of $\eL$
of $l$--highest weight $(\mu,\bfb)$.
\end{enumerate}
\end{thm}

\begin{pf}
(1) follows from Proposition \ref{prop: gauge1} and Theorem \ref{thm: class1}.\\

(2): Fix a subset $\slice\subset\C$ as in Section \ref{ssec: choice-slice}.
Using Propositions \ref{prop: gauge1}, \ref{prop:lambda-independence} 
and Theorem \ref{thm:inverse} we know that there exist (necessarily irreducible)
representations $\V, \W\in\Rloopnc$ such that $\eV$ and $\eW$ are isomorphic
to $\sfTh(\V)$ and $\sfTh(\W)$ respectively. If $l$--weights of $\eV$ 
and $\eW$ are the same as stated in (2), then $\V$ and $\W$ have the same
Drinfeld polynomials and the result follows from the classification 
theorem for irreducible representations of quantum loop algebras
(Theorem \ref{thm: qla-dp}).\\

(3): Now let $(\mu,\bfb)$ be as given above. Note that by (1)
the data of $l$--highest weight for an irreducible object of $\eL$ depends only
on the equivalence class of $b^{(i)}_k$ modulo $\Lambda_{\tau}$. Thus
we can assume without loss of generality that $b^{(i)}_k\in \slice$
for each $i\in\bfI$ and $1\leq k\leq N_i = \mu(\alpha_i^{\vee})$.
Let $\V$ be the irreducible representation of $\qloop$ of $l$--highest
weight given by $(\mu,\{\PP_i(w)\})$ where the Drinfeld polynomials $\PP_i$
are given by
\[
\PP_i(w) = \prod_{k=1}^{N_i} \lp w-e^{2\pi\iota b^{(i)}_k}\rp
\]

Let $\vac
\in\V_{\mu}$ be the (unique up to scalar) highest weight vector. 
Thus the action of the commuting elements $\{\Psi_i(z)\}$ on 
$\vac$ is given by
\[
\Psi_i(z)\vac = \prod_{k=1}^{N_i}\frac{q_iz-q_i^{-1}\beta^{(i)}_k}
{z-\beta^{(i)}_k}
\]
where $\beta^{(i)}_k = \exp\lp 2\pi\iota b^{(i)}_k\rp$.\\

Now we have
the following calculation for the action of $\Phi_i(u)$ on $\vac$, using
the definition \eqref{eq:functorPhi} and identity \eqref{eq: Jacobi-triple}.

\begin{align*}
\Phi_i(u)\vac &= \prod_{k=1}^{N_i} \lp 
\lp \prod_{n\geq 1} \frac{1-(q_i^2z/\beta^{(i)}_k)p^n}{1-(z/\beta^{(i)}_k)p^n}\rp
\frac{q_iz-q_i^{-1}\beta^{(i)}_k}{z-\beta^{(i)}_k}
\lp \prod_{n\geq 1} \frac{1-(\beta^{(i)}_k/q_i^2z)p^n}{1-(\beta^{(i)}_k/q_i^2z)}
\rp\rp\vac\\
&= \prod_{k=1}^{N_i} \frac{\theta(u-b^{(i)}_k+d_i\hbar)}{\theta(u-b^{(i)}_k)}\vac
\end{align*}

Thus it only remains to show that $\eV$ is irreducible. Assuming the contrary,
let $\eV_1\subset \eV$ be a non--zero proper subspace which is stable under
the action of $\{\Phi_i(u), \eX_i^{\pm}(u,\lambda)\}$. Then $\V_1 = \sfhT(\eV_1)$
is a non--zero proper subrepresentation of $\V$ contradicting its irreducibility.
\end{pf}

\subsection{Schur's Lemma}\label{ssec:schur}

\begin{prop}\label{prop:schur}
Let $\eV$ and $\eW$ be two irreducible objects of $\eL$. Then either
$\eV$ and $\eW$ are isomorphic, or $\Hom_{\eL}(\eV,\eW) = \{0\}$.
\end{prop}

\begin{pf}
By Corollary \ref{cor:secondgauge} we can assume that both $\eV$ and $\eW$
are in $\eLskeleton$. Let $\varphi(\lambda) : \eV \to \eW$ be a morphism.
Let $\Omega_{\eV} \in \eV_{\mu}$ be the (up to a scalar) highest--weight
vector of $\eV$. For every $i\in\bfI$, let $N_i = \mu(\alpha_i^{\vee})$ and let
$b^{(i)}_1,\cdots,b^{(i)}_{N_i} \in \slice$ be such that
\[
\Phi_i(u)_{\eV}\Omega_{\eV} = \prod_{k=1}^{N_i} \frac{\theta(u-b^{(i)}_k + d_i\hbar)}
{\theta(u-b^{(i)}_k)} \Omega
\]
Below we will just write $A_i(u)$ for the eigenvalue of $\Phi_i(u)$ on $\Omega_{\eV}$.\\

Let us define $\Omega(\lambda) = \varphi(\lambda)(\Omega_{\eV})$
in $\eW_{\mu}$. According to the definition of morphisms (see Section \ref{ssec: categoryL})
we obtain the following two identities for each $i\in \bfI$:
\begin{equation}\label{eq:schurpf1}
\Phi_i(u)_{\eW}\Omega(\lambda) = A_i(u) \Omega(\lambda+\hbar\alpha_i)
\end{equation}

\begin{equation}\label{eq:schurpf2}
\eX_i^+(u,\lambda)\Omega\lp\lambda-\frac{\hbar}{2}(\mu+\alpha_i)\rp = 0
\end{equation}

Recall that by Proposition \ref{prop:lambda-independence} we know the following
expression for $\eX_i^+(u,\lambda)$ acting on $\eW_{\mu}$:
\[
\eX_i^+(u,\lambda) = \sum_{\begin{subarray}{c} a\in \slice \\  n\in \N \end{subarray}}
\frac{\partial_u^n}{n!}\lp \ttheta{u-a+\lambda_i}{u-a}{\lambda_i}\rp \partialX^+_{i;a,n}
\]
where the sum is finite and $\partialX^+_{i;a,n}$ do no depend on the dynamical variable.
Now we can multiply equation \eqref{eq:schurpf2} by $(u-a)^n$ and integrate over
a small contour around $a$ to get that $\partialX^+_{i;a,n}$ annihilate $\Omega(\lambda)$
for all $\lambda$. Hence, for $\lambda_0\in\hres$ such that $\varphi(\lambda_0)$ is 
defined on $\eV_{\mu}$, $\Omega(\lambda_0)\in \eW_{\mu}$ is annihilated by all
the raising operators.\\

By irreducibility of $\eW$, if $\mu$ is not the highest--weight of $\eW$, we get
that $\varphi(\lambda)\Omega_{\eV} = 0$. By Theorem \ref{thm: eqg-pbw}, this means that
$\varphi(\lambda) \equiv 0$.\\

Now assume that the highest--weight of $\eW$ is also $\mu$ and let $\Omega_{\eW} \in \eW_{\mu}$
be the unique (up to scalar) highest--weight vector of $\eW$. Thus we have a scalar function
$\phi(\lambda)$ such that $\Omega(\lambda) = \phi(\lambda)\Omega_{\eW}$. 
Let us write $B_i(u)$ for the eigenvalue of $\Phi_i(u)$ acting on $\Omega_{\eW}$. Then
by \eqref{eq:schurpf1} we get

\[
\phi(\lambda)B_i(u) = \phi(\lambda+\hbar\alpha_i)A_i(u)
\]

Thus $B_i(u) = C_iA_i(u)$ for some $C_i\in \nC$. We know from Proposition \ref{prop:lambda-independence}
that $C_i = 1$. This proves that $\eV$ and $\eW$ have the same $l$--highest weight, and hence
so must $\V = \sfhT(\eV)$ and $\W = \sfhT(\eW)$. Let $\psi : \V \to \W$ be the isomorphism
in $\Rloopslice$. Then $\psi$ is also an isomorphism between $\eV = \sfTh(\V)$ and
$\eW = \sfTh(\W)$ and we are done.

\end{pf}

\subsection{}\label{ssec:generalinjective}

\begin{cor}\label{cor:generalinjective}
Let $\eV$ and $\eW$ be two objects of $\eL$. Assume $\eW$ is irreducible
and $\eV \neq \{0\}$.
Let there be an injective morphism $\varphi(\lambda) : \eV \to \eW$.
Then $\varphi(\lambda)$ is an isomorphism.
\end{cor}

\begin{pf}
Again using Corollary \ref{cor:secondgauge} we can assume that $\eV$
and $\eW$ are in $\eLskeleton$. Let $\eV_1\subset \eV$ be a minimal (\wrt inclusion)
non--zero subspace which is stable under the action of $\{\Phi_i(u),\eX_i^{\pm}(u,\lambda)\}$.
Using Proposition \ref{prop:schur} we know that either $\varphi(\lambda)$ restricted
to $\eV_1$ is zero; or $\eV_1$ and $\eW$ are isomorphic, via an isomorphism which
is independent of $\lambda$. In the former case we contradict the injectivity
of $\varphi(\lambda)$. In the latter case, let us use the isomorphism between $\eV_1$
and $\eW$ to identify the two. \\  

This brings us to the situation where we
have $j: \eW \subset \eV$ as a subspace (stable under the action of $\{\Phi_i(u),
\eX_i^{\pm}(u,\lambda)\}$) and $\varphi(\lambda) : \eV \to \eW$ such that
$\varphi(\lambda)$ restricted to $\eW$ is identity. In other words, $\psi(\lambda) = j \circ \varphi(\lambda)$
is injective (being a composition of two injective morphisms) and a projection
(that is, $\psi \circ \psi = \psi$). Now it is a general fact (also not so hard to prove)
that an injective projection has to be identity, which finishes the proof
of the corollary.

\end{pf}

\appendix

\section{Serre relations}\label{sec: serre-relns}

We formulate below an analogue of the Serre relations for elliptic quantum groups,
conjecture that they hold on any \fd representation, and prove that this is the case
if $\g$ is simply--laced. These relations are not used anywhere in this paper, and
are included solely for completeness.

\subsection{}

Fix $i\neq j\in \bfI$, and let $m = 1-a_{ij}$. For an integer $n\in \Z$, we use
the following notation
\begin{gather*}
[n]_{\hbar,\tau} = \frac{\theta(n\hbar)}{\theta(\hbar)}\\[.5ex]
[n]_{\hbar,\tau}! = [n]_{\hbar,\tau}[n-1]_{\hbar,\tau}\cdots [1]_{\hbar,\tau} \qquad
\qbin{n}{k}{\hbar,\tau} = \frac{[n]_{\hbar,\tau}!}{[k]_{\hbar,\tau}![n-k]_{\hbar,\tau}!}
\end{gather*}

In this section use the following abbreviation (where recall that $\lambda_i = (\lambda,\alpha_i)$)
\[
F_i(u,\lambda) = \theta(\lambda_i)\eX_i^-(u,\lambda) \qquad
F_j(v,\lambda) = \theta(\lambda_j)\eX_j^-(v,\lambda)
\]

Define
\begin{equation}\label{eq: serre1}
\mathcal{S}(u,v) = \sum_{k=0}^m (-1)^{m-k}\qbin{m}{k}{d_i\hbar,\tau} 
\frac{\theta\lp u-v+\lambda_j+\hbar\lp k-\frac{1}{2}\rp d_ia_{ij} + \hbar d_j\rp}
{\theta(u-v)\theta\lp \lambda_j+\hbar\lp k-\frac{1}{2}\rp d_ia_{ij} + \hbar d_j\rp}
T_k(u,v)
\end{equation}

where the terms $T_k(u,v)$ are given by

\begin{align*}
T_k(u,v) &= 
F_i(u,\lambda) \cdots F_i(u,\lambda+(k-1)\hbar\alpha_i) \cdot\\
& \qquad \cdot F_j\lp v+\lambda_i-\lambda_j,\lambda+(k-1)\hbar\alpha_i + \frac{\hbar}{2}(\alpha_i+\alpha_j)\rp \cdot\\
& \qquad \cdot F_i(u,\lambda+k\hbar\alpha_i+\hbar\alpha_j)\cdots F_i(u,\lambda+(m-1)\hbar\alpha_i+\hbar\alpha_j)
\end{align*}

\begin{defn}
The Serre relations are the relations
\begin{equation}\label{eq:ell-serre}
\tag{\EQ 6}\mathcal{S}(u,v) \text{ is independent of } v
\end{equation}
\end{defn}

\begin{conj}\label{co:Serre}
The Serre relations \eqref{eq:ell-serre} hold on any $\V\in\ReL$.
\end{conj}

\subsection{}

We now explain the analogy of relation \eqref{eq:ell-serre}
with the Serre relations for quantum loop algebras.
Recall that the Serre relations (QL6) for 
the quantum loop algebra take the following form
\begin{multline*}
\sum_{\sigma\in\Sym_m} \sum_{k=0}^m (-1)^{m-k}
\left[\begin{array}{c} m \\ k \end{array}\right]_{q_i}
\X^-_i(z_{\sigma(1)})\cdots \X^-_i(z_{\sigma(k)})
\X^-_j(w)\\
\cdot \X^-_i(z_{\sigma(k+1)})\cdots \X^-_i(z_{\sigma(m)}) = 0
\end{multline*}
It is easy to show, using the relation (\QL 3) of Definition \ref{defn: rloop} above,
that this relation is equivalent to its special case
\[
\sum_{k=0}^m (-1)^{m-k}
\left[\begin{array}{c} m \\ k \end{array}\right]_{q_i}
\X^-_i(z)^{k}\X^-_j(w)
\X^-_i(z)^{m-k} = 0
\]

We notice that this relation can be written in an apparent weaker form

\[
\sum_{k=0}^m (-1)^{m-k}
\left[\begin{array}{c} m \\ k \end{array}\right]_{q_i}
\X^-_i(z)^{k}\X^-_j(w)
\X^-_i(z)^{m-k} \text{ is independent of } w .
\]
This is because, as $w\to 0$, $\X^-_j(w)\to 0$, hence the two assertions above
are equivalent. The latter one, however, appears to have a generalization
to the elliptic world.

\begin{rem}\label{rem: howto-serre}\hfill
\begin{enumerate}
\item The function $\mathcal{S}(u,v)$ is doubly--periodic in $v$. Thus showing that 
it is independent of $v$ is equivalent to proving that it has no poles in $v$.

\item Clearing the denominator $\theta(u-v)$ from the definition of $\mathcal{S}(u,v)$ given in
\eqref{eq: serre1} and setting $v=u$ yields a more familiar version of the Serre relations

\[
\sum_{k=0}^m (-1)^{m-k}\qbin{m}{k}{d_i\hbar,\tau} 
T_k(u,u)
=0
\]
\end{enumerate}

\end{rem}

\subsection{}

The remainder of this appendix is devoted to proving the following

\begin{thm}
Conjecture \ref{co:Serre} holds if $\g$ is simply--laced.
\end{thm}
\begin{pf}
For the $\mathsf{A}_1\times\mathsf{A}_1$ root subsystem, 
the proof follows from the relation (\EQ 4) (see Section \ref{ssec: categoryL}),
as we demonstrate in \ref{asec:a1xa1} below.

By Corollary \ref{cor:secondgauge} and Theorem \ref{thm:inverse}, we may
assume that $\eV = \Th(\V)$, where $\V$ is a (non--congruent) \fd representation
of the quantum loop algebra $\qloop$. We prove the theorem
for $\mathsf{A}_2$ subsystem, in Section \ref{asec:a2} below. This is achieved
by adapting the analoguous proof for Yangians, as given in \cite[\S 7.11]{sachin-valerio-1}.
\end{pf}

\subsection{}\label{asec:a1xa1}

Assume that $a_{ij}=0$. Then, the following expression is independent of $v$
\begin{multline*}
\frac{\theta(u-v+\lambda_j+\hbar d_j)}{\theta(u-v)\theta(\lambda_j+\hbar d_j)}
\lp F_i(u,\lambda)F_j\lp v+\lambda_i-\lambda_j,\lambda+\frac{\hbar}{2}(\alpha_i+\alpha_j)\rp
\right.
\\ \left. - F_j\lp v+\lambda_i-\lambda_j,\lambda+\frac{\hbar}{2}(\alpha_j-\alpha_i)\rp
F_i(u,\lambda+\hbar\alpha_j)\rp
\end{multline*}

Thus the function above is equal to its evaluation at $v=u+\lambda_j+\hbar d_j$, where it is zero.
This implies, after replacing $\lambda$ by $\lambda - \hbar\alpha_j/2$, and renaming $v$

\[
\eX_i^-\lp u,\lambda -\frac{\hbar}{2}\alpha_j\rp \eX_j^- \lp v,\lambda+\frac{\hbar}{2}\alpha_i\rp = 
\eX_j^-\lp v,\lambda -\frac{\hbar}{2}\alpha_i\rp \eX_j^- \lp u,\lambda+\frac{\hbar}{2}\alpha_j\rp 
\]

This is precisely the relation (\EQ 4) for the case when $a_{ij}=a_{ji}=0$
(see Section \ref{ssec: categoryL}).

\subsection{$\mathsf{A}_2$ case.}\label{asec:a2}

Let us assume that $a_{ij}=a_{ji}=-1$ and $d_i=d_j=1$. 
Now we present a proof of the special case of the Serre relations, as stated in
Remark \ref{rem: howto-serre} (2) above. 
As mentioned in Remark \ref{rem: howto-serre} (1), the general case amounts to proving
that, for every $b\in \C$ and $n\in \N$, the following integral over a small circle
around $b$ is zero
\[
\oint (v-b)^n \mathcal{S}(u,v) = 0
\]

It is straightforward (by flipping the order of integration) to adapt the proof
of $\mathcal{S}(u)=0$ given below in order to show the equation given above.
This is the reason why we only focus on the special case.\\

We need to prove that the function $\mathcal{S}(u)=0$
where
$\mathcal{S}(u) = T_1 - \frac{\theta(2\hbar)}{\theta(\hbar)} T_2 + T_3$. The three terms
involved here are given by

\begin{align*}
T_1 &= 
F_i(u,\lambda)F_i\lp u,\lambda+\hbar\alpha_i\rp F_j\lp \wt{u},\lambda + \frac{3}{2}\hbar \alpha_i + \frac{\hbar}{2} \alpha_j\rp\ ,\\
T_2 &= 
F_i(u,\lambda) F_j\lp \wt{u},\lambda + \frac{\hbar}{2}(\alpha_i + \alpha_j)\rp F_i\lp u,\lambda+\hbar(\alpha_i+\alpha_j)\rp\ ,\\
T_3 &= 
F_j\lp \wt{u},\lambda - \frac{\hbar}{2} \alpha_i + \frac{\hbar}{2} \alpha_j\rp
F_i(u,\lambda+\hbar\alpha_j)F_i\lp u,\lambda+\hbar(\alpha_i+\alpha_j)\rp \ .
\end{align*}

and $\wt{u} = u+\lambda_i-\lambda_j$.
Using Theorem \ref{thm:inverse} and Corollary \ref{cor:secondgauge} above, it suffices
to work with $\eV = \Th(\V)$ for some
representation $\V$ of the quantum loop algebra $\qloop$.\\

\subsection{}

Recall the definition of the functor $\Th$ from Section \ref{ssec: functor-defn}. The terms
appearing in $\mathcal{S}(u)$ above can be written as follows, 
using the commutation relations stated in Proposition \ref{prop: prep}

\begin{align*}
T_1 &= \oint f_1(u; u_1,u_2;\lambda) G_i^-(u_1)G_i^-(u_2)G_j^-(v) \X_i^-(u_1)\X_i^-(u_2)\X_j^-(v)\, du_1\, du_2\, dv \\
T_2 &= \oint f_2(u; u_1,u_2;\lambda) G_i^-(u_1)G_i^-(u_2)G_j^-(v) \X_i^-(u_1)\X_j^-(v)\X_i^-(u_2)\, du_1\, du_2\, dv \\
T_3 &= \oint f_3(u; u_1,u_2;\lambda) G_i^-(u_1)G_i^-(u_2)G_j^-(v) \X_j^-(v)\X_i^-(u_1)\X_i^-(u_2)\, du_1\, du_2\, dv
\end{align*}

The functions appearing within the integral written above are

\begin{align*}
f_1(u; u_1, u_2; \lambda) &=
\thetaratio{u-u_1+\lambda_i}{u-u_1}\thetaratio{u-u_2+\lambda_i+2\hbar}{u-u_2}
\thetaratio{\wt{u}-v+\lambda_j-\frac{\hbar}{2}}{\wt{u}-v} \\
& \phantom{=} \cdot\ 
\thetaratiop{u_1-u_2-\hbar}{u_1-u_2+\hbar}
\thetaratiop{u_1-v+\frac{\hbar}{2}}{u_1-v-\frac{\hbar}{2}}
\thetaratiop{u_2-v+\frac{\hbar}{2}}{u_2-v-\frac{\hbar}{2}}\ \ , \\
f_2(u; u_1, u_2; \lambda) &=
\thetaratio{u-u_1+\lambda_i}{u-u_1}\thetaratio{u-u_2+\lambda_i+\hbar}{u-u_2}
\thetaratio{\wt{u}-v+\lambda_j+\frac{\hbar}{2}}{\wt{u}-v} \\
& \phantom{=} \cdot\ 
\thetaratiop{u_1-u_2-\hbar}{u_1-u_2+\hbar}
\thetaratiop{u_1-v+\frac{\hbar}{2}}{u_1-v-\frac{\hbar}{2}}
\thetaration{u_2-v-\frac{\hbar}{2}}{u_2-v+\frac{\hbar}{2}}\ \ ,\\
f_3(u; u_1, u_2; \lambda) &=
\thetaratio{u-u_1+\lambda_i-\hbar}{u-u_1}\thetaratio{u-u_2+\lambda_i+\hbar}{u-u_2}
\thetaratio{\wt{u}-v+\lambda_j+\frac{3\hbar}{2}}{\wt{u}-v} \\
& \phantom{=} \cdot\ 
\thetaratiop{u_1-u_2-\hbar}{u_1-u_2+\hbar}
\thetaration{u_1-v-\frac{\hbar}{2}}{u_1-v+\frac{\hbar}{2}}
\thetaration{u_2-v-\frac{\hbar}{2}}{u_2-v+\frac{\hbar}{2}}\ .
\end{align*}

From this point on, the proof proceeds exactly as in \cite[\S 7.11]{sachin-valerio-1}.
We briefly recall the idea of the proof

The analogue of Lemma 7.8 of \cite{sachin-valerio-1}, stated below (section \ref{ssec: serre-lemma}),
implies that we can replace $f_1$ (and $f_3$ resp.) by a function symmetric in $u_1, u_2$, denoted by
$f_1'$ (resp. $f_3'$), provided $f_1-f_1' = (z_1-q_i^{-2}z_2)f_1''$, where $f_1''$ is again symmetric
in $u_1,u_2$. Here $z_\ell = e^{2\pi\iota u_\ell}$ for $\ell = 1,2$. Similarly for $f_3$.
A lengthy, but straightforward, calculation using the identities of the theta function
given in Section \ref{ssec: ell-funs} gives the following result

\vspace*{-0.1in}

\begin{align*}
f_1'(u;u_1,u_2;\lambda) &= \frac{1}{[2]_{q}}
\thetaratio{2\hbar}{\hbar}
\frac{\theta(u_1-u_2)}{z_1-z_2}
\frac{q^{-1}z_1-qz_2}{\theta(u_1-u_2-\hbar)}
\thetaratiop{u_1-u_2-\hbar}{u_1-u_2+\hbar} \\
&\phantom{=}\ \cdot \thetaratio{u-u_1+\lambda_i+\hbar}{u-u_1}
\thetaratio{u-u_2+\lambda_i+\hbar}{u-u_2} \\
& \phantom{=}\ \cdot
\thetaratiop{u_1-v+\frac{\hbar}{2}}{u_1-v-\frac{\hbar}{2}}
\thetaratiop{u_2-v+\frac{\hbar}{2}}{u_2-v-\frac{\hbar}{2}}
\thetaratio{\wt{u}-v+\lambda_j-\frac{\hbar}{2}}{\wt{u}-v}
\end{align*}

\begin{align*}
f_3'(u;u_1,u_2;\lambda) &= \frac{1}{[2]_{q}}
\thetaratio{2\hbar}{\hbar}
\frac{\theta(u_1-u_2)}{z_1-z_2}
\frac{q^{-1}z_1-qz_2}{\theta(u_1-u_2-\hbar)}
\thetaratiop{u_1-u_2-\hbar}{u_1-u_2+\hbar} \\
&\phantom{=}\ \cdot \thetaratio{u-u_1+\lambda_i}{u-u_1}
\thetaratio{u-u_2+\lambda_i}{u-u_2} \\
& \phantom{=}\ \cdot
\thetaration{u_1-v-\frac{\hbar}{2}}{u_1-v+\frac{\hbar}{2}}
\thetaration{u_2-v-\frac{\hbar}{2}}{u_2-v+\frac{\hbar}{2}}
\thetaratio{\wt{u}-v+\lambda_j+\frac{3\hbar}{2}}{\wt{u}-v}
\end{align*}

Now $[2]_{\hbar,\tau}$ is common in all three terms of $\mathcal{S}(u)$. Scaling by it
reduces the proof to showing that
\[
T_1' - [2]_{q} T_2' + T_3' = 0
\]
where $T_2' = T_2$ and $T_1',T_3'$ are defined exactly as before, but with
kernels $[2]_{q}/[2]_{\hbar,\tau} f_\ell'$ ($\ell = 1,3$).\\

Again using Lemma \ref{lem: serre-lemma} (3) below, the identity to prove simplifies to:
$T_2''- T_3'' = 0$, where the kernels defining $T_2''$ and $T_3''$ are

\[
f_2'' = f_2 - \frac{[2]_{q}}{[2]_{\hbar,\tau}} f_1' \qquad
f_3'' = \frac{1}{[2]_{\hbar,\tau}}(f_3' - f_1')
\]

The remainder of the proof is the exact same verification as in Steps 2 and 3 of
\cite[\S 7.11]{sachin-valerio-1}. Namely, the proof reduces to proving that the following
function is symmetric in $u_1,u_2$

\[
\frac{1}{z_1-q^{-2}z_2}\lp \frac{qz_1-w}{z_1-qw} f_2'' - f_3''\rp
\]
where, as usual, we have $z_{\ell} = e^{2\pi\iota u_{\ell}}$ and $w = e^{2\pi\iota v}$.
Again, we leave this direct verification to an interested reader.


\subsection{}\label{ssec: serre-lemma}

We retain the notations of the previous section. Namely, we are working with a representation
$\V\in \Rloop$, and under the assumption that $i,j\in\bfI$ are such that $a_{ij} = a_{ji} = -1$,
and $d_i = d_j=1$.

Assume we are given meromorphic functions of complex variables $A(z,w)$ and $B(z_1,z_2,w)$
taking values in the commutative subalgebra of $\End(\V)$ generated by the commuting generators
$\{\Psi_k(z)\}_{k\in\bfI}$, such that $B(z_1,z_2)$ is symmetric in $z_1,z_2$.

\begin{lem}\label{lem: serre-lemma}
The following relations hold in $\End(\V)$. The integrals below are over arbitrary
contours such that the given functions $A(z,w)$ and $B(z_1,z_2,w)$ are holomorphic
in the interior of, and on the contour.

\begin{enumerate}
\item
\[
\oint (z_1-q^{-2}z_2) B(z_1,z_2) \X_i^-(z_1)\X_i^-(z_2)\, dz_1\, dz_2 = 0
\]

\item

\[
\oint (z-qw) A(z,w) \X_i^-(z)\X_j^-(w)\, dz\, dw = 
\oint (qz-w) A(z,w) \X_j^-(w)\X_i^-(z)\, dz\, dw
\]

\item

\begin{multline*}
\oint B(z_1,z_2,w) \lp \X_i^-(z_1)\X_i^-(z_2)\X_j^-(w) - 
[2]_q \X_i^-(z_1)\X_j^-(w)\X_i^-(z_2)\right. \\
\left.+ \X_j^-(w)\X_i^-(z_1)\X_i^-(z_2)\rp\, dz_1\, dz_2\, dw = 0
\end{multline*}

\end{enumerate}

\end{lem}

\bibliographystyle{amsplain}
\bibliography{LambOfGod}

\providecommand{\bysame}{\leavevmode\hbox to3em{\hrulefill}\thinspace}
\providecommand{\MR}{\relax\ifhmode\unskip\space\fi MR }
\providecommand{\MRhref}[2]{%
  \href{http://www.ams.org/mathscinet-getitem?mr=#1}{#2}
}
\providecommand{\href}[2]{#2}
\begin{thebibliography}{10}

\bibitem{aganagic-okounkov}
M.~Aganagic and A.~Okounkov, \emph{Elliptic stable envelope}, {\sf
  arXiv:160400423}.

\bibitem{bbb}
O.~Babelon, D.~Bernard, and E.~Billey, \emph{{A quasi--Hopf algebra
  interpretation of quantum $3j$ and $6j$--symbols and difference equations}},
  Phys. Lett. B \textbf{375} (1996), 89--97.

\bibitem{baxter-8v}
R.~J. Baxter, \emph{Partition function of the eight--vertex lattice model},
  Ann. Physics \textbf{70} (1972), 193--228.

\bibitem{baxter}
\bysame, \emph{Exactly solved models in statistical mechanics}, Academic Press,
  Inc. [Harcourt Brace Jovanovich, Publishers], London, 1982.

\bibitem{beck-kac}
J.~Beck and V.~G. Kac, \emph{{Finite--dimensional representations of quantum
  affine algebras at roots of unity}}, J. Amer. Math. Soc. \textbf{9} (1996),
  no.~2, 391--423.

\bibitem{belavin}
A.~A. Belavin, \emph{{Dynamical symmetry of integrable quantum systems}},
  Nuclear Physics B \textbf{180} (1981), 189--200.

\bibitem{belavin-drinfeld-book}
A.~A. Belavin and V.~G. Drinfel$\prime$d, \emph{Triangle equations and simple
  {L}ie algebras}, Mathematical physics reviews, {V}ol. 4, Soviet Sci. Rev.
  Sect. C Math. Phys. Rev., vol.~4, Harwood Academic Publ., Chur, 1984,
  Translated from the Russian, pp.~93--165.

\bibitem{chari-pressley-yangian}
V.~Chari and A.~Pressley, \emph{{Yangians and $R$--matrices}}, Enseign. Math.
  \textbf{36} (1990), 267--302.

\bibitem{chari-pressley-qaffine}
\bysame, \emph{Quantum affine algebras}, Comm. Math. Phys. \textbf{142} (1991),
  261--283.

\bibitem{chari-pressley}
\bysame, \emph{A guide to quantum groups}, Cambridge University Press, 1994.

\bibitem{chari-pressley-qaffine-rep}
\bysame, \emph{{Quantum affine algebras and their representations}},
  {Representations of Groups}, CMS Conf. Proc., vol.~16, 1995, pp.~59--78.

\bibitem{cherednik-sklyanin}
I.~Cherednik, \emph{{Some finite--dimensional representations of generalized
  Sklyanin algebras}}, Funct. Anal. Appl. \textbf{19} (1985), 77--79.

\bibitem{drinfeld-yangian-qaffine}
V.~G. Drinfeld, \emph{{A new realization of Yangians and quantum affine
  algebras}}, Soviet Math. Dokl. \textbf{36} (1988), no.~2, 212--216.

\bibitem{felder-enriquez}
B.~Enriquez and G.~Felder, \emph{{Elliptic quantum group
  $E_{\tau,\eta}(\mathfrak{sl}_2)$ and quasi--Hopf algebras}}, Comm. Math.
  Phys. \textbf{195} (1998), 651--689.

\bibitem{etingof-latour}
P.~I. Etingof and F.~Latour, \emph{The dynamical {Y}ang--{B}axter equation,
  representation theory and quantum integrable systems}, Lecture series in
  mathematics and its applications - 29, Oxford University Press, 2005.

\bibitem{etingofmoura-KL}
P.~I. Etingof and A.~A. Moura, \emph{On the quantum {K}azhdan-{L}usztig
  functor}, Math. Res. Lett. \textbf{9} (2002), 449--463.

\bibitem{etingof-schiffmann-link}
P.~I. Etingof and O.~Schiffmann, \emph{{A link between two elliptic quantum
  groups}}, Asian J. Math. \textbf{2} (1998), no.~2, 345--354.

\bibitem{etingofschiffmann-elliptic}
\bysame, \emph{On highest weight modules over elliptic quantum groups}, Lett.
  Math. Phys. \textbf{47} (1999), 179--188.

\bibitem{etingof-varchenko-classical}
P.~I. Etingof and A.~Varchenko, \emph{Geometry and classification of solutions
  of the classical dynamical {Y}ang--{B}axter equation}, Comm. Math. Phys.
  \textbf{192} (1998), 77--120.

\bibitem{etingof-varchenko-quantum}
\bysame, \emph{Solutions of the quantum dynamical {Y}ang--{B}axter equation and
  dynamical quantum group}, Comm. Math. Phys. \textbf{196} (1998), 591--640.

\bibitem{etingof-varchenko-exchange}
\bysame, \emph{Exchange dynamical quantum groups}, Comm. Math. Phys.
  \textbf{205} (1999), 19--52.

\bibitem{frt-quantization}
L.~D. Faddeev, N.~Reshetikhin, and L.~A. Takhtajan, \emph{Quantization of {L}ie
  groups and {L}ie algebras}, Algebra and Analysis \textbf{11} (1989),
  118--206.

\bibitem{konno-farghly-oshima}
R.~M. Farghly, H.~Konno, and K.~Oshima, \emph{Elliptic algebra
  {$U_{q,p}(\widehat{\mathfrak{g}})$} and quantum {$Z$}--algebras}, Algebr.
  Represent. Theory \textbf{18} (2015), no.~1, 103--135.

\bibitem{felder-icm}
G.~Felder, \emph{Conformal field theory and integrable systems associated to
  elliptic curves}, Proceedings of the {I}nternational {C}ongress of
  {M}athematicians, {V}ol.\ 1, 2 ({Z}\"urich, 1994), Birkh\"auser, Basel, 1995,
  pp.~1247--1255.

\bibitem{fv-reps}
G.~Felder and A.~Varchenko, \emph{{On representations of the elliptic quantum
  group $E_{\tau,\eta}(\mathfrak{sl}_2)$}}, Comm. Math. Phys. \textbf{181}
  (1996), 741--761.

\bibitem{felder-wieczerkowski}
G.~Felder and C.~Wieczerkowski, \emph{Conformal blocks on elliptic curves and
  the {K}nizhnik--{Z}amolodchikov--{B}ernard equation}, Comm. Math. Phys.
  \textbf{176} (1996), 133--161.

\bibitem{feigin-odesskii1}
B.~L. Fe\u{\i}gin and A.~V. Odesski\u{\i}, \emph{Sklyanin's elliptic algebras},
  Funktsional. Anal. i Prilozhen. \textbf{23} (1989), no.~3, 45--54, 96.

\bibitem{feigin-odesskii2}
\bysame, \emph{Vector bundles on an elliptic curve and {S}klyanin algebras},
  Topics in quantum groups and finite-type invariants, Amer. Math. Soc. Transl.
  Ser. 2, vol. 185, Amer. Math. Soc., Providence, RI, 1998, pp.~65--84.

\bibitem{fronsdal1}
C.~Fr$\o$nsdal, \emph{Generalization and exact deformations of quantum groups},
  Publ. RIMS, Kyoto \textbf{33} (1997), 91--149.

\bibitem{fronsdal2}
\bysame, \emph{{Quasi--Hopf deformation of quantum groups}}, Lett. Math. Phys.
  \textbf{40} (1997), 117--134.

\bibitem{sachin-rtt}
S.~Gautam, \emph{{Drinfeld-type presentation of the elliptic quantum group
  $\mathcal{E}_{\hbar,\tau}(\mathfrak{sl}_n)$}}, {in preparation}.

\bibitem{sachin-valerio-1}
S.~Gautam and V.~Toledano~Laredo, \emph{Yangians and quantum loop algebras},
  Selecta Math. (N.S.) \textbf{19} (2013), no.~2, 271--336.

\bibitem{sachin-valerio-2}
\bysame, \emph{Yangians, quantum loop algebras, and abelian difference
  equations}, J. Amer. Math. Soc. \textbf{29} (2016), no.~3, 775--824.

\bibitem{hernandez-affinizations}
D.~Hernandez, \emph{{Representations of quantum affinizations and fusion
  product}}, Transform. Groups \textbf{10} (2005), no.~2, 163--200.

\bibitem{hernandez-q-toroidal}
\bysame, \emph{Quantum toroidal algebras and their representations}, Selecta
  Math. (N.S.) \textbf{14} (2009), no.~3-4, 701--725.

\bibitem{jkos}
M.~Jimbo, H.~Konno, S.~Odake, and J.~Shiraishi, \emph{Quasi-{H}opf twistors for
  elliptic quantum groups}, Transform. Groups \textbf{4} (1999), no.~4,
  303--327.

\bibitem{kac}
V.~G. Kac, \emph{Infinite dimensional {L}ie algebras}, Cambridge University
  Press, 1990.

\bibitem{knight}
H.~Knight, \emph{{Spectra of tensor products of finite dimensional
  representations of Yangians}}, J. Algebra \textbf{174} (1995), 187--196.

\bibitem{kojima-konno}
T.~Kojima and H.~Konno, \emph{{Elliptic algebra $U_{q,p}(\widehat{sl_N})$ and
  the Drinfeld realization of the elliptic quantum group
  $B_{q,\lambda}(\widehat{\sl_N})$}}, Comm. Math. Phys. \textbf{239} (2003),
  405--447.

\bibitem{konno-first}
H.~Konno, \emph{An elliptic algebra {$U_{q,p}(\widehat{\rm sl}_2)$} and the
  fusion {RSOS} model}, Comm. Math. Phys. \textbf{195} (1998), no.~2, 373--403.

\bibitem{konno2}
\bysame, \emph{{Elliptic quantum groups $U_{q,p}(\gl_N)$ and
  $E_{q,p}(\gl_N)$}}, Representation Theory, Special Functions and Painlev\'{e}
  equations -- RIMS 2015, Advanced Series in Pure Mathematics, Mathematical
  Society of Japan, Tokyo, Japan, 2018, pp.~347--417.

\bibitem{maulik-okounkov-qgqc}
D.~Maulik and A.~Okounkov, \emph{{Quantum groups and quantum cohomology}},
  (2012), {\sf arXiv:1211.1287}.

\bibitem{molev-yangian}
A.~Molev, \emph{Yangians and classical {L}ie algebras}, Mathematical Surveys
  and Monographs, vol. 143, A.M.S., 2007.

\bibitem{nakajima-qaffine}
H.~Nakajima, \emph{Quiver varieties and finite dimensional representations of
  quantum affine algebras}, J. Amer. Math. Soc. \textbf{14} (2001), 145--238.

\bibitem{sauloy}
J.~Sauloy, \emph{{Syst\`{e}mes aux $q$--diff\'{e}rences singuliers
  r\'{e}guliers: classification, matrices de connexion et monodromie}}, Ann.
  Inst. Fourier (Grenoble) \textbf{50} (2000), no.~4, 1021--1071.

\bibitem{sklyanin}
E.~K. Sklyanin, \emph{{Some algebraic structures connected with the
  Yang--Baxter equation}}, Funct. Anal. Appl. \textbf{16} (1982), 263--270.

\bibitem{sklyanin2}
\bysame, \emph{{Some algebraic structures connected with the Yang--Baxter
  equation. Representations of a quantum algebra.}}, Funct. Anal. Appl.
  \textbf{17} (1983), 273--284.

\bibitem{varagnolo-yangian}
M.~Varagnolo, \emph{Quiver varieties and {Y}angians}, Lett. Math. Phys.
  \textbf{53} (2000), no.~4, 273--283.

\bibitem{whittaker-watson}
E.~T. Whittaker and G.~N. Watson, \emph{A course of modern analysis}, 4th ed.,
  Cambridge University Press, 1927.

\bibitem{yang-zhao-elliptic}
Y.~Yang and G.~Zhao, \emph{{Quiver varieties and elliptic quantum groups}},
  {\sf arXiv:1708.01418}.

\bibitem{young}
C.~Young, \emph{{Quantum loop algebras and $\ell$--root operators}}, Transform.
  Groups \textbf{20} (2015), 1195--1226.

\end{thebibliography}
\end{document}